\documentclass[12pt]{article}
\usepackage{amsfonts, amsthm, amsmath, amssymb}
\usepackage{amsxtra, amsopn}
\usepackage{array} 

\theoremstyle{plain}
\newtheorem{theorem}{Theorem}[section]
\newtheorem{corollary}[theorem]{Corollary}
\newtheorem{lemma}[theorem]{Lemma}
\newtheorem{proposition}[theorem]{Proposition}

\newtheorem{conjecture}{Conjecture}[section]
\newtheorem{assumption}{Assumption}[section]

\theoremstyle{definition}
\newtheorem{definition}{Definition}[section]
\newtheorem{remark}{Remark}[section]
\newtheorem{example}{Example}[section]

\newcommand{\limfunc}[1]{\mathop{\rm #1}}
\def\wt{\widetilde}
\def\dint{\displaystyle \int }
\DeclareMathOperator{\sgn}{sgn}
\DeclareMathOperator{\im}{Im}
\DeclareMathOperator{\re}{Re}
\DeclareMathOperator{\supp}{supp}
\DeclareMathOperator{\dist}{dist}
\DeclareMathOperator{\sn}{sn}
\DeclareMathOperator{\dom}{dom}
\DeclareMathOperator{\ran}{ran}
\DeclareMathOperator{\Span}{span}
\DeclareMathOperator{\diag}{diag}
\DeclareMathOperator{\Lat}{Lat}
\DeclareMathOperator{\Ext}{Ext}
\DeclareMathOperator{\vraisup}{vrai\,sup}
\def\R{\mathbb{R}}
\def\C{\mathbb{C}}
\def\N{\mathbb{N}}
\def\Z{\mathbb{Z}}
\newcommand{\F}{\mathcal{F}}
\newcommand{\la}{\lambda}
\newcommand{\ep}{\varepsilon}
\def\LimInMed{\limfunc{l.i.m.}}
\newcommand{\Rs}{\mathcal{R}}
\def\intl{\int\limits}
\def\suml{\sum\limits}

\def\prodl{\prod\limits}
\def\cupl{\bigcup\limits}
\def\sp{\sigma}
\def\ra{\rightarrow}
\def\wh{\widehat}

\def\I{\mathcal{I}}
\def\H{\mathcal{H}}

\def\K{\mathcal{K}}
\def\Np{\mathcal{N}^+}
\def\wheta{\overset{r}{\mu}}
\def\cheta{\overset{l}{\mu}}
\def\Mc{\overset{*}{M}}
\def\Si{\Sigma}
\def\si{\sigma}
\def\m{\mathfrak{m}}
\def\m{\mathfrak{m}}
\def\J{\mathcal{J}}
\def\Aess{A_{ess}}
\def\Adisc{A_{disc}}
\def\Hsp{\mathfrak{H}}
\def\PA{\mathfrak{P}}

\numberwithin{equation}{section}

\title{Indefinite Sturm-Liouville operators
$ (\sgn x) ( - \frac{d^2}{dx^2} +q(x))$
with finite-zone potentials}

\author{I. M. Karabash
\footnote{The Institute of Applied Mathematics and Mechanics, Donetsk, Ukraine. E-mail:karabashi@mail.ru,
karabashi@yahoo.com}
,
M. M. Malamud
\footnote{The Institute of Applied Mathematics and Mechanics, Donetsk, Ukraine. 
E-mail:mmm@telenet.dn.ua} 
}

\date{}

\begin{document}

\maketitle

Subj-class: Spectral Theory

MSC-class: 47E05, 34B24, 34B09 (Primary) 34L10, 47B50  (Secondary)

\begin{abstract} 
The indefinite Sturm-Liouville operator $A = (\sgn x)( -d^2/dx^2+q(x))$ is studied. It is proved that similarity of $A$ to a selfadjoint operator is equivalent to integral estimates of Cauchy integrals. Also similarity conditions in terms of Weyl functions are given. For operators with a finite-zone potential, the components $\Aess$ and $\Adisc$  of $A$ corresponding to essential and discrete spectrums, respectively, are considered. A criterion of similarity of $\Aess$ to a selfadjoint operator is given in terms of Weyl functions for the Sturm-Liouville operator  $-d^2/dx^2+q(x)$ with a finite-zone potential $q$. Jordan structure of the operator $\Adisc$ is described. We present an example of the operator 
$A = (\sgn x)( -d^2/dx^2+q(x))$ such that $A$ is nondefinitizable and $A$ is similar to a normal operator.
\end{abstract}


\section{Introduction}

The main object of the paper is a nonselfadjoint indefinite Sturm-Liouville operator
       \begin{equation}\label{0.1}
A = (\sgn x)\left( - \frac{d^2}{dx^2} + q(x)\right) =: JL, \qquad \dom A =\dom L,
       \end{equation}
where $J:\  f\to \sgn x \cdot f(\cdot)$ and $L :=  - \frac{d^2}{dx^2} + q(x)$ is a
selfadjoint  Sturm-Liouville operator on $L^2({\R})$ with a real continuous potential
$q(\cdot)$. Differential operators with indefinite weights have intensively been investigated during two last decades
(see
\cite{KKLZ84,B85,CL89,Pyat89,Sh93,CN95,F96,Vol96,
FN98,KarMFAT00,FSh00,Par03,KarMal04}).
The operator \eqref{0.1} on a finite interval subject to selfadjoint boundary conditions has discrete spectrum. The Riesz basis property of Dirichlet and other boundary value
problem for Sturm-Liouville operators with indefinite weights has been investigated in \cite{KKLZ84,B85,CL89,Pyat89,Sh93,
Par03}.

In general, the operator \eqref{0.1}    considered on $L^2({\R})$ has continuous
spectrum. In this case in place  of the Riesz basis property one  considers  the property
of similarity  either to  a normal or  to  a selfadjoint operator.

Let us recall that two
closed operators $T_1$ and $T_2$ in a Hilbert space $\Hsp$ are
called similar if there exist bounded operator $V$ with the
bounded inverse $V^{-1}$ in $\Hsp$ such that $V\dom(T_1) =
\dom(T_2)$ and $T_2 = V T_1 V^{-1}$.

Using the Krein-Langer  technique of definitizable operators in Krein spaces
\'Curgus and Langer \cite{CL89} have obtained the first result in this
direction. In particular, their result yields that
the $J$-selfadjoint  operator \eqref{0.1} is similar to a selfadjoint operator if
$L$ is a uniformly positive operator (i.e., $L \geq \delta >0$). Similarity of the operator $(\sgn x)\frac {d^2}{dx^2}$ to selfadjoint one was proved by \'Curgus and Najman \cite{CN95}. Later on, one of the
authors \cite{KarKROMSH00,KarMFAT00} reproved this result using another
approach.
More precisely, using the resolvent criterion of similarity to a selfadjoint operator
\cite{NabCr,MMMCr}
 (see also Theorem \ref{t SimCr}  below)  he proved in  \cite{KarKROMSH00,KarMFAT00} that
 the operator  $A= (\sgn x) \cdot p(-i\frac{d}{dx})$ is similar to a selfadjoint operator if and only if
the polynomial $p$ is nonnegative.

 Further, Faddeev and Shterenberg \cite{FSh00} investigated operator \eqref{0.1} with decaying
 potential. They shown, that $A$ is similar to a selfadjoint operator if $L\ge 0$ and
 $\int_{{\R}}(1+x^2)|q(x)|dx<\infty$.

The paper under consideration consists of two parts. In the first part we investigate the
operator $A$ assuming only that $q(\cdot)$ is continuous. We investigate this operator in
the framework of extension theory considering it as a (nonselfadjoint) extension of the
minimal symmetric operator
$$
A_{\min}=A^+_{\min}\oplus A^-_{\min}=L^+_{\min}\oplus(-L^-_{\min}),
$$
where $L^+_{\min}$ and $L^-_{\min}$ are minimal Sturm-Liouville operators generated by
the differential expression $L$ in $L^2({\R}_+)$ and $L^2({\R}_-)$, respectively.
Here  $\dom L^{\pm}_{\min} :=\{f\in\dom L: \  P_{\pm}f\in\dom L\}$, where $P_{\pm}$ is
the orthoprojection in $L^2({\R})$ onto $L^2({\R}_{\pm})$.

With operators $L^{\pm}_{\min}$ one associates the Weyl functions $m_{\pm}(\lambda)$
corresponding to the extensions $L^{\pm}_N$  of $L^{\pm}_{\min}, \  \dom L^{\pm}_N =
\{f\in\dom L:\  f'(\pm 0)=0\}$.

We obtain necessary and sufficient conditions of similarity in terms of the Weyl
functions $M_+(\la):=m_+(\la)$ and $M_-(\la)=-m_-(\la)$. Note, that $M_{\pm}$  are R-functions (Nevanlinna-Herglotz functions), hence the limit values  $M_{\pm}(x):= M_{\pm}(x+i0)$ exist
a.e. on $\R.$

It is worth to note that the similarity problem for the operator $A$ gives rise to two
weight estimates for the Hilbert transform in $L^2({\R})$. If fact, we show that the
following  estimate
         \begin{gather}\label{0.1Hilbert+}
\int\limits_\R \frac{\im M_{\pm}(t) + \im M_{\mp}(t) }{|M_+ (t) - M_- (t)|^2}
\left|
  g^{\pm}(t)\Sigma_{ac\pm}^\prime (t) + (H (g^{\pm}\cdot d\Sigma_{\pm})(t) \right|^2 dt
\leq K_1 \int_\R | g^{\pm} (t)|^2 d{\Sigma_{\pm}(t)} ,
          \end{gather}
(see Theorem \ref{th MNessH}) gives two necessary conditions for the operator $A$ to be
similar to a selfadjoint operator.

We conjecture that, under the assumption $\sigma_{disc} (A) = \emptyset$,
these estimates are also sufficient for similarity  to a selfadjoint
operator. 

We show that the condition  \eqref{0.1Hilbert+} yields the following necessary condition
for similarity
       \begin{equation} \label{0Mac2}
\left(\frac{1}{|\I\cap E_{\pm}|}\int_{\I}\frac{\im
M_{\pm}(t)}{|M_+(t)-M_-(t)|^2}dt\right)\cdot \left(\frac{1}{|\I\cap E_{\pm}|}\int_{\I}
\im M_{\pm}(t)dt\right) < C,
    \end{equation}
where  $\I (\subset \R)$ is any interval,  $E_{\pm}$ stand for the topological supports
of  $\im M_{\pm}, \  E_{\pm}=\supp (M_{\pm}),$  and $C$ does not depend on $\I.$

In turn, \eqref{0Mac2} implies the following weaker (and simpler) necessary condition of
similarity
       \begin{equation}\label{0.3}
\frac{ \im M_{+}(t) + \im M_- (t)}{M_+ (t)-M_-(t)}  \in L^{\infty } (\R).
       \end{equation}
 Moreover, we show that the stronger condition
      \begin{equation}\label{0.4}
\sup_{\la\in \C_+} \dfrac{|M_+ (\la ) + M_- (\la)|}{| M_+ (\la ) - M_- (\la)|} \ < \
\infty
       \end{equation}
is sufficient for similarity to a selfadjoint operator.

The second part of the paper is devoted to the spectral analysis of the operator
\eqref{0.1} with a finite-zone potential $q(\cdot)$.

Recall that a quasiperiodic (in particular, periodic) potential
$q(\cdot)={\overline{q(\cdot)}}$ is called a finite-zone potential if the spectrum
$\sigma(L)$ of the operator  $L$  has a finite number of bands (equivalently, the
resolvent set $\rho(L)$ has a finite number of gaps = forbidden zones).

We show that the operator $A=(\sgn x)\left( - \frac{d^2}{dx^2} + q(x)\right)$ with a finite-zone potential $q$ has a finite number of complex eigenvalues, and  $A$ has no
(embedding) eigenvalues on the essential spectrum $\sigma_{ess}(A)$, that is
$\sigma_p(A)\cap{ \sigma_{ess}(A)}=\emptyset$ (equivalently, the essential spectrum of $A$ coincides with  purely continuous spectrum).
Moreover, we show that the operator $A$  admits the following direct sum decomposition:
       \begin{equation*}
A=A_{disc}\dotplus A_{ess},
     \end{equation*}
where $A_{ess}$ is a part of the operator $A$ corresponding  to essential spectrum
$\sigma_{ess}(A)$ of $A$.

We summarize our main results (Theorem  \ref{t FZ Ae} and Corollary \ref{c FZ>0}) as follows:

{\it If the potential $q$ is finite-zone, then the part
$A_{ess}$ of the operator $A$ is similar to a selfadjoint operator if and only if
condition \eqref{0.3} is satisfied. Moreover,  in this case $A_{ess}$ is similar to a
selfadjoint operator with absolutely continuous spectrum.}

The main results of the paper have been announced  in our short communication
\cite{KarMal04}.

{\bf Notations.}

Throughout the paper we use the following notation.
Let $T$ be a linear operator
in a Hilbert space $\mathfrak{H}$. In what follows $\dom (T)$,
$\ker (T)$, $\ran (T)$ are the domain, kernel, range of $T$,
respectively. By $\Lat T$ we denote the set of invariant subspaces of a linear operator $T$. $ \Span \{ f_1,f_2, \dots \}$ is the closed linear hull of vectors $f_1$, $f_2$, \dots. We denote by $\sigma(T)$ the spectrum of $T$.
The discrete spectrum $\sigma_{disc} (T)$
is the set of isolated eigenvalues of finite algebraic
multiplicity; the essential spectrum is defined
by $\sigma_{ess} (T):= \sigma (T) \setminus \sigma_{disc} (T)$;
$\sigma_p (T)$ stands for the set of eigenvalues; $\rho (T)$ is the
resolvent set of $T$; $R_T \left( \lambda \right)$
is the resolvent of $T$,
\[
 R_T \left( \lambda
\right):=\left( T-\lambda I\right)^{-1} , \quad \lambda \in \rho(T) .
\]
The continuous spectrum is defined by
\[
\sigma_{c} (T) := \{ \lambda \in \C \setminus \sigma_{p} (T): \ran (T-\lambda) \neq \overline{\ran (T-\lambda)} = \Hsp \ \}.
\]
Let $\sp_{ac} (T)$ and $\sp_s (T)$ denote
the absolutely continuous and singular spectra of a selfadjoint
operator $T$ (see, for example, \cite{AG}).

Let $\mathcal{I}$ be an interval in $\R$. Let $d\Sigma$ be
a Borel measure on $\mathcal{I}$. $L^2 (\mathcal{I}, d\Sigma)$
is the Hilbert space of measurable functions $f$ on
$\mathcal{I}$ which satisfy $\int_\mathcal{I} |f|^2 d\Sigma < \infty$.
If $\mathcal{I}$ or $d\Sigma$ is fixed, we will write
$L^2 (d\Sigma)$ or $L^2 (\mathcal{I})$.
The topological support $\supp d\Sigma$ of $d\Sigma$ is the smallest closed set $S$ such that $d\Sigma (\R \setminus S) = 0$.
We denote the indicator function of a set
$S$ by $\chi_{S} ( \cdot )$;
$\ \chi_{\pm} ( t ) := \chi_{\R_{\scriptstyle \pm}} ( t )$.

We say $ f \in H(\mathcal{D})$ if $f( \cdot )$ is a holomorphic function on
a domain $\mathcal{D}$. By $\Np (\C_+)$ we denote the Smirnov class on $\C_+$
(see Subsection \ref{ss Hp}).
Suppose $\mathcal{I}$ be an interval in $\R$; then by
$\mathrm{Lip}^\alpha (\mathcal{I})$, $\alpha \in (0,1]$, we denote the
Lipschitz classes on $\mathcal{I}$ (see, for example, \cite{Gar}).

We write $f(x) \asymp g (x) \quad (x \ra x_0)$, if
the functions $\frac fg$ and $\frac gf$ are bounded in a sufficiently small
neighborhood of the point $x_0$; $f(x) \asymp g (x) \quad (x \in D)$
means that $\frac fg$ and $\frac gf$ are bounded on the set $D$.


\section{Preliminaries}

\subsection{Indefinite Sturm-Liouville operators
$ (\sgn x) ( - \frac{d^2}{dx^2} +q(x))$} \label{ss OpDef}

Denote by $J$  the multiplication operator by $\sgn x$ in the Hilbert space $L^2 (\R)$,
$J:\ f(x)\to \sgn x f(x)$. Next we consider in $L^2 (\R)$ the  differential expression
          \begin{equation} \label{e L}
L = - \frac {d^2}{dx^2} + q(x) \ ,
            \end{equation}
with a real continuous   potential $q.$ Suppose additionally that  the minimal operators
$L_{min}^+$, $L_{min}^-$ (see \cite{Naimark1}, \cite{Naimark2}) associated with \eqref{e
L} in $L^2 (\R_+)$ and  $L^2 (\R_-)$, respectively, have the deficiency indices $(1,1)$.
Denote also by $L$  the Sturm-Liouville operator generated in $L^2 (\R)$ by the
differential expression  \eqref{e L}. It is clear that $L$ is selfadjoint in $L^2 (\R).$

The main object of our paper is  an indefinite Sturm-Liouville operator
      \begin{gather}\label{e A}
A:=JL=(\sgn x) \left( - \frac {d^2}{dx^2}+ q(x) \right) \ , \qquad \dom(A):=\dom(L),  
     \end{gather}
in $ L^2( \R )$. It is easy to see that $A \neq A^*$. Indeed, the operator $A^*=LJ$ is
defined by the same differential expression \eqref{e A} on the domain, $\ \dom (A^*) = J
\dom{L} \neq \dom (A)$, containing functions discontinuous at zero together with the
first derivative.

        \begin{definition} \label{d Jsa}
Let $J$ be an signature operator on a Hilbert space $\Hsp$, $\ J=J^*=J^{-1}$. An operator
$T$ in $\Hsp$ is called J-selfadjoint if $JT = (JT)^*$.
         \end{definition}

It is clear that A is a J-selfadjoint operator. We will investigate  the operator $A$ in
the framework of extension theory of symmetric operators. For this purpose we recall the
following
            \begin{definition}[\cite{AG}] \label{d qsa}
Let $S$ be a closed symmetric operator with equal finite deficiency indices $(n,n)$, $ n
< \infty$. A closed operator $\widetilde{S}$ is called a quasi-selfadjoint extension of
$S$ if
\begin{equation*}      
S \subset \widetilde S  \subset S^* \qquad \text{and} \qquad \dim \left(
\dom(\tilde{S})/\dom(S) \right) = n.
\end{equation*}
         \end{definition}

Let $A_{min}:= A \cap A^* $, $A_{min}^\pm := \pm L_{min}^\pm $. Then
         \begin{equation}\label{2.2A}
 A_{min}= A_{min}^- \oplus A_{min}^+ , \qquad   \dom (A_{min}) := \{ y \in \dom(L):\
y(0)=y'(0) = 0 \} .
       \end{equation}
It is clear that $A_{min}$ is a simple symmetric operator  with deficiency indices
$(2,2)$ and  $A$ is its quasi-selfadjoint extension.
Indeed,
      \begin{gather} \label{e domA}
\dom (A) := \left\{ y \in \dom \left( (A_{min}^+)^* \right) \oplus \left( (A_{min}^-)^*
\right): \ y(+0) = y(-0), \ y'(+0) = y' (-0) \right\},
       \end{gather}
and $\dim (\dom(A)/\dom(A_{min}))=2.$

Note in conclusion that if  $q$ is bounded, then  $\dom(A):=\dom(L)= W^2_2(\R),$  the
Sobolev space, and  $\dom (A_{min})= W^{2,0}_2(\R) := \{ y \in W^2_2(\R):\  y(0)=y'(0) =
0 \}.$

\subsection{Weyl functions} \label{ss WF}

Recall  definition  of the Weyl functions of the Sturm-Liouville operator \eqref{e L}
assuming as before, the limit point cases at $\pm\infty.$  Denote by  $s (x,\lambda)$ and
 $c (x,\lambda)$ the solutions of
\[
-y^{\prime \prime} (x)+q(x)y(x)=\lambda y(x)
\]
obeying the following initial conditions
\[
s( 0, \lambda) = \frac{d}{dx} c ( 0, \lambda)=0, \qquad
\frac{d}{dx} s (0, \lambda) = c (0, \lambda)=1 .
\]
According to Weyl theory (see \cite{LevSar}) there exists the  function $m_{\pm}
(\lambda)$ on $\C_+ \cup \C_-$ such that
        \begin{equation} \label{e defML}
s (\cdot,\lambda) \mp m_{\pm} (\lambda) c (\cdot,\lambda)   \in L^2 (\R_\pm).
         \end{equation}
 The function $m_{\pm}$ is called \emph{the Weyl function of $L_{min}^\pm $}
corresponding to the initial condition $y'(0)=0$. The functions
       \begin{gather}
M_\pm (\lambda) := \pm m_{\pm} (\pm \lambda) \label{e Mpm}
        \end{gather}
is said to be \emph{the Weyl function
of $A_{min}^\pm $} (corresponding to
the initial condition $y'(0)=0$).

Define
       \begin{equation} \label{e def psi}
 \psi_\pm (\cdot ,\la ) :=
          \begin{cases}
- \left( s_\pm (\cdot , \pm \la ) - M_\pm (\la ) c( \cdot , \pm \la ) \right), & x \in \R_{\pm }, \\  
    0  & x \in \R_ {\mp}.
          \end{cases}
         \end{equation}
It is easily  seen that   $\psi_\pm (\cdot ,\la )   \in L^2 (\R_\pm)$  for $\la \in \C_+
\cup \C_-$  and $(A_{min}^\pm)^* \psi_\pm (x,\la ) = \la \psi_\pm (x,\la )$.

Recall that a function $m(\la)$ is called  \emph{an $R$-function (Herglotz or Nevanlinna function)}  \cite{AG,KacKr} if it is holomorphic in $\C_+ \cup \C_-$,
$$
\im \la \cdot \im m(\la) >0 \quad \text{ for }\    \la \in \C_+ \cup \C_-  \qquad \text
{and} \qquad  m(\bar{\la}) = \overline {m(\la)}.
$$
The set of all \emph{R-functions} is denoted by $(R)$ (see \cite{KacKr}).

The functions $m_{\pm}$,  as well as $M_{\pm}$ are R-functions \  (see\ \cite{LevSar}).
Moreover, it follows from \eqref{e Mpm} and the  known integral representation of $m_\pm
(\lambda )$ (see \cite{Lev, {Naimark2}}) that $M_\pm (\lambda )$ admits the following
integral representation
        \begin{gather} \label{e intM}
M_\pm (\lambda ) = \int_{\R} \frac{d \Sigma_\pm (t)}{t-\lambda}  \qquad \text{and} \qquad
\int_{\R} \frac{d \Sigma_\pm (t)}{1+|t|} < \infty.
         \end{gather}
with a (nonunique) nondecreasing scalar function $ \Sigma_\pm (t).$ Note  that
$\Sigma_\pm (t)$ in  \eqref{e intM} is uniquely determined by the following normalized
conditions: 
\[  2\Sigma_\pm (t) =  {\Sigma_\pm (t+0) + \Sigma_\pm (t-0)} \  ,
 \quad \Sigma_\pm (0)=0.
 \]
Note also that \eqref{e intM} gives a holomorphic continuation of $\ m_{\pm} (\lambda)$
to $\ \C \setminus \supp d\Sigma_{\pm}$.

Moreover, the known  asymptotic relations for $m_\pm (\cdot)$ (see \cite{Lev}) yield
          \begin{gather} \label{e asM}
M_\pm (\la ) = \pm \frac{i}{\sqrt{\pm\la}} + O\left( \frac{1}{\la}\right) , \quad (\la
\ra \infty , \
0<\delta<\arg \la < \pi - \delta ) \\
\label{e asS} \Sigma_\pm (t ) = \pm \frac{2}{\pi}\sqrt{\pm t} \ \pm \Sigma_\pm (\pm
\infty ) \ + o(1), \quad t \ra \pm\infty.
          \end{gather}
Here and below $\sqrt{z}$ is the branch
of the multifunction on the complex plane
$\C$ with the cut along $\R_+$, singled out by the condition $\sqrt{-1}=i$.
We assume that $\sqrt{\la} \geq 0$ for $\la \in [0, +\infty)$.

Consider the operator
       \begin{equation} \label{e DA2}
A_0^\pm :=(A_{min}^\pm)^*  \upharpoonright \dom (A_0^\pm), \qquad \dom (A_0^\pm) = \{ y
\in \dom \left( (A_{min}^\pm)^* \right):\  y^\prime (\pm 0) = 0 \} .
        \end{equation}
Clearly, $A_0^\pm = (A_0^\pm)^*$. The function $ \Sigma_\pm $ is
       \emph{the spectral function} of $A_0^\pm$
\cite{LevSar, {Naimark2}}. It means that the generalized Fourier transform $\F_\pm$,
defined by
        \begin{equation} \label{e F}
(\F_\pm f) (t) := \LimInMed \limits_{x_1 \rightarrow \pm \infty}
\pm \int_0^{x_1} f(x) c (x, \pm t ) dx ,
       \end{equation}
is an isometric operator from $L^2 (\R_\pm)$ onto $ L^2 (\R , d\Sigma_\pm )$. Here
$\LimInMed $ denotes the strong limit in $L^2 (\R,d\Sigma_\pm)$.

The operator $\widehat A_0^\pm := \F_\pm A_0^\pm \F_\pm^{-1}$ is the operator of
multiplication by $t$ in 
$L^2 (\R,d\Sigma_\pm (t))$, \ $ \ {\widehat A_0^\pm} :  \
g (t ) \to  \ t g(t) $ (see \cite{LevSar, {Naimark2}}). Note that $\sigma(A_0^\pm) = \supp d\Sigma_\pm$.

Suppose $f \in L^2 (\R)$. Let  $f_\pm := P_{\pm}f \in L^2 (\R_\pm)$ where $P_{\pm}$ is the
orthoprojection  in  $L^2 (\R)$ onto  $L^2 (\R_\pm).$ The following two representations
of the resolvent $\Rs_{A_0^\pm}$ are  known (see \cite{ LevSar, Naimark2}):
\begin{gather}
(\Rs_{A_0^\pm} (\la ) f_\pm ) (x)= \int_\R \frac{c(x,\pm t) \, (\F_\pm f_\pm) (t) \,
d\Sigma_\pm (t)}{t-\la} , \label{e ResRepF} \\
(\Rs_{A_0^\pm} (\la ) f_\pm ) (x)= \mp \psi_\pm (x,\la ) \int_0^{\pm x} c(s,\pm \la) f
(s) ds \mp c(x,\pm \la) \int_{\pm x}^{\pm \infty} \psi_\pm (s, \la) f (s) ds. \label{e
ResRepPsi}
\end{gather}

\subsection{Definitizable operators}

The spectral theory of linear operators
in Kre\u{\i}n spaces can be found in \cite{AJ89}, \cite{Lan82}.
Here we give some basic definitions.

Consider a Hilbert space $\Hsp$ with a scalar product
$( \cdot, \cdot )$. Let $J$ be an operator in $\Hsp$ such that
$J=J^{-1}=J^{*}$. By $[ \cdot, \cdot ]$ we define
a Hermitian sesquilinear form (s.f.) \ $ (J \cdot ,\cdot )$.
Then the pair $\K = ( \Hsp , [\cdot, \cdot])$ is a Kre\u{\i}n space
(see the literature cited above).
If $J \neq I$, then the s.f. $[\cdot, \cdot]$ is indefinite.

Let $T$ be a closely defined operator in $\Hsp$. Then
J-adjoint operator $T^{[*]}$ is defined by
\[
[Tf,g] = [f,T^{[*]}g] , \quad \ f\in D(T), \ g\in D(T^{[*]}) .
\]
Clearly, $\ T^{[*]} = J T^* J $,
where $T^*$ is the adjoint operator with respect to the scalar product
$ (\cdot, \cdot )$.
An operator $T$ is called \emph{J-selfadjoint} if $ T=T^{[*]}$. Evidently, this definition is equivalent to Definition \ref{d Jsa}
and 
\begin{gather*} \label{e T=JTJ}
T=T^{[*]} \quad \Longleftrightarrow T = J T^* J \ .
\end{gather*}

\begin{definition}[{\cite{Lan82}}]
A J-selfadjoint operator $T$ is called definitizable if
$ \rho (T) \neq \emptyset $ and there exist a real polynomial
$p$ such that
\[
[p(T)f,f] \geq 0  \quad \text{for} \quad f \in \dom (p(T)) \ .
\]
\end{definition}

Definitizable operators have spectral functions with critical points.
Thus theirs spectral properties are close to
spectral properties of selfadjoint operators in some sense
(see \cite{Lan82}).

Operators of the form \eqref{e A} are J-selfadjoint.
In this case, $\Hsp=L^2 (\R)$ and $J$
is a multiplication operator by $\sgn x$.
Such operators can be nondefinitizable.
The following theorem gives a criterion of definitizability.

\begin{theorem}[{\cite{KarDis, KarKROMSH05}}] \label{t CrDef}
Let $A=(\sgn x)(-d^2/dx^2 +q(x))$ be an operator of
the form \eqref{e A}. Then $A$ is definitizable if and only if
the sets $ \supp d\Sigma_+ $ and
$ \supp d\Sigma_- $ (see Subsection \ref{ss WF} for definitions)
are separated by a finite number of points, i.e.,
there exists a finite ordered set
\[
\{ \alpha_j \}_{j=1}^{2n-1} , \quad
-\infty = \alpha_0 < \alpha_1 \leq \alpha_2 \leq \cdots
\leq \alpha_{2n-1} < \alpha_{2n} = +\infty ,
\]
such that
\[ \supp d\Sigma_- \subset \bigcup_{k=0}^{n-1} [\alpha_{2k} , \alpha_{2k+1} ],
\qquad
\supp d\Sigma_+ \subset \bigcup_{k=0}^{n-1} [\alpha_{2k+1} , \alpha_{2k+2} ] .
\]
\end{theorem}

Several conditions of definitizability in abstract terms was given in
\cite{JL79} and \cite{JT02}.

Spectral properties of some classes of differential definitizable operators was studied in \cite{CL89,F96,CN98}; see also references
in \cite{CN98}.

\begin{definition}
An operator $T$ is called J-nonnegative if
\[
[Tf,f] \geq 0  \quad \text{for} \quad f \in \dom (T) \ .
\]
\end{definition}

Denote the \emph{root subspace} (the algebraic eigensubspace) of $T$ for $\la$ by
$\mathfrak{L}_\la(T)$, that is
\[
\mathfrak{L}_\la(T) := 
\Span \{ \ker (T-\la )^k:  \ k\in \Z_+\}.
\]

         \begin{proposition}[\cite{P55}, see also \cite{AJ89}] \label{p spJ>0}
Let $T$ be a J-nonnegative operator. Then
\begin{description}
\item[(i)] $ \sigma_p (A) \cap (\C_+ \cup \C_-) = \emptyset $.
\item[(ii)] If $\la \in \sigma_p (T) $ and $\la \neq 0$, then
the eigenvalue $\la$ is semisimple, i.e.,
$\mathfrak{L}_\lambda = \ker (T-\la)$.
\item[(iii)] If $0 \in \sigma_p (T) $, then
$\mathfrak{L}_0 = \ker T^2$
(generally, $\mathfrak{L}_0 \neq \ker T$).
\end{description}
\end{proposition}

\subsection{Finite-zone potentials} \label{ss P FZ}

Following \cite{Lev} we recall a definition of  Sturm-Liouville operator with a
finite-zone potential. Let $N \in \Z_+ :=\N \cup \{0\}$. Consider sets of real numbers
 $\{ \cheta_j \}_{j=0}^{N+1}$, $\{ \wheta_j \}_0^N$,
$\{ \xi_j \}_1^N $ such that
\[
-\infty = \cheta_0 < \wheta_0 < \cheta_1 < \wheta_1 < \dots <
\cheta_N < \wheta_N < \cheta_{N+1} = +\infty ,
\]
$ \xi_j \in [\cheta_j , \wheta_j ] $,  $j =1, \dots, N$ .
Define polynomials $R(\lambda)$, $P(\lambda)$ by
\[
P(\la ) = \prodl_{j=1}^N (\la - \xi_j ) ,
\qquad
R(\la ) = (\la - \wheta_0) \prodl_{j=1}^N
(\la - \cheta_j ) (\la- \wheta_j) .
\]
Then  there exist (see \cite{Lev}) real polynomials $S(\lambda)$ and $Q(\lambda)$ of
degrees $\deg S =N+1$ and $\deg Q = N-1$ respectively and such that
\begin{gather*}
S(\la) = \prodl_{j=0}^N (\la - \tau_j ) , \quad \tau_0 \in ( - \infty , \wheta_0 ] ,
\quad \tau_j \in [\cheta_j , \wheta_j ] ,  \quad j\in \{1, \dots, N\},
\end{gather*}
and such that the following identity  holds
       \begin{gather} \label{e PSQR}
P(\la) S(\la) - Q^2 (\la ) = R(\la ) .
\end{gather}

The function
     \begin{gather}
m_{\pm}(\la)  :=  \pm  \frac {P(\la)}{ Q(\la) \mp i \sqrt{R( \la)}} \label{e ML FZ+}
       \end{gather}
is the Weyl function 
corresponding to the Neumann boundary value problem on  $\R_{\pm}$  for some
Sturm-Liouville operator  $L=-d^2/dx^2 +q(x)$ with a bounded quasi-periodic potential
$q=\bar q$. Here the branch of the
multifunction is chosen in such a way  that  $m_{\pm}(\cdot)$  
is R-function.

\begin{definition} \label{d FZp}
A (quasi-periodic)\ potential $q=\bar q$ is called a finite-zone potential if the
Weyl functions $m_{\pm}$  of  $L_{\pm}$  defined by \eqref{e defML} admit representations
 \eqref{e ML FZ+}.
       \end{definition}

Assume  $q$ to be  a finite-zone potential. Then $q$ is an  analytic function, and the
nth derivative $\frac {d^n}{dx^n}q$ is bounded on $\R$ for any $n\in\N.$ Moreover, the
spectrum of $L=-d^2/dx^2 +q(x)$ is absolutely continuous, and
\[
\sigma (L) = \sigma_{ac} (L) = [\wheta_0 , \cheta_1 ] \cup [\wheta_1 , \cheta_2 ]
\cup \cdots \cup [\wheta_N , +\infty ) .
\]
Combining  \eqref{e ML FZ+} with \eqref{e Mpm}, we get
       \begin{gather}
M_\pm (\la ) = \frac {P(\pm \la)}{ Q(\pm \la) \mp i \sqrt{R(\pm \la)}} .
\label{e M FZ}
\end{gather}
Using \eqref{e PSQR}, we rewrite   \eqref{e M FZ}  as
    \begin{gather}
M_\pm (\la ) 
=  \frac {Q(\pm \la) \pm i \sqrt{R(\pm \la)}  }{ S (\pm \la)} . \label{e M FZ2}
       \end{gather}

\subsection{Boundary triplets and abstract Weyl functions
\label{ss BT and aWF}}

\subsubsection{Weyl  functions and spectra of proper extensions.}

Let  $\Hsp$ and $\H$ be separable Hilbert spaces.

\begin{definition}\label{dI.1}
A closed linear relation $\Theta$ in $\H$ is a closed subspace of $\H \oplus \H$.
\end{definition}

\begin{example}\label{eI.1}
For any  closed operator $B$ in $\H$  its graph $G(B)$  is a closed relation in $\H$.
\end{example}

Let $S$ be a closed densely defined symmetric operator in $\Hsp$ with equal
deficiency indices $n_+ (S) = n_- (S)$, where $n_{\pm}(S):=\dim \mathfrak{N}_{\pm i}$ and
$\mathfrak{N}_{\lambda}:=\ker(S^*-\lambda).$

\begin{definition}[\cite{AG}]\label{dI.2}
A closed extension $\tilde{S}$ of $S$ is called a proper extension if $S\subset
\tilde{S}\subset S^*$ . The set of all proper extensions is denoted by $\Ext_S$.
\end{definition}

Recall the definition of a boundary triplet.  

         \begin{definition}[\cite{GG}]\label{dI.3}
A triplet $\Pi = \{\H, \Gamma_0,\Gamma_1\}$ consisting of an auxiliary Hilbert space $\H$
and linear mappings
$\Gamma_j: \  \dom(S^*) \longrightarrow \H,\quad j\in \{0,1\},$  is called a boundary
triplet for the  operator $S^*$ if the following  conditions are satisfied:

$(i)$ The second Green  formula
\begin{equation} \label{15}
\ (S^*f,g)-(f,S^*g)=(\Gamma _1 f,\Gamma _0g)_{\H}-
                         (\Gamma_0f,\Gamma _1g)_{\H}, \qquad f,g \in \dom(S^*),
\end{equation}
holds;

$(ii)$ The mapping
$\Gamma : \  \dom(S^*) \longrightarrow \H \oplus \H, \  \Gamma f := \{\Gamma_0 f,\Gamma_1
f\}\  $
is surjective.
           \end{definition}

Definition \ref{dI.3}  allows one to describe the set $\Ext_S$ in the following way (see
\cite{DM87,DM91}).
\begin{proposition}[\cite{DM87,DM91}]\label{pI.1}
Let $\Pi=\{\H, \Gamma_0, \Gamma_1\}$ be a boundary triplet for $S^*$. Then the mapping
$\Gamma$ establishes a bijective correspondence $\tilde{S}\rightarrow
\Theta:=\Gamma(\dom(\tilde{S}))$ between the set $\Ext_S$ and the set of closed linear
relations in $\H$.
        \end{proposition}
By Proposition \ref{pI.1} the following definition is natural.
\begin{definition}\label{dI.4}
Let $\Pi=\{\H, \Gamma_0,\Gamma_1\}$ be a boundary triplet for the operator $S^*$.

$(i)$ Denote $S_{\Theta}=\tilde{S}$, if
$\Theta=\Gamma(\dom(\tilde{S}))$ that is
\begin{multline}\label{17}
S_{\Theta}:=S^*|D_{\Theta}, \quad \text{where} \quad
\dom(S_{\Theta})=D_{\Theta}:=\{f\in\dom(S^*): \{\Gamma_0f,
\Gamma_1f\}\in \Theta\}.
\end{multline}

$(ii)$ If $\Theta = G(B)$ is the graph of $B\in\mathcal{C}(\H)$ then $\dom(S_{\Theta})$
determined by the equation $\dom(S_B)=D_B:=D_{\Theta}=\ker(\Gamma_1-B\Gamma_0).$ We set
$S_B:=S_{\Theta}$.
 \end{definition}
  Let us make the following remarks.
\begin{remark}\label{rI.1}
$1)$ The deficiency indices $n_{\pm}(S)$ are equal to the
dimension of $\mathcal{H}$, i.e.,
$\dim(\mathcal{H})=n_{\pm}(S)$.

$2)$ There exist two self-adjoint extensions $S_j:=S^*|\ker(\Gamma_j)$ which are
naturally associated to a boundary triplet. According to Definition \ref{dI.4} \ 
$S_j=S_{\Theta_j}, j\in\{0, 1\}$, where $\Theta_0=\{0\}\times\H, \
\Theta_1=\H\times\{0\}$. Conversely, if $S_0$ is a self-adjoint extension of $A$, then
there exists a boundary triplet $\Pi=\{\H, \Gamma_0,\Gamma_1\}$ such that
$S_0=S^*|\ker(\Gamma_0)$.

$3)$ $\Theta $ is the graph of a closed operator $B$ iff
$\tilde{S}$ and $S_0$ are disjoint,
i.e., $\dom(\tilde{S})\cap \dom(S_0)=\dom(S)$.

$4)$ $\Theta = G(B)$ with $B\in[\H]$ iff $\tilde{S}$ and $S_0$ are transversal, i.e.,
$\tilde{S}$ and $S_0$ are disjoint and $\dom(\tilde{S}) + \dom(S_0)=\dom(S^*)$.
\end{remark}

\begin{definition}[\cite{DM92}]\label{dI.5}
A proper extension $\tilde{S}\in \Ext_S$ is called an almost solvable if there exists a
boundary triplet $\Pi=\{\H, \Gamma_0, \Gamma_1\}$ and an operator $B\in[\H]$ such that
      \begin{equation}\label{18}
\dom(\tilde{S})=\dom(S_B):=\ker(\Gamma_1-B\Gamma_0).
\end{equation}
         \end{definition}
The set of almost solvable extensions is denoted by $\mathcal{A}s_S$. Note that the class
$\mathcal{A}s_S$ is sufficiently wide. Proper extensions having two regular points
$\lambda_1, \lambda_2 \in\C $ such that $\im\lambda_1\cdot\im\lambda_2<0$ belong to
$\mathcal{A}s_S$. All quasiselfadjoint extensions are in $\mathcal{A}s_S$.

In \cite{DM87, DM91} the concept of Weyl function was generalized to an arbitrary
symmetric operator $T$ with infinite deficiency indices $n_+(A) = n_-(A)$.  Recall some
basic facts about Weyl functions.
        \begin{definition}[\cite{DM87,DM91}]\label{dI.6}
Let $\Pi=\{\H, \Gamma_0,\Gamma_1\}$ be a boundary triplet for $S^*$. The Weyl function of
$T$ corresponding to the boundary triplet $\{\H, \Gamma_0,\Gamma_1\}$ is a unique mapping
\begin{equation}\label{19}
M(\cdot):\   \rho(T_0) \longrightarrow [\H]
        \end{equation}
satisfying
\begin{equation}\label{110}
\Gamma_1f_{\lambda} = M(\lambda)\Gamma_0f_{\lambda}, \quad f_{\lambda} \in
\mathfrak{N}_{\lambda}=\ker(S^* - \lambda I), \quad \lambda \in \rho(S_0).
\end{equation}
\end{definition}
It is well known (see \cite{DM87,DM91}) that the above implicit definition of the Weyl
function is correct and $M(\cdot)$ is an operator-valued R-function obeying $0\in \rho(\im(M(i)))$(see \cite{DM95}). The
Weyl function immediately provides some information about the ''spectral properties'' of
proper extensions. We confine ourselves to the case of almost solvable extensions of the
symmetric operator $S$.

         \begin{proposition}[\cite{DM91,DM92}]\label{pI.2}
Suppose that $\Pi=\{\H, \Gamma_0, \Gamma_1\}$ is a
boundary triplet for $S^*,\ M(\cdot)$ is the corresponding Weyl
function, $\lambda \in \rho(S_0)$ and $B\in[\H]$. Then:

   $1)\ \lambda \in \rho(S_B)\ \ \text{if and only if}\ \ 0\in
            \rho(B-M(\lambda))$;

   $2)\ \lambda \in \sigma_i(S_B)\ \  \text{if and only if}\ \ 0\in
            \sigma_i(B-M(\lambda)), \ \ i\in\{p, r, c\}$.
       \end{proposition}
We demonstrate  applicability  of Proposition   \ref{pI.2}  by  describing   a discrete
spectrum of the operator $A$.
        \begin{proposition}\label{prop2.4} Let   $S := A_{\min}$  be a  (minimal) symmetric operator  defined by \eqref{2.2A} and
let  $M_{\pm}(\cdot)$  be  defined by   \eqref{e Mpm}. Then
       \begin{description}
\item[(i)] \  $\Pi=\{{\C}^2,  \Gamma_0, \Gamma_1\}$  defined by
       \begin{equation}\label{boundtriple}
\Gamma_0, \ \Gamma_1:\ \dom (A^*_{\min}) \to{\mathcal H}={\C}^2,\qquad \Gamma_0 f=
\begin{pmatrix}
f(+0),\\
f'(-0)
\end{pmatrix}, \qquad
\Gamma_1 f=
\begin{binom}
{f'(+0)}
 {-f(-0)}
\end{binom},
                    \end{equation}
forms  a boundary triplet  for the operator $S^*= A^*_{\min}$;
\item[(ii)] \  The
corresponding Weyl function  is
           \begin{equation}
M(\lambda):=M_{\Pi}(\lambda)  = \diag \bigl(- M^{-1}_+(\la), M_-(\la)\bigr);
         \end{equation}
\item[(iii)] \  The operator  $A=JL$ defined by \eqref{e A} is a quasi-selfadjoint
extension of $S$ and it is determined by
    \begin{equation}
A= S^*|\dom A, \qquad  \dom A=\ker(\Gamma_1-B\Gamma_0), \qquad \text{where} \qquad B=
\begin{pmatrix}
0&1\\
-1&0
\end{pmatrix},
    \end{equation}
that is $A=S_B$;

\item[(iv)] \  $\rho (A)\not =\emptyset$  and  $\la_0 \in \rho (A)\cap \C_{\pm}$ if and
only if $M_+ (\la_0) \neq M_- (\la_0 ).$ Moreover, $\rho (A)\cap\R = \cup_j(\alpha_j,
\beta_j)$ where $(\alpha_j, \beta_j)$ is such an interval that both $M_+$ and  $M_-$
admit holomorphic continuation trough  $(\alpha_j, \beta_j)$ and  $M_+(x+i0)\not =
M_-(x+i0),\  x\in  (\alpha_j, \beta_j).$

\item[(v)] \  The sets  $\sigma_p(A)\cap \C_{\pm}$ are at most countable with possible
limit points belonging to $\R\cup \{\infty\}$. Moreover, $\la_0 \in \sigma_p (A)\cap
\C_{\pm}$ if and only if $M_+ (\la_0) = M_- (\la_0 ).$                In   the latter
case $\dim \mathfrak{L}_{\lambda_0}(A)=\m(\la_0),$ where $\m(\la_0)$ is the multiplicity of
$\la_0 $ as a zero of the analytic function  $M_+ (\la ) - M_- (\la );$

\item[(vi)] \  The spectrum $\sigma(A)$ is symmetric with respect to the real line, that
is $\la_0 \in \sigma_p (A) \Longleftrightarrow {\overline \la_0} \in \sigma_p (A)$ and
$\dim \mathfrak{L}_{\lambda_0}(A)= \dim \mathfrak{L}_{\overline \la_0}(A)$ (equivalently
$\la_0 \in \sigma (A) \Longleftrightarrow { \la_0} \in \sigma (A^*)$  and $\dim
\mathfrak{L}_{\lambda_0}(A)= \dim \mathfrak{L}_{ \la_0}(A^*)$.
       \end{description}
              \end{proposition}
       \begin{proof}
(i)-(iii) These statements are obvious.

(iv) \  By Proposition \ref{pI.2}  
$\la_0 \in \rho (A)$ if and only if  $0 \in \rho (B-M(\la_0),$ that is
          \begin{equation}\label{3.17}
\det\bigl(B-M(\lambda)\bigr)=\det
\begin{pmatrix}
M^{-1}_+(\lambda)&1\\
-1&-M_-(\lambda)
\end{pmatrix}
=M^{-1}_+(\la)\cdot[M_+(\la)-M_-(\la)] \not =0.
           \end{equation}
Note that due to \eqref{e asM} $M_+(\cdot)$ and $M_-(\cdot)$ have different asymptotic
behavior along any semi-axes $t\cdot e^{i\varphi}, \ t>0$ with  $\varphi\in(0,\pi/2)$.
Hence $M_+ -M_- \not \equiv 0$, that is  the determinant $\det\bigl(B-M(\lambda)\bigr)$
does not vanish identically  and  $\rho (A)\not =\emptyset$.

The last statement follows from  Proposition  \ref{pI.2} and  the identity
$$
(B-M(\la))^{-1}= \frac{1}{M_+(\la)-M_-(\la)}\begin{pmatrix}
1& M_+(\la)\\
-M_+(\la)& - M_+(\la)M_-(\la)
\end{pmatrix}
$$

(v)\   By Proposition \ref{pI.2} $\sigma(S_B)\cap{\C}_{\pm}$ coincides with the set
of zeros of the determinant  $\det\bigl(B-M(\la)\bigr)$  in   ${\C}_{\pm}.$ Due to
\eqref{3.17} $\sigma(S_B)\cap{\C}_{\pm}$ coincides with the set of zeros of
$M_+(\lambda)-M_-(\lambda)$ in  ${\C}_{\pm}$ since $M_+(\la)$ has no zeros  in ${\C}_{\pm}$.  The analytic function  $M_+(\lambda)-M_-(\lambda)$  does not vanish
identically, hence it has at most countable set of zeros in both ${\C}_+$ and ${\C}_-$. The rest statements follow from analyticity of  $M_+ -M_-$ and Proposition
\ref{pI.2}.

(vi) \ Note that   $M_+(\la_0) -M_-(\la_0)=0$ yields  $M_+({\overline\la_0})
-M_-({\overline\la_0})={\overline {M_+(\la_0) -M_-(\la_0)}}=0.$ Similar implication is
valid for $j$th  derivative. This completes the proof.
       \end{proof}

\subsubsection{A functional model of a symmetric operator.}

Next  we recall  construction of a  functional model of a symmetric operator following
\cite{DM95}, \cite{MM03}.  We need only the case of the deficiency indices $(1,1)$.

Let $\Sigma (t) $ be a nondecreasing scalar function obeying the conditions
       \begin{gather}
\int_\R \frac{1}{1+t^2} d\Sigma (t) < \infty ,
\quad  \int_\R d\Sigma (t) = \infty \ , \quad
\Sigma (t) = \frac 12 (\Sigma(t-0) + \Sigma (t+0)), \quad
\Sigma (0) = 0.      \label{e S assump}
         \end{gather}
The operator of multiplication $Q_{\Sigma} : f(t) \rightarrow tf(t)$ is selfadjoint  in
$L^2 (\R,d\Sigma)$. Consider  its restriction
\[
\widehat T_\Sigma = Q_\Sigma \upharpoonright \dom (\wh T_\Sigma), \qquad
\dom (\wh T_\Sigma) = \{ f \in \dom Q_\Sigma : \int_\R f(t)d\Sigma(t) =0\}.
\]
Then $\wh T_\Sigma$ is a simple densely defined symmetric operator in
$L^2 (\R, d\Sigma )$ with deficiency indices (1,1). The adjoint operator
$\wh T_\Sigma^*$ has the form
\begin{gather*}
\dom (\wh T_\Sigma^*) = \{ f=f_Q + t(t^2+1)^{-1}h \ : \
f_Q \in \dom (Q_\Sigma), \ h \in \C \}, \quad
\wh T_\Sigma^* f = t f_Q - (t^2 +1)^{-1} h .           
\end{gather*}
Let $C\in \R.$  Define linear mappings  $\Gamma_0^{\Sigma}$, $\Gamma_1^{\Sigma,C}:\
\dom(\wh T_\Sigma^*)\to \C $ by
\begin{gather}
\Gamma_0^{\Sigma} f = h , \qquad
\Gamma_1^{\Sigma,C} f = Ch+\int_\R f_Q (t) d\Sigma(t), \label{e GammaS} \\
\text{where } \qquad f= f_Q + t (t^2+1)^{-1} h \in \dom(\wh T_\Sigma^*),
\qquad f_Q \in \dom (Q_\Sigma), \qquad h\in \C.
\notag
\end{gather}
Then $\{ \C, \Gamma_0^{\Sigma} , \Gamma_1^{\Sigma,C}\}$ is a boundary triple
for $\wh T_\Sigma^*$. The function
\begin{gather} \label{e Msc}
 M_{\Sigma,C} (\lambda) := C+\int_\R \left( \frac 1{t-\lambda} -
\frac t{1+t^2} \right) d\Sigma(t) , \quad \lambda
\in \C \setminus \supp d\Sigma,
\end{gather}
is the corresponding Weyl function of $ \wh T_\Sigma $.

\subsection{Some facts of Hardy spaces theory\label{ss Hp}}

\subsubsection{The Hilbert transform in weighted spaces}

Let us recall some facts of Hardy spaces theory
following \cite{Gar} and \cite{Koosis}.

Let  $\mu$ be a Borel  measure on  $\R$  obeying  $\int_{\R}(1 + t^2)^{-1}d\mu(t)
<\infty.$
As usual we denote by $u(\la) = {\mathcal P}_{\la}(\mu)$ its  harmonic extension (the
Poisson integral) at the point $\la=x+iy\in{\C}_+$,
   \begin{equation}\label{2.7.1}
 u(x+iy):=  {\mathcal P}_{\la}(\mu) :=(P_y * \mu)(x) := \frac{1}{\pi}
\int_{\R}\frac {y}{(x-t)^2 + y^2}d\mu(t).
       \end{equation}
For any  function $\varphi\in L^1(dt/1+t^2)$
we put ${\mathcal P}_{\la}(\varphi) := {\mathcal P}_{\la}(\mu)$
where $\mu =\varphi dx.$

Moreover, assuming that  $\int_{\R}(1+ |t|)^{-1}d\mu(t) <\infty$ one introduces  the
harmonic conjugate $\widetilde u(\cdot)$  of $u(\cdot)$ by setting
       \begin{equation}
\widetilde u(x+iy) := \frac{1}{\pi} \int_{\R}\frac {x-t}{(x-t)^2 + y^2}d\mu(t).
         \end{equation}
Here  we require the normalization $\lim_{y\ra +\infty} \wt u(x+iy) = 0.$ 
By Fatou theorem   for a.e. $x\in \R$  the limit $\lim_{y\to 0} u(x,y)=:u(x+i0)$ exists
and $u(x+i0)=\mu'(x).$ Moreover, the limit $\lim_{y\to 0} \wt u(x+iy) =: \wt u(x+i0) $ exists  
a.e. and coincides with the Hilbert transform of $\mu$, that is
            \begin{equation}  \label{2.7.2}
\wt u(x+i0) = (H \mu) (x) 
:= \frac 1\pi \lim_{\delta \ra 0} \intl_{|x - t| > \delta} \frac {1}{x - t} d\mu(t)  .
          \end{equation}
If $f \in L^p (\R)$ with  $p\in [1, \infty)$, then by definition $(Hf)(x) := (H\mu)(x)$
with $\mu=fdx.$ The operator $H$ is a unitary operator on $L^2 (\R)$.

Recall the  Helson - Szeg\"o theorem \cite{HS60} (see also \cite{Gar}).
       \begin{theorem}[Helson, Szeg\"o] \label{t HS}
Let $d\mu$ be a positive Borel measure on $\R$, finite on compact sets.  There is a
constant $K$ such that
\[
\int_\R |H f(x)|^2 d\mu (x) \leq K \int_\R |f(x)|^2 d\mu (x)
\]
for all $f \in L^2 (\R ) \cap L^2 (\R , d\mu)$ if and only if
$\mu$ is absolutely continuous, $d\mu(x) = w(x)dx$,
and
\begin{gather}
\log w(x) = u + H v, \quad u \in L^\infty (\R), \quad \|v \|_{L^\infty (\R)} < \pi /2 .
\label{e HSth}
\end{gather}
\end{theorem}

Theorem \ref{t HS}, the Helson-Szeg\"o theorem, provides a necessary and sufficient
condition  for  the Hilbert transform to  be bounded on $L^2 (d\mu)$.

Another solution  to this problem has been obtained by  Muckenhoupt \cite{Muck} and Hunt,
Muckenhoupt and Wheeden \cite{HMW}.

\begin{theorem}[Hunt, Muckenhoupt, Wheeden] \label{t HMW}
Let $d\mu$ be a positive Borel measure on $\R$, finite on compact sets.  
Then the inequality
\[
\int_\R |H f(x)|^2 d\mu (x) \leq K_2 \int_\R |f(x)|^2 d\mu (x)
\]
with $K_2$ independent of $f \in L^2 (\R ) \cap L^2 (\R , d\mu)$ holds if and only if
$d\mu(x) = w(x)dx$ and the density $w(x)$ satisfies the following condition
      \begin{gather} \label{e (Ap)}
\qquad \sup_{\I} \left( \frac 1{|\I|} \int_\I w(t) dt \right) \left( \frac 1{|\I|}
\int_\I \left( \frac 1{w(t)} \right) dt \right)
 < \infty .
\end{gather}
        \end{theorem}
Here, in  \eqref{e (Ap)}  $\sup$ is taken over the set of all (closed) intervals
$\I\subset \R.$

Condition \eqref{e (Ap)} is called the $(A_2)$-condition; we will write $w \in (A_2)$ if
\eqref{e (Ap)} is satisfied.

It is well known that the necessary part of the condition  \eqref{e (Ap)} remains valid
(with the same proof) for two-weight estimates of Hilbert transform.

More precisely,  suppose that $w_1(\cdot)$ and $w_2(\cdot)$ are two nonnegative functions
(weights) and  $E= \supp w_2={\overline E}$  is a topological support of $w_2.$  Then the
two-weight  inequality
        \begin{equation}\label{Hilbert2weight}
\int_\R |H f(x)|^2\cdot w_1(x)dx \leq K_2 \int_\R |f(x)|^2 \cdot w_2(x)dx
         \end{equation}
implies the estimate
    \begin{equation}\label{twoA2}
\qquad \sup_{\I} \left( \frac 1{|\I\cap E|} \int_\I w_1(t) dt \right) \left( \frac
1{|\I\cap E|} \int_\I \left( \frac 1{w_2(t)} \right) dt \right)
 < \infty .
     \end{equation}
In turn,  inequality \eqref{twoA2} yields
    \begin{equation}\label{twoweightsen}
\vraisup\limits_{t\in E}[w_1(x)\cdot w_2(x)^{-1}] = C < \infty.
    \end{equation}
In fact, inequalities  \eqref{twoA2} and  \eqref{twoweightsen} are not equivalent, that
is \eqref{twoA2} is stronger than \eqref{twoweightsen}.

Following \cite{NikTre} we mention  one more consequence  of  two-weight estimate
\eqref{Hilbert2weight}.
      \begin{proposition}\label{prop2.5B}
Let $w_1,w_2\ge 0$ be two nonnegative measurable functions on ${\R}$ and
$w^{-1}_2(\cdot)$ is finite a.e. on ${\R}$. Then for the two-weight estimate
\eqref{Hilbert2weight} to be valid it is necessary that
        \begin{equation}\label{PoissonEstimate}
\sup_{\lambda\in{\C}_+}{\mathcal P}_{\lambda}(w_1)\cdot{\mathcal
P}_{\lambda}(w^{-1}_2) = C <\infty.
            \end{equation}
      \end{proposition}

D. Sarason has conjectured that the converse is also true, that is  condition
\eqref{PoissonEstimate} is also sufficient for the two-weight estimate to be hold. Later
on  F. Nazarov (see \cite{NikTre}) shown that it is false.

It is easily seen (and well known) that condition \eqref{PoissonEstimate} is stronger
than   \eqref{twoA2}. Indeed, if $x$ is a middle of ${\mathcal I},\  y=|{\mathcal I}|/2$
and $\la=x+iy,$ then $|{\mathcal I}|^{-1}\chi_{{\mathcal I}}(t)\le\pi P_y(x-t)$ (cf.
\cite[Theorem VI.1.2]{Gar}).  Hence for any nonnegative $\varphi\in L^1_{loc}({\R})$
        \begin{equation}\label{2.40C}
\frac{1}{|{\mathcal I}|}\int_{{\mathcal I}}\varphi(t)dt\le\int_{{\mathcal
I}}P_y(x-t)\varphi(t)dt={\mathcal P}_{\la}(\varphi).
        \end{equation}
        
Also we will use the following result.

\begin{proposition}[cf. Theorem 4 in \cite{HS60}] \label{p A2}
Let $\{ t_j \}_{j=1}^{N}$ be a finite set of real numbers.
Assume that a (positive) weight function $w(t)$, $ t\in \R$, has the following properties: 
\begin{gather}
w(t) \asymp t^{\alpha_\infty} \qquad (|t| \ra \infty), \qquad \text{where} \quad -1 <\alpha_\infty <1, \label{e aw1}\\ 
w(t) \asymp |t-t_j|^{\alpha_j}  \qquad (t \ra t_j), \qquad \text{where} \quad -1<\alpha_j <1, \quad j=1, \dots, N,  
\label{e aw2}\\
w(t) \asymp 1  \qquad (t \ra t_0) \qquad \forall t_0 \in \R
\setminus \{ t_j \}_{j=1}^{N}.
\label{e aw3}
\end{gather}

Then $w \in (A_2) $, i.e., the weight function $w$ satisfies \eqref{e (Ap)} with $p=2$.
\end{proposition}

\begin{proof}
In this proof the letter $C$ will be used to denote a positive constant not necessarily the same at each occurrence.

If $w \not\in (A_2)$, then 
there exists a sequence of intervals $\I_n=[a_n,b_n]$, $n \in \N$, with the following properties:
\begin{description}
\item[(S1)] $\{ a_n \}_{n=1}^{\infty}$ and $\{ b_n \}_{n=1}^{\infty}$ are monotone;
\item[(S2)] there exist limits 
$
a = \lim a_n, \quad b = \lim b_n, \quad 
-\infty \leq a \leq b \leq +\infty ;
$ 
\item[(S3)] $ \displaystyle
\lim_{n \ra \infty} \left( \frac 1{|\I_n|} \int_{\I_n} w(t) dt \right)
\left( \frac 1{|\I_n|} \int_{\I_n}
\frac 1{w(t)}  dt \right)
= \infty $.
\end{description}

Let us suppose now that assumptions \eqref{e aw1}-\eqref{e aw3} hold true and let the sequences $\{ a_n \}_{n=1}^{\infty}$ and $\{ b_n \}_{n=1}^{\infty}$ have properties (S1), (S2). We will prove that property (S3) does not hold in this case, i.e.,
\begin{gather} \label{e PA}
\PA_n := \left( \frac 1{|\I_n|} \int_{\I_n} w(t) dt \right)
\left( \frac 1{|\I_n|} \int_{\I_n}
\frac 1{w(t)}  dt \right)
< C \qquad \text{for all} \quad n \in N.
\end{gather}

First note that assumptions \eqref{e aw1}-
\eqref{e aw3} yields that 
$w(\cdot) \in L^1_{loc} (\R)$ and $\frac 1{w(\cdot)} \in L^1_{loc} (\R)$. 
Hence it suffices to show \eqref{e PA} for sufficiently large $n$. 

We should consider 7 cases.

\textbf{Case 1.} Let $a=b=+\infty$ (the case $a=b=-\infty$ is similar).

By \eqref{e aw1}, 
$w(t) < C |t|^{\alpha_{\infty}}$ and 
$\frac 1{w(t)} < C |t|^{-\alpha_{\infty}}$
for sufficiently large $t>0$. 
Hence, for $n$ large enough, we have
\[
\PA_n  = \frac 1{(b_n - a_n)^2} \int_{a_n}^{b_n} w(t) dt \int_{a_n}^{b_n} \frac 1{w(t)} dt < 
C \frac 1{(b_n - a_n)^2} \int_{a_n}^{b_n} t^{\alpha_{\infty}} dt \int_{a_n}^{b_n} t^{-\alpha_{\infty}} dt. 
\]
Since $\alpha_{\infty} \in (-1,1)$, we have
\[
\PA_n < 
C \frac {(b_n^{1+\alpha_{\infty}}- a_n^{1+\alpha_{\infty}}) 
(b_n^{1-\alpha_{\infty}}- a_n^{1-\alpha_{\infty}})}
{(b_n - a_n)^2 (1+\alpha_{\infty} ) (1-\alpha_{\infty})}  < 
C \frac {
b_n^2 + a_n^2 - b_n^{1-\alpha_{\infty}} a_n^{1+\alpha_{\infty}} - 
b_n^{1+\alpha_{\infty}}a_n^{1-\alpha_{\infty}}
}
{b_n^2 + a_n^2 -2 b_n a_n} \]
(it is assumed that $a_n,b_n>0$).
By the Cauchy inequality, 
\[
b_n^{1-\alpha_{\infty}}a_n^{1+\alpha_{\infty}} + 
b_n^{1+\alpha_{\infty}}a_n^{1-\alpha_{\infty}}
> 2 b_n a_n. 
\]
Thus $\PA_n < C$ for $n$ large enough.

\textbf{Case 2.} 
Let $a=- \infty$, $b=+\infty$.

By \eqref{e aw1}, there exist a constants $a_0 <0$ and $b_0 >0$ such that 
\[
w(t) < C |t|^{\alpha_{\infty}} \quad \text{and} \quad
\frac 1{w(t)} < C |t|^{-\alpha_{\infty}} \qquad \text{for} \quad
t \in (-\infty,a_0) \cup (b_0,+\infty) .
\]
Therefore, 
\begin{gather*} 
\PA_n  < 
C  \frac 1{(b_n - a_n)^2} \left( \int_{a_n}^{a_0} |t|^{\alpha_{\infty}} dt+ \int_{a_0}^{b_0} w(t) dt + \int_{b_0}^{b_n} t^{\alpha_{\infty}} dt \right) \times \\ \times
\left( 
\int_{a_n}^{a_0} |t|^{-\alpha_{\infty}} dt+ \int_{a_0}^{b_0} \frac 1{w(t)} dt + \int_{b_0}^{b_n} t^{-\alpha_{\infty}} dt 
\right)
\end{gather*} 
for $n$ large enough.
Taking into account the fact that $\int_{a_0}^{b_0} w(t) dt < \infty$ and 
$\int_{a_0}^{b_0} \frac 1{w(t)} dt < \infty$, we get 
\begin{gather*} 
\begin{split}
\PA_n &< 
C \frac {
\left( |a_n|^{1+\alpha_{\infty}} \!- |a_0|^{1+\alpha_{\infty}} \! + C + 
b_n^{1+\alpha_{\infty}} \! - b_0^{1+\alpha_{\infty}} \right) 
\left( |a_n|^{1-\alpha_{\infty}} \!- |a_0|^{1-\alpha_{\infty}} \!+ C + 
b_n^{1-\alpha_{\infty}} \!- b_0^{1-\alpha_{\infty}} \right)
}
{(b_n - a_n)^2}  \\
&< 
C \frac {
\left( |a_n|^{1+\alpha_{\infty}} \, + \, 
b_n^{1+\alpha_{\infty}} \right) 
\left( |a_n|^{1-\alpha_{\infty}} \, + \,
b_n^{1-\alpha_{\infty}} \right)
}
{(b_n - a_n)^2} < C .
\end{split}
\end{gather*}

\textbf{Case 3.} Let $-\infty<a= b<+\infty$, $a_n \uparrow a$, and 
$b_n \downarrow a(=b)$. 

By \eqref{e aw2}-\eqref{e aw3}, there exist $\alpha \in (-1,1)$ such that
\begin{gather*}
w(t) \asymp |t-a|^{\alpha}, \quad  
\frac 1{w(t)} \asymp |t-a|^{-\alpha}, 
\qquad (t \ra a).
\end{gather*}
So, for $n$ large enough, 
\begin{gather*} 
\begin{split}
\PA_n  &< 
C  \frac 1{(b_n - a_n)^2} \left( \int_{a_n}^{a} |t-a|^{\alpha} dt+ \int_{a}^{b_n} (t-a)^{\alpha} dt \right) 
\left( 
\int_{a_n}^{a} |t-a|^{-\alpha} dt+ \int_{a}^{b_n} (t-a)^{-\alpha} dt 
\right) \\
&< C \, \frac { \left( 
(a - a_n)^{1+\alpha} + 
(b_n -a)^{1+\alpha} \right) 
\left( (a - a_n)^{1-\alpha} + 
(b_n -a)^{1-\alpha} \right) }
{\left( (b_n - a) + (a - a_n) \right)^2} \\
& = 
C \, \frac { 
(a - a_n)^2 + (b_n -a)^2 + 
(a - a_n)^{1-\alpha} (b_n -a)^{1+\alpha} 
+ (a - a_n)^{1+\alpha} (b_n -a)^{1-\alpha} }
{(a - a_n)^2 + (b_n -a)^2 + 2 (a - a_n)  (b_n -a)} \\
&< C +C \, \frac { (a - a_n)^{1-\alpha} (b_n -a)^{1+\alpha} 
+ (a - a_n)^{1+\alpha} (b_n -a)^{1-\alpha} }
{\max \{(a - a_n)^2, (b_n -a)^2 \} } <C.
\end{split}
\end{gather*} 

\textbf{Case 4.} Let $-\infty<a< b=+\infty$ and $a_n \downarrow a$ (the case  $-\infty=a< b<+\infty$, $b_n \uparrow b$ is similar). 

By \eqref{e aw1}-\eqref{e aw3}, 
\begin{gather} \label{e c4w1}
w(t) < C t^{\alpha_{\infty}}, \quad  
\frac 1{w(t)} < C t^{-\alpha_{\infty}} 
\qquad \text{for} \quad t \in (b_0 , +\infty),
\end{gather}
where $b_0$ is a certain positive constant.
Since 
\[
\int_{a_n}^b  w(t) dt \leq \int_{a}^b  w(t) dt<C 
\qquad \text{and} \qquad 
\int_{a_n} ^b  \frac 1{w(t)} dt\leq 
\int_{a} ^b  \frac 1{w(t)} dt  < C
\] 
for all $n \in N$, 
we clearly have  
\begin{gather*} 
\begin{split}
\PA_n  &< 
\frac 1{(b_n - a_n)^2} \left( \int_{a_n}^{b_0} w(t) dt+ 
\int_{b_0}^{b_n} w(t) dt \right) 
\left( \int_{a_n}^{b_0} \frac 1{w(t)} dt 
+ \int_{b_0}^{b_n} \frac 1{w(t)} dt \right)  \\
& <
C  \frac 1{(b_n - a_n)^2} \left(  C +
\int_{b}^{b_n} t^{\alpha_{\infty}} dt \right) 
\left(  C+
\int_{b}^{b_n} t^{-\alpha_{\infty}} dt \right)  \\
&< 
C  \frac 
{(b_n ^{1+\alpha_{\infty}} - b_0 ^{1+\alpha_{\infty}})(b_n ^{1-\alpha_{\infty}} - b_0 ^{1+\alpha_{\infty}})}
{b_n^2 - 2b_n a_n +a_n^2} .
\end{split}
\end{gather*}  
It follows from $\lim b_n = +\infty$ that 
$\PA_n < C$ for $n \in N$.
 
In the same way one can treat the following cases: 
\begin{description}
\item[Case 5:] $-\infty<a= b<+\infty$, $a_n \downarrow a$, and 
$b_n \downarrow a(=b)$ (
the case $a_n \uparrow a$, 
$b_n \uparrow a$ is similar);
\item[Case 6:] $-\infty<a<b=+\infty$, $a_n \uparrow a$ (
the case $-\infty=a<b<+\infty$, $b_n \downarrow b$ is analogous);
\item[Case 7:] $-\infty<a<b<+\infty$.

Thus property (S3) does not hold. This shows that $w \in (A_2)$.
\end{description}
\end{proof}

\subsubsection{The Smirnov class}

We denote by $\Np (\C_+)$ the Smirnov class on $\C_+$.
Recall that $\Np (\C_+)$ consists of holomorphic on $\C_+$ functions $U(z)$
such that $U(z)$ admits the factorization
\[
U(z) = c \, B(z) \, F(z) \, S(z) , \quad z\in \C_+,
\]
where $B$ is a Blaschke product, $F$ is an outer function,
$S$ is a singular function, $c$ is a constant, $|c|=1$
(see  \cite[Corollary II.5.6 and Theorem II.5.5]{Gar}).

The following lemmas are well known.

\begin{lemma} \label{e N++N+}
If $f,g \in \Np (\C_+)$, then $f+g \in \Np (\C_+)$.
\end{lemma}

\begin{lemma} \label{l N+}
Let $\{ t_j \}_{j=1}^{N}$ be a finite set of real numbers.
Let $U(z)$ be a holomorphic function on $\C_+$ such that
\begin{gather*}
U(z) = O(z^{\alpha_\infty}) \quad (z \ra \infty), \\
U(z - t_j) \asymp |z-t_j|^{\alpha_j}  \quad (z \ra t_j),
\ j=1, \dots, N, \\
U(z-z_0) = O (1)  \quad (z \ra z_0) \quad \forall z_0 \in (\C_+ \cup \R)
\setminus \{ t_j \}_{j=1}^{N}, \\
\end{gather*}
where $\alpha_\infty\in \R_+$, $\alpha_j \in \R_-$ , $j=1,\dots,N$.
Then $U(z) \in \Np (\C_+)$.
\end{lemma}

The proofs of these lemmas are standard.

\section{Similarity conditions}

\subsection{Characteristic functions and similarity}

Let $S$ be a symmetric operator in a Hilbert space $\Hsp$ with finite deficiency indices
$(n,n)$, $n\in \N$. Let $T$ be a quasi-selfadjoint extension of $S$. Then (see Subsection
\ref{ss BT and aWF} and \cite{DM95}) there exists a boundary triple $\{\H, \Gamma_0,
\Gamma_1\}$ for $T_{min}^*$ such that $\dom T = \ker (\Gamma_1 -B\Gamma_0)$ with some $B
\in [\H ]$, that is $T=S_B$. Let $M(\cdot )$ be the Weyl function associated with the
boundary triple $\{\H, \Gamma_0, \Gamma_1\}$. The characteristic function $\theta_T
(\cdot )$ of  almost solvable extension $T (\in \Ext_S)$ is determined and investigated in
\cite{DM92}, \cite{DM95}. In the sequel we need the following formula for the
characteristic   function   $\theta_T (\cdot )$  obtained in \cite{DM92}. It express the
$\theta_T (\cdot )$  by means of a boundary operator $B$ and the corresponding Weyl
function $M(\la)$.
         \begin{theorem}[\cite{DM92}] \label{t DM92}
Let   $\Pi= \{\H, \Gamma_0, \Gamma_1\}$  be  a boundary triple  for $S^*,$  $M(\cdot)$
the corresponding Weyl function,  $B\in [\mathcal H]$, and $E$ an auxiliary Hilbert space.
Then for any factorization $B_I:=(B-B^*)/2i=K\J K^*$  of $B_I$ with  $K\in[E,{\mathcal H}]$
and $\J =\J^*=\J^{-1}\in[E],$  the characteristic function $\theta(\la):=\theta_{A_B}(\la)$
of the extension $A_B(\in \Ext_S), \  \dom S_B=\ker(\Gamma_1-B\Gamma_0),$ admits the
following representation
        \begin{equation}\label{charfuncformula}
\theta_T(\lambda)=I + 2iK^*\bigl(B^*-M(\lambda)\bigr)^{-1}K\J.
      \end{equation}
\end{theorem}
It is shown in \cite{DM92} that if $\ker (B-B^*) = \{0\}$, then
       \begin{gather*}
\theta_T (\lambda) = \left( B-M(\lambda) \right)
\left( B^* - M(\lambda) \right)^{-1} .
       \end{gather*}
It is well known that the characteristic function $\theta_T (\la)$ obeys  the following
properties  ($\J$-properties):
     \begin{equation}\label{j-properties}
       \begin{cases}
\omega_\theta(\la):=\J-\theta_T(\la)\J \theta^*_T(\la) \  >\  0, \qquad \la \in \C_+,  \\
\omega_\theta(\la):=\J -\theta_T(\la)\J \theta^*_T(\la) \ < \  0, \qquad \la \in \C_-.
            \end{cases}
         \end{equation}
The second $\J$-form $\omega_{\theta^*}(\la):=\J-\theta_T^*(\la)\J\theta_T(\la)$ has the
same properties.

Next we recall some (sufficient) conditions of  similarity to a selfadjoint operator in
terms of the characteristic function $\theta_T (\la)$ and the corresponding $\J$-forms
$\omega_\theta(\cdot)$ and $\omega_{\theta^*}(\cdot)$.
          \begin{theorem}[\cite{MMM00}] \label{t MMM00}
Let $T$ be a solvable  extension of  $S,$  that is  $\dom T = \ker (\Gamma_1
-B\Gamma_0),\ $ with  $B\in [\H],$    $B_I:=(B-B^*)/2i=K\J K^*$ where $\J:= \sgn B_I$ and
$\pi_{\pm}:=(I\pm \J)/2$.
Suppose that $\sigma (T)\subset \R$  and  at least one of the following two conditions is
satisfied
        \begin{gather} \label{e JCond1}
(i)\qquad \max \left\{ \sup\limits_{\la \in  \C_-} \|\pi_+\theta_T^* (\la) \J \theta_T
(\la)\pi_{+} \|, \qquad  \sup\limits_{\la \in \C_+ } \|\pi_{-}\theta_T (\la) \J \theta_T^*
(\la)\pi_{-} \| \right\}  <\  \infty .
   \end{gather}
        \begin{gather} \label{e JCond2}
(ii)\qquad  \max \left\{ \sup\limits_{\la \in  \C_+} \|\pi_-\theta_T^* (\la) \J \theta_T
(\la)\pi_{-} \|, \qquad  \sup\limits_{\la \in \C_- } \|\pi_{+}\theta_T (\la) \J \theta_T^*
(\la)\pi_{+} \| \right\}  <\  \infty .
     \end{gather}
Then $T$ is similar to a selfadjoint operator $T_0$.  Moreover, if $T$ is completely
non-selfadjoint then $T_0$ has purely absolutely continuous spectrum.
         \end{theorem}
The next result  has originally been obtained in \cite{Sah}. It is immediate from Theorem
\ref{t MMM00},  other proofs can be found in \cite{NabCr, MMMCr, MMM00}.
        \begin{theorem}[\cite{Sah}] \label{t Sah}
     Let $T$ be a quasi-selfadjoint extension of $S$  and the spectrum $\sigma(T)$ is real,
$\sigma (T)\subset \R.$  If
\begin{gather} \label{e CharCond}
\ \sup\limits_{\la \in \C_+ \cup \C_-}
\|\theta_T (\la) \| <\infty ,
\end{gather}
then T is similar to a selfadjoint operator  $T_0$.  Moreover, if $T$ is completely
non-selfadjoint then $T_0$ has purely absolutely continuous spectrum.
\end{theorem}

According to the   B.S. Nagy and C. Foias  result (see \cite{SF}) condition \eqref{e
CharCond} is also necessary for a  dissipative operator $T$  to be similar to a
selfadjoint operator.

To the best of our knowledge  the most stronger sufficient condition of similarity of a
non-dissipative operator to a selfadjoint one in terms of characteristic functions, is
contained in Theorem \ref{t MMM00}. Some previous results in this direction can be found
in \cite{SF}, \cite{Sah}, \cite{MMMCr},  and  \cite{MMM00} (see also references in
\cite{MMM00}). We mention also recent publication \cite{KisFad} and \cite{Kap}.

Note that under the conditions of  all mentioned results  a  completely nonselfadjoint
part of $T$ is similar to a selfadjoint  operator $T_0 = T_0^*$ with absolutely
continuous spectrum. In this connection we mention that Kapustin \cite{Kap}  found some
sufficient conditions for an almost unitary operator $T$ to be similar to an operator
$U_{ac}\oplus T_s$ where $U_{ac}$ is an absolutely continuous unitary operator and $T_s$
is some singular almost unitary operator. Recall, that  $T$ is called  an almost
unitary operator,   if $\sigma(T) \not\supset {\mathbb D}$ and (at least one of)
non-unitary defects $I-T^*T$ and $I-TT^*$  are trace class operators.

         \begin{proposition}\label{prop2.5abscont}
Let a closed operator $T$ on $\Hsp$ be  similar to a selfadjoint operator $T_0=T_0^*,$
$\  VTV^{-1}=T_0,$ and let  $E_{T_0}(\cdot)$ be the spectral measure of   $T_0$. Then

(i) For any Borel subset $\delta\subset {\R}$  the subspace $\Hsp_T(\delta):=V^{-1} \Hsp_{T_0}(\delta),\  $ where $\Hsp_{T_0}(\delta):=E_{T_0}(\delta) \Hsp$ is a regularly and  ultra-invariant invariant
subspace for $T;$

(ii) The operator  $T(\delta) := T\lceil \Hsp_T(\delta),$\  $\dom T(\delta)= V^{-1}\dom
T_0(\delta)$ is similar to the operator  $T_0(\delta) := E_{T_0}(\delta)T;$

(iii) Suppose  additionally that  $T$ is completely non-selfadjoint,
$\sigma_p(T)=\emptyset$ and there exists a closed at most countable  set
$\{a_j\}_1^N\subset \R, \ N\le \infty,\ $ such that 
for any domain 
\begin{gather*}
\mathcal{D} := \{\la\in \C:\ |\la|>  \frac 1{\varepsilon_\infty}\} \ \cup \ \bigcup_1^N 
\{\la\in \C:\ |\la-a_j|<\ {\varepsilon}_j\}
\end{gather*}
with sufficiently small
$\varepsilon_\infty$, $\varepsilon_1$, $\varepsilon_2$, \dots, the following inequality holds
         \begin{gather} \label{e CharCond2}
\quad 
\sup\limits_{
\la \in \C_+\cup \C_- \setminus \mathcal{D}
}
\|\omega_{\theta} (\la) \| =  \sup\limits_{\la \in \C_+\cup \C_-
\setminus\mathcal{D}} \|\J - \theta_T
(\la) \J \theta_T^* (\la) \|  
\  <\ \infty .
                \end{gather}

Then the spectrum of $T_0$   is
purely absolutely continuous, that is   $T$ is similar to the  selfadjoint operator $T_0$
with absolutely continuous spectrum.
         \end{proposition}
         \begin{proof}
(i) \ It is clear that $\Hsp_T(\delta)\in \Lat T$, that is $\Hsp_T(\delta)$ is
invariant for  $T$. Moreover, $\Hsp_T(\delta)\in \Lat T$ is regularly invariant, that
is $(T-\la)^{-1}\Hsp_T(\delta)=\Hsp_T(\delta)$ since
        \begin{equation}
E_{T_0}(\delta) \Hsp = (T_0-\la)^{-1}E_{T_0}(\delta) \Hsp=V(T-\la)^{-1}V^{-1}E_{T_0}(\delta) \Hsp = V(T-\la)^{-1} \Hsp_T(\delta).
       \end{equation}
The last statement is a partial case of Proposition 5.1 from \cite{SF}, part II.

(ii)\ \  It follows from the identity  $VTV^{-1} = T_0$ that $V(T-\la)^{-1}V^{-1} =
(T_0-\la)^{-1}.$  Introducing block matrix representations of the operators $V,\
T(\delta)$ and $T_0(\delta)$ with respect to the orthogonal decompositions $\Hsp=\Hsp_T(\delta)\oplus\Hsp_T(\delta)^{\perp}= \Hsp_{T_0}(\delta)\oplus \Hsp_{T_0}({\R}\setminus\delta)$ we rewrite the above identity
in the block-matrix form
     \begin{equation}
\begin{pmatrix}
V_{11}&V_{12}\\
V_{21}&V_{22}
\end{pmatrix}
\cdot     \begin{pmatrix}
\bigl(T(\delta)-\lambda\bigr)^{-1}&T_{12}\\
0&T_{22}
    \end{pmatrix}
=
\begin{pmatrix}
\bigl(T_0(\delta)-\lambda\bigr)^{-1}&0\\
0&\bigl(T_0({\R}\setminus \delta)-\lambda\bigr)^{-1}
\end{pmatrix}
\cdot
\begin{pmatrix}
V_{11}&V_{12}\\
V_{21}&V_{22}
\end{pmatrix},
        \end{equation}
        where $V_{ij} = P_iV \lceil \Hsp_j ,\   i,j\in \{1,2\},$ \  $P_1$ is the
orthoprojection in  $\Hsp$ onto  $\Hsp_T(\delta)$ and $P_2 :=I-P_1.$ 
Hence $V_{11}\bigl(T(\delta)-\la)^{-1}=\bigl(T_0(\delta)-\la)^{-1}V_{11}$. To complete
the proof it remains to note that $\dom V_{11} = \Hsp_T(\delta),\  \ran V_{11}=  \Hsp_{T_0}(\delta)$ and  $\ker V_{11}=\{0\}$  by definition of  $V_{11}.$

(iii)  First we prove  that the operator $T_2:=P_2 T\lceil \Hsp_T(\delta)^\perp$ is
similar to the operator $T_0({\R}\setminus\delta)$. Note that  $T_2^*=T^*\lceil \Hsp_T(\delta)^{\perp}$ and
     \begin{equation}
(V^{-1})^*T^* V^*=T_0=T^*_0.
     \end{equation}
By statement (ii)  the operator $T_2^*$ is  similar the  operator $T_0({\R}\setminus\delta)$  since  $\Hsp_T(\delta)^{\perp}=V^* \Hsp_{T_0}({\R}\setminus\delta)\in 
\Lat T^*.$
Hence $T_2$ is similar to the operator  $T_0({\R}\setminus\delta) = T^*_0({\R}\setminus\delta)$ too.

Now, let $(a,b)$ be any component interval of the (open) set ${\R}\setminus\{a_j\}^{N}_1$ and $\delta=(a+\varepsilon, b-\varepsilon),\ $ $\varepsilon
> 0$. It is clear that $T$ is a coupling (see \cite{Bro, DM92, DM95}) of $T_1=T(\delta)$
and $T_2=P_2 T\lceil \Hsp_T (\delta)^{\perp}$. Therefore $\theta_T(\cdot)$ admits a
factorization (see \cite{DM92, DM95})
       \begin{equation}\label{factorization}
\theta_T(\lambda) = \theta_{T_1}(\la)\cdot\theta_{T_2}(\la) =:
\theta_1(\la)\cdot\theta_2(\la), \qquad \la\in{\C}_+\cup{\C}_-,
\end{equation}
where $\theta_j(\cdot) :=\theta_{T_j}(\cdot)$ is the corresponding characteristic
function of the operator $T_j, j\in\{1,2\}$. Since $T_2$ is similar to $T_0({\R}\setminus\delta)$, then $\theta_2(\cdot) = \theta_{T_2}(\cdot)$ admits a holomorphic
continuation through $(a+\varepsilon, b-\varepsilon)$.


It easily follows from \eqref{factorization} and the first $\J$-property of $\theta_{T_1}$
and $\theta_{T_2}$ (see \eqref{j-properties}) that
       \begin{gather}\label{2.40A}
\omega_\theta(\la) = \J-\theta_T(\la)\J\theta^*_T(\la)=
\J-\theta_{T_1}(\la)\J\theta^*_{T_1}(\la)   + \theta_{T_1}(\la)\cdot
(\J-\theta_{T_2}(\la)\J\theta^*_{T_2}(\la))\cdot \theta^*_{T_1}(\la) \\  \nonumber \ge \
\J-\theta_{T_1}(\la)\J\theta^*_{T_1}(\la)   \ge 0 \qquad \text{for} \qquad \la \in \C_+\cup
(a+\varepsilon, b-\varepsilon).
          \end{gather}
In turn, it follows from  \eqref{e CharCond2} that $\omega_{\theta}(\cdot)$ is bounded in
a small neighborhood $G_\delta^+ \ (\subset \C_+ )$ of $\delta=(a+\varepsilon,
b-\varepsilon).$ Therefore \eqref{2.40A} yields  the estimate $\sup\limits_{\la \in \,
G_\delta^+ } \|\omega_{\theta_1} (\la) \| 
\le \sup\limits_{\la \in \, G_\delta^+ } \|\omega_{\theta_{T}} (\la)\|
<\ \infty .$

On the other hand, $\theta_1 (\la)=\theta_{T_1} (\la)$ is bounded at infinity since $T_1$
is bounded. Therefore $C_+ := \sup\limits_{\la \in \C_+ } \|\omega_{\theta_1} (\la)\|
< \ \infty. $

Similarly, starting  with   \eqref{factorization} and using the second $\J$-property
\eqref{j-properties} of $\theta_{T_1}$ and $\theta_{T_2}$  we get
       \begin{gather}\label{2.40B}
\omega_\theta(\la) = \J -\theta_T(\la)\J \theta^*_T(\la)=
\J -\theta_{T_1}(\la)\J\theta^*_{T_1}(\la)   + \theta_{T_1}(\la)\cdot
(\J-\theta_{T_2}(\la)\J\theta^*_{T_2}(\la))\cdot \theta^*_{T_1}(\la) \\  \nonumber \le \
\J-\theta_{T_1}(\la)\J\theta^*_{T_1}(\la)   \le 0 \qquad \text{for} \qquad \la \in \C_-\cup
(a+\varepsilon, b-\varepsilon).
          \end{gather}
By \eqref{e CharCond2}  $\omega_{\theta}(\cdot)$ is bounded in a small neighborhood
$G_\delta^- \ (\subset \C_- )$ of $\delta=(a+\varepsilon, b-\varepsilon)$  and due to
\eqref{2.40B} so is  $\omega_{\theta_1}(\cdot).$ \  Since $\theta_1 (\la)$ is bounded at
infinity  we have  $C_- := \sup\limits_{\la \in \C_- } \|\omega_{\theta_1} (\la)\| < \
\infty. $ Summing up we get
         \begin{gather} \label{e CharCond3}
\quad          \sup\limits_{\la \in \C_+\cup \C_-} \| \theta_{T_1}(\la) \J \theta_{T_1}^*
(\la) \| \  <\ \infty .
                \end{gather}
Note that $T_1$ is completely nonselfadjoint  because  so is $T$. Since $T_1=T(\delta)$
is completely nonselfadjoint and it  is similar to the  selfadjoint operator
$T_0(\delta)$, then condition     \eqref{e CharCond3} imply absolute continuity of the
operator   $T_0(\delta)$ (see \cite{MMM00}, Theorem 1.4).  Since $(a,b)$ is any component
interval of  ${\R}\setminus\{a_j\}^{N}_1$,\  $\delta=(a+\varepsilon, b-\varepsilon),\
$ and  $\varepsilon >0$ is arbitrary, then the singular spectrum $\sigma_s(T_0)$ of $T_0$
is supported on $\{a_j\}_1^N$,  that is  $\sigma_s(T_0) \subset\{a_j\}_1^N$. Thus,
$\sigma_s(T_0)$ is at most countable, hence $\sigma_s(T_0)=\sigma_p(T_0)$. But according
to our assumption $\sigma_p(T_0)=\emptyset$ and $T_0$  is  purely absolutely continuous.
         \end{proof}
         \begin{corollary}\label{cor2.5abscont}
Let a closed operator $T$ on $\Hsp$ be  similar to a selfadjoint operator $T_0=T_0^*.$
Suppose  additionally that  $T$ is completely non-selfadjoint, $\sigma_p(T)=\emptyset$
and there exists a closed at most countable  set $\{a_j\}_1^N\subset \R, \ N\le \infty,$\
such that 
for any domain 
\begin{gather*}
\mathcal{D} := \{\la\in \C:\ |\la|>  \frac 1{\varepsilon_\infty}\} \ \cup \ \bigcup_1^N 
\{\la\in \C:\ |\la-a_j|<\ {\varepsilon}_j\}
\end{gather*}
with sufficiently small
$\varepsilon_\infty$, $\varepsilon_1$, $\varepsilon_2$, \dots, the following inequality holds
         \begin{gather} \label{e CharCond2B}
\quad    \sup\limits_{\la \in \C_+\cup \C_- \setminus \mathcal{D}}\|{\theta}_T (\la) \|
  <\ \infty .
                \end{gather}
Then $T_0$   is purely absolutely continuous, that is $T$ is
similar to the  selfadjoint operator $T_0$ with absolutely continuous spectrum.
         \end{corollary}

    \begin{remark}
It is shown in \cite{MMM00}  that  conditions  \eqref{e JCond1} and  \eqref{e JCond2} are
equivalent to each other and even are equivalent to similar conditions obtaining by
dropping  the  corresponding  orthoprojections  $\pi_{\pm}.$  Note, however that in
general condition \eqref{e CharCond3} is weaker than each  of the  (equivalent)
conditions \eqref{e JCond1} \eqref{e JCond2} and it is not sufficient for similarity to a
selfadjoint operator (cf. \cite{MMM00}).
      \end{remark}

\subsection{Characteristic functions and similarity of $J$-selfadjoint operators}

In the case of  $J$-selfadjoint operators  conditions \eqref{e JCond1},  \eqref{e JCond2}
and \eqref{e CharCond}  can be weaken.   The following  two  results are immediate from
Theorem \ref{t MMM00} and Theorem   \ref{t Sah} respectively.
     \begin{proposition}\label{corJself1}
     Suppose  additionally to the conditions of Theorem \ref{t MMM00} \; 
that  $T$ is a  $J$-selfadjoint operator.   Assume also that  $\sigma (T)\subset \R$  and
at least one of the following four conditions is satisfied
        \begin{gather} \label{e JCondJself1}
(i)\quad  
C_1 := \sup\limits_{\la \in \C_+ } \|\theta_T^* (\la) \J \theta_T (\la) \|  <\ \infty ,
\qquad  (ii)\quad C_2 := \sup\limits_{\la \in \C_- } \|\theta_T^* (\la) \J \theta_T (\la)
\| <\ \infty ,
     \end{gather}
       \begin{gather} \label{e JCondJself2}
(iii)\quad   C_3 :=   \sup\limits_{\la \in \C_- } \|\theta_T (\la) \J \theta_T^* (\la) \|
<\ \infty , \qquad (iv)\quad   C_4 :=   \sup\limits_{\la \in \C_+ } \|\theta_T (\la) \J
\theta_T^* (\la) \| <\ \infty ,
              \end{gather}
Then  $T$ is similar to a selfadjoint operator $T_0$.  Moreover, if $T$ is completely
non-selfadjoint then $T_0$ has purely absolutely continuous spectrum.
             \end{proposition}
       \begin{proof}
If two operators $T_1$ and $T_2$ are unitarily equivalent, then any characteristic
function $\theta_{T_1}(\cdot)$ of $T_1$ is at the same time the characteristic function
of $T_2$.

We prove only  that conditions $(i)$  and  $(iii)$ are equivalent and  $C_1=C_3.$  The
equivalence $(ii) \Longleftrightarrow  (iv)$ and the equality $C_2=C_4$ can be proved in
just the same way.

Since $T$ is $J$-selfadjoint it is unitarily equivalent to $T^*, \quad T^*=J T J^{-1}$.
Hence $\theta_{T} (\la) = \theta_{T^*}(\la).$
On the other hand, it easily follows from    \eqref{charfuncformula}, that
    \begin{equation*}
\theta^*_T({\overline\la})=\J\theta_{T^*}(\la)\J \  \bigl(=\J\theta_T(\la)^{-1}\J\bigr),
\qquad   \la\in\rho(T).
     \end{equation*}
This relation yields
        \begin{equation}\label{2.46A}
\theta_T^* (\overline\la) \J \theta_T (\overline\la)  =  \theta_{T^*} (\la) \J
\theta_{T^*}^* (\la)  =   \theta_T (\la) \J \theta_T^* (\la).
        \end{equation}
It follows that  $C_1=C_2.$ To complete the proof it suffices to apply Theorem  \ref{t
MMM00}.
     \end{proof}
     \begin{corollary}\label{corJself2}
     Suppose  additionally to the conditions of Theorem \ref{t Sah}
that  $T$ is a  $J$-selfadjoint operator. If  $\sigma (T)\subset \R$  and
         \begin{gather} \label{e CharCondJself}
\ \sup\limits_{\la \in \C_+} \|\theta_T (\la) \| <\infty ,
       \end{gather}
then  $T$ is similar to a selfadjoint operator $T_0$.  Moreover, if\  $T$ is completely
non-selfadjoint then $T_0$ has purely absolutely continuous spectrum.
             \end{corollary}
           \begin{remark}
Note, that four conditions (i),\  (ii),\  (iii),\  (iv)  in Proposition  \ref{corJself1}
are equaivalent. This   statement is  implied by combining identity  \eqref{2.46A} with
Proposition 1.4 from \cite{MMM00}.

 In fact, it can proved using some reasonings from \cite{MMM00} based on the
resolvent criterion (see below) that for $J$-selfadjoint operator $T$ only "half" of
either conditions \eqref{e JCond1}  or conditions  \eqref{e JCond2} is sufficient for $T$
to be similar to a selfadjoint operator. Say, the condition $\sup\limits_{\la \in \C_-}
\|\pi_+\theta_T^* (\la) \J \theta_T (\la)\pi_{+} \|\ <\ \infty$ is sufficient for $T$ to
be similar to a selfadjoint operator.
        \end{remark}

Next combining Proposition \ref{prop2.5abscont} with   Proposition \ref{corJself1} we
arrive at the following result showing that in the case of $J$-selfadjointness of the
operator $T$ condition  \eqref{e CharCond2}  can also be weaken.

         \begin{proposition}\label{prop2.5Jabscont}
Let a closed $J$-selfadjoint operator $T$ on $\Hsp$ be  similar to a selfadjoint
operator $T_0=T_0^*.$  Suppose  additionally that  $T$ is completely non-selfadjoint,
$\sigma_p(T)=\emptyset$ and there exists a closed at most countable  set
$\{a_j\}_1^N\subset \R, \ N\le \infty,\ $ such that for any domain 
\begin{gather*}
\mathcal{D} := \{\la\in \C:\ |\la|>  \frac 1{\varepsilon_\infty}\} \ \cup \ \bigcup_1^N 
\{\la\in \C:\ |\la-a_j|<\ {\varepsilon}_j\}
\end{gather*}
with sufficiently small
$\varepsilon_\infty$, $\varepsilon_1$, $\varepsilon_2$, \dots, the following inequality holds
         \begin{gather} \label{e CharCond3+}
\quad    \sup\limits_{\la \in \C_+ \setminus \mathcal{D}}\|\omega_{\theta} (\la) \| =  \sup\limits_{\la \in \C_+
\setminus\mathcal{D} } \|\J - \theta_T
(\la) \J \theta_T^* (\la) \|  
  <\ \infty .
                \end{gather}
Then $T_0$   is
purely absolutely continuous, that is   $T$ is similar to the  selfadjoint operator $T_0$
with absolutely continuous spectrum.
         \end{proposition}
     \begin{proof}
Since $T$ is $\J$-selfadjoint, then combining condition \eqref{e CharCond3} with identity
\eqref{2.46A} we get
                \begin{gather} \label{e CharCond3J}
\quad    \sup\limits_{\la \in \C_- \setminus\mathcal{D}}\|\omega_{\theta^*} (\la) \| =  \sup\limits_{\la \in \C_-
\setminus\mathcal{D}} \|\J - \theta_T^*
(\la) \J \theta_T (\la) \|  
\  <\ \infty ,
        \end{gather}
Following \cite{MMM00} it can easily be shown that both conditions  \eqref{e CharCond3+}
and   \eqref{e CharCond3J} together  yield   condition \eqref{e CharCond2}. It remains to
apply Proposition \ref{prop2.5abscont}.
               \end{proof}

        \begin{proposition}\label{prop2.5}
Let   $S := A_{\min}$  be a  (minimal) symmetric operator  defined by \eqref{2.2A} and
$A=JL$. Suppose that conditions of Proposition   \ref{prop2.4}  are satisfied and $B=
\begin{pmatrix}
0&1\\
-1&0
\end{pmatrix}$.
Then

(i) $B_I=
\begin{pmatrix}
0&-i\\
i&0
\end{pmatrix}=:\J$
and the characteristic function $\theta_A(\cdot)$ of the operator $A$ admits the
following representation
       \begin{equation}\label{formulaforcharfunc}
\theta_A(\lambda)=\frac{1}{M_-(\lambda)-M_+(\lambda)}
\begin{pmatrix}
M_+(\lambda)+M_-(\lambda)&2M_+(\lambda)M_-(\lambda)\\
2&M_+(\lambda)+M_-(\lambda)
\end{pmatrix}
    \end{equation}

(ii) The corresponding $\J$-forms are
        \begin{multline}
\omega_\theta(\lambda):=\J-\theta_A(\lambda)\J\theta^*_A(\lambda)  \\
=\J-\frac{1}{|M_+-M_-|^2}
\begin{pmatrix}
4 \cdot \im \bigl({\overline {M_+ M_-}} \cdot(M_+ +M_-)\bigr) & 4iM_+ M_- -i|M_+ +M_-|^2  \\
i|M_+ + M_-|^2-4 i {\overline {M_+ M_- }} & 4\cdot \im {\overline {(M_+ +M_-)}}
\end{pmatrix},
             \end{multline}
        \begin{multline}
\omega_{\theta^*}(\lambda):=\J-\theta_A^*(\lambda)\J\theta_A(\lambda)  \\
=\J-\frac{1}{|M_+-M_-|^2}
\begin{pmatrix}
 4\cdot \im {\overline {(M_+ +M_-)}}  & 4iM_+ M_- -i|M_+ +M_-|^2  \\
i|M_+ + M_-|^2-4 i {\overline {M_+ M_-}} &  4 \cdot \im \bigl({\overline {M_+ M_-}}
\cdot(M_+ +M_-)\bigr)
\end{pmatrix}.
             \end{multline}
(iii) The determinant $\det\theta_A(\lambda)$ defined originally on $\rho(A^*)$, admits
holomorphic continuation to the complex plane ${\C}$ and
     \begin{equation}\label{determinant}
\det\theta_A(\lambda)=1, \qquad   \lambda\in{\C}.
      \end{equation}
                         \end{proposition}

Combining Theorem \ref{t MMM00} with Proposition  \ref{prop2.5}  we arrive at the
following statement.

       \begin{corollary}\label{2.4}
Let $M_{\pm}$ be as above. Then the operator $A$ is similar to a selfadjoint operator
with absolutely continuous spectrum if the following two conditions hold
     \begin{equation}\label{2.55A}
(a)\ \  \sup_{\la\in \C_{+}}\frac{\im\bigl(M_+(\la) + M_-(\la)\bigr) + |M_+(\la)|^2
\cdot\im M_-(\la) +  |M_-(\la)|^2 \cdot \im M_+(\la)} {|M_+(\la)-M_-(\la)|^2}\ <\ \infty,
       \end{equation}
              \begin{equation}\label{2.56A}
(b) \ \  \sup_{\la\in \C_{+}}\frac{ \im M_+(\la)\cdot \im M_-(\la)}
{|M_+(\la)-M_-(\la)|^2} \  < \ \infty.
                  \end{equation}
       \end{corollary}
           \begin{proof}
Note that $\pi_{\pm}:=(I\pm \J)/2=\frac{1}{2}
\begin{pmatrix}
1&\mp i\\
\pm i&1
\end{pmatrix}$.
Setting for brevity
$\begin{pmatrix}
a&b\\
c&d
\end{pmatrix}
:=\J -\omega_{\theta^*}(\la)$ and noting that $\J - \omega_\theta(\la)=
\begin{pmatrix}
d&b\\
c&a
\end{pmatrix}$ we easily get
      \begin{equation}
\pi_+\omega_{\theta^*}(\la)\pi_+=\pi_+-\frac{1}{4}
\begin{pmatrix}
k_+&-ik_+\\
ik_+&k_+
\end{pmatrix}, \qquad
\pi_-\omega_\theta(\la)\pi_-=-\pi_- - \frac{1}{4}
\begin{pmatrix}
k_-&ik_-\\
-ik_-&k_-
\end{pmatrix},
\end{equation}
where $k_+=a-ic+ib+d$ and $k_-=a+ic-ib+d$. Hence both $k_+$ and $k_-$ are bounded in
${\C}_+$ if and only if so are $a+d=k_++k_-$ and $b-c=i(k_--k_+)$. Note that
$$
\frac{c-b}{2i} =  \frac{ |M_+(\la) + M_-(\la)|^2 - 4 Re
\bigl(M_+(\lambda)\cdot M_-(\la)\bigr)} {|M_+(\la)-M_-(\la)|^2} 
= 1 + \frac{ 8 \im M_+(\la)\cdot \im M_-(\la)}{|M_+(\la)-M_-(\la)|^2}\  ,
$$
$$
 \im\bigl(M_+(\la)\cdot M_-(\la)\cdot \overline{(M_+(\la)+M_-(\la)\bigr)} =
|M_+(\la)|^2 \cdot\im M_-(\la) +  |M_-(\la)|^2 \cdot \im M_+(\la).
$$
To complete the proof it remains to apply Theorem \ref{t MMM00}.
          \end{proof}

\begin{remark}
(i)  A weaker sufficient condition of similarity is implied by Theorem  \ref{t Sah}.
Namely, combining  Theorem  \ref{t Sah} with  formula \eqref{formulaforcharfunc}
 we conclude that the condition
      \begin{equation}\label{2.58A}
\max \left\{{\sup_{\lambda\in{\C}_+}\frac{|M_+ + M_-|}{|M_- - M_+|}\ , \qquad
\sup_{\la\in{\C}_+} \frac{1}{|M_- - M_+|}\ , \qquad \sup_{\la\in{\C}_+}
\frac{|M_+ M_-|}{|M_- - M_+|}}\right\} \ < \ \infty
    \end{equation}
is sufficient for the operator $A$ to be similar to a selfadjoint operator with
absolutely continuous spectrum.

(ii)  A counter part of identity  \eqref{determinant} for a discrete part  $A_{disc}$ of
the operator  $A$, \\ $\det\theta_{A_{disc}}(\lambda)=1,$  is immediate from  symmetry of
its spectrum (see Proposition \ref{prop2.4}(vi)). However,  identity  \eqref{determinant}
is not predictable for operators with absolutely continuous spectrum. In the latter case
$\theta_{A}(\cdot)$  is  $j$-outer function while  $\det\theta_{A}(\lambda)=1.$
                 \end{remark}

Alongside the operator   $A$ we consider  its  "dissipative and accumulative parts". More
precisely, we consider extensions $A_{\pm}$  of  $S=A_{\min}$ determined by
      \begin{gather} \label{e domA+-}
\dom (A_{\pm}) := \   \left\{ y \in \dom ( (S^*):   2y'(+0) = y'(-0) \pm iy(+0),\  2y(-0)
= y(+0) \mp iy'(-0) \right\}.
       \end{gather}

        \begin{proposition}\label{prop2.9}
Let   $S := A_{\min}$  be a  (minimal) symmetric operator  defined by \eqref{2.2A} and
let  $M_{\pm}(\cdot)$  be  defined by   \eqref{e Mpm}. Let also $\Pi=\{{\C}^2,
\Gamma_0, \Gamma_1\}$  be a boundary triplet defined by  \eqref{boundtriple}. Then
       \begin{description}
\item[(i)] \ The operators  $A_{\pm}$ defined by \eqref{e A} are  quasi-selfadjoint
extensions  of $S$ and they are  determined by
    \begin{equation}\label{A-+}
A_{\pm}= S^*|\dom A_{\pm}, \qquad  \dom A_{\pm} =  \ker(\Gamma_1-B_{\pm}\Gamma_0), \quad
\text{and} \quad B_{\pm} := \pi_{\pm}B = \frac{1}{2}
\begin{pmatrix}
\pm i&1\\
-1&\pm i
\end{pmatrix},
    \end{equation}
that is $A_{\pm}=S_{B_{\pm}},$ where
       \begin{equation}\label{B-+}
B=\begin{pmatrix}
0&1\\
-1&0
\end{pmatrix},  \quad
\J= -iB=\begin{pmatrix}
0&-i\\
i&0
\end{pmatrix}, \quad  \text{and} \quad  \pi_{\pm} := (I\pm \J)/2 = \frac{1}{2} \begin{pmatrix}
1& \mp i\\
\pm i&1
\end{pmatrix}.
        \end{equation}

\item[(ii)] \  Some of the characteristic functions  of the operators $A_{\pm}$
 are
      \begin{equation}\label{chfuncA+}
\theta_{A_{\pm}}(\la)= I- \frac{1 - M_+(\la)M_-(\la)}{\Delta_{\pm}(\la)}
\begin{pmatrix}
1& \mp i\\
\pm i&1
\end{pmatrix},
      \end{equation}
where  $\Delta_{\pm}(\la) := 1 - M_+(\la)M_-(\la)\mp 2iM_-(\la).$

\item[(iii)] \  The operator $A_+$ (resp $A_-$) is similar to a  selfadjoint operator  if
and only if
      \begin{equation}\label{2.51}
\inf_{\la \in \C_-} |1- i \Phi(\la)| =:  \varepsilon \ >\  0,
           \end{equation}
where
$$
\Phi(\cdot) := 2( M_-^{-1} (\cdot ) - M_+ (\cdot ))^{-1} \in (R).
$$
       \end{description}
              \end{proposition}
        \begin{proof}
(i) \  This statement is obvious.

 (ii)\  This statement is implied by combining  formula \eqref{charfuncformula} with  \eqref{A-+}
and  \eqref{B-+}.

(iii)\  First we note that  by \eqref{2.51}
$$
\sup _{\la\in \C_-}\left |\frac{1 - M_+(\la)M_-(\la)}{\Delta_{\pm}(\la)} \right | = \left
|\frac{1}{1- i \Phi(\la)}\right | = \frac{1}{\varepsilon}\ <\ \infty.
$$
Therefore it follows from  \eqref{chfuncA+}  that condition \eqref{2.51} is equivalent to
the boundedness of the characteristic function  $\theta_{A_+}(\cdot)$  in $\C_-.$.

Now the result is immediate from the  B.S. Nagy and Foias \cite{SF} criterion.
          \end{proof}

\subsection{Resolvent criterion}

It turns out, that in general conditions \eqref{e CharCond}, \eqref{e JCond1}, \eqref{e
JCond2} are not satisfied  for the operators of type \eqref{e A}, though such operators
may  be similar to a selfadjoint operator (see \cite{MMM00}).

Our approach is  based on the resolvent similarity criterion obtained in \cite{NabCr} and
 \cite{MMMCr} (  under an additional
assumption  this criterion was obtained in \cite{vCas}, another proof has also  been
obtained in \cite{Gom} ).

     \begin{theorem}[\cite{NabCr,MMMCr}] \label{t SimCr}
A closed operator $T$ on a Hilbert space $\Hsp$ is similar to a selfadjoint operator
if and only if $\sigma (T) \subset \R$ and for all $f \in \Hsp$ the inequalities
       \begin{equation}
\sup_{\varepsilon >0}\  \varepsilon \cdot \int_{\R} \left\| \Rs_T ( \eta +i\ep ) f
\right\|^2 d\eta \leq K_1 \left\| f\right\| ^2,          \qquad \sup_{\varepsilon
>0}\  \varepsilon \cdot \int_{\R} \left\| \Rs_{T^{*}} ( \eta
+i\varepsilon ) f\right\|^2 d\eta \leq K_{1*} \left\| f\right\| ^2,
\label{ine1}
\end{equation}
hold with constants $K_1$ and $K_{1*}$ independent of $f.$
         \end{theorem}

The following proposition is immediate from Theorem \ref{t SimCr}.
           \begin{proposition} \label{p SimCrJ}
A J-selfadjoint operator $T$ on a Hilbert space $\Hsp$ is similar to a selfadjoint
operator
 if and only if $\sigma (T) \subset \R$ and  the following inequality holds
              \begin{equation}  \label{e ine1J}
\sup_{\varepsilon >0}\ \varepsilon\cdot  \int_{\R} \left\| \Rs_T ( \eta +i\ep ) f
\right\|^2 d\eta \leq K_1 \left\| f\right\| ^2,  \qquad  f \in \Hsp,
        \end{equation}
with a constant $K_1$ independent of $f.$
          \end{proposition}
     \begin{proof}
If T is a J-selfadjoint operator, then $T^* = JTJ$ and the second inequality  in
\eqref{ine1} is equivalent to  the first one.
          \end{proof}

In the case of a bounded operator $T$  we can slightly clarify  Theorem  \ref{t SimCr} in
the following way.
       \begin{corollary} \label{corCrit}
Let  $T=T_1 + iT_2$  where $T_1=T_1^*$ and  $T_2=T_2^*\in [\Hsp].$ Then $T$  is
similar to a selfadjoint operator  if and only if $\sigma (T) \subset \R$ and for all $f
\in \Hsp$ the inequalities
        \begin{equation} \label{2.31B}
\sup_{0< \varepsilon <2\|T_2\|} \varepsilon \int_{\R }\left\| \Rs_T ( \eta +i\ep ) f
\right\|^2 d\eta \leq K_1 \left\| f\right\| ^2,          \qquad \sup_{0< \varepsilon
<2\|T_2\|}\varepsilon \int_{\R} \left\| \Rs_{T^{*}} ( \eta +i\varepsilon ) f\right\|^2
d\eta \leq K_{1*}
\left\| f\right\| ^2, 
     \end{equation}
hold  with constants $K_1$ and $K_{1*}$ independent of $f\in \Hsp.$

In particular,   a  bounded operator $T$ on  $\Hsp$ with $\sigma (T) \subset \R$  is
similar to a selfadjoint operator  if and only if  inequalities  \eqref{2.31B} are valid
with $2\|T_2\|$  replaced for any $\varepsilon_0>0.$
         \end{corollary}
      \begin{proof}
(i)  It is clear that
     \begin{equation*}\label{2.31C}
(T-z)^{-1}  =  (T_1-z)^{-1}  - (T_1-z)^{-1}\cdot T_2\cdot (T-z)^{-1} , \qquad   z \in
\C_+.
       \end{equation*}
It follows that
       \begin{gather*}
 \|(T-z)^{-1}f\|^2 \le 2 \|(T-z)^{-1}f\|^2  +  2 \|(T_1-z)^{-1}\cdot
T_2\cdot (T-z)^{-1}f\|^2  \\
\le  \|(T_1-z)^{-1}f\|^2 + \frac{2\|T_2\|^2}{|\im z|} \|(T-z)^{-1}f\|^2  \qquad z \in
\C_+ , \quad  f\in \Hsp.
      \end{gather*}
In turn, this inequality yields $\  \|(T-z)^{-1}f\| \le  2 \|(T_1-z)^{-1}f\| \ $  for
$\im z > 2\|T_2\|.$  Hence

It is known that $\|(T-z)^{-1}\|\le(|z|-\|T\|)^{-1}$ for $|z|>\|T\|$. Hence
      \begin{equation}\label{2.32B}
 \sup_{ \varepsilon \ge 2\|T_2\|}\  \varepsilon \cdot \int_{\R }\left\| \Rs_T ( \eta +i\ep ) f
\right\|^2 d\eta \leq  4 \varepsilon \cdot  \int_{\R} \left\| \Rs_{T_1} ( \eta +i\ep ) f
\right\|^2 d\eta  = 4\pi \left\| f\right\| ^2.
        \end{equation}
Combining this inequality  with the first of inequalities \eqref{ine1} we arrive at the
first of inequalities \eqref{2.31B}. The second one can be  proved similarly.
             \end{proof}



\begin{remark}
If  $T$ is a closed unbounded  operator, then conditions   \eqref{ine1} and \eqref{2.31B}
are not equivalent, in general. In fact,  there exists an operator $T$ such that:
\begin{description}
\item[(i)] \  $\sigma(T) \subset \R$ ;
\item[(ii)] \ conditions \eqref{2.31B}  are fulfilled;
\item[(iii)] conditions \eqref{ine1} do not hold and, consequently, $T$ is not similar to a self-adjoint operator.
\end{description}

Setting
       \begin{equation*}
\begin{cases} w(x) = 1, &  x \in (-\infty, -1)\cup (1, +\infty),  \\
               w(x)= -1,  &  x \in (-1,1) ,
\end{cases}
         \end{equation*}
consider the operator $ D_w = -i \frac 1{w(x)}\, \frac {d}{dx}$ in $L^2 (\R).$ It was shown in
\cite{KarKROMSH00} that the operator $D_w$ has the properties (i) and (iii). 
It is not difficult to check that conditions \eqref{2.31B} are fulfilled for $D_w$.
       \end{remark}

\section{Eigenvalues and their multiplicities}

In \cite{KarDis, KarKROMSH05} the functional model for
J-selfadjoint
quasiselfadjoint operator was given.
The model is based on classical Sturm-Liouville spectral theory
and the functional model for a symmetric operator given in
Section \ref{ss BT and aWF}.

Let $\Sigma_\pm$ be
the spectral functions of $A_0^\pm$ (see \eqref{e intM}).
It follows from \eqref{e asS} that they satisfy
\eqref{e S assump}. Let $C_{\pm} := \displaystyle
\int_\R \frac {t}{1+t^2} d\Sigma_{\pm}$.
We denote $\wh \Gamma_0^\pm := \Gamma_0^{\Sigma_\pm}$,
$\wh \Gamma_1^\pm := \Gamma_1^{\Sigma_\pm , C_\pm}$.
From the definition of $C_\pm$ and \eqref{e GammaS}, we get
\begin{gather*}
\wh \Gamma_1^\pm f = C_\pm \wh \Gamma_0^\pm f +
\int_\R \left( f (t) -
\frac {t \, \wh \Gamma_0^\pm f }{t^2+1} \right) d\Sigma_\pm (t) =
\int_\R f(t) d\Sigma_\pm (t)
\end{gather*}
for $ f \in \dom ( \wh T_{\Sigma_\pm})$.
Consider the operator $\widehat A $
in $L^2 (d\Sigma_+) \oplus L^2 (d\Sigma_-)$ defined by
\begin{gather}
\widehat A =
\wh T_{\Sigma_+}^* \oplus \wh T_{\Sigma_-}^* \upharpoonright
\dom (\widehat A) \label{e A hat} , \\
\dom (\widehat A) = \{ \ f=f_+ + f_- \ : \ f_\pm \in
\dom (\wh T_{\Sigma_\pm}^*) , \
\wh \Gamma_0^+ f_+ = \wh \Gamma_0^- f_- , \
\wh \Gamma_1^+ f_+ = \wh \Gamma_1^-  f_- \ \}
\notag
\end{gather}
(for the definition of $\wh T_{\Sigma_\pm}$ see Section
\ref{ss BT and aWF}).

\begin{proposition}[\cite{KarDis, KarKROMSH05}]   \label{p SLModel}
The operator $A$ of type \eqref{e A} is unitary equivalent
to the operator $\widehat A $.
Moreover,
\begin{gather} \label{e FAF=whA}
(\F_- \oplus \F_+ ) A (\F_-^{-1} \oplus \F_+^{-1} ) =
\widehat A  \ .
\end{gather}
\end{proposition}

Note that we can write the Weyl functions of $A$ in the form
\[
\  M_{\pm} (\lambda) = M_{ \Sigma_{\scriptstyle \pm},
C_{\pm}} (\lambda) \ ,  \quad
 \lambda\in \C \setminus \supp d\Sigma_{\pm} \
\]
(see \eqref{e Msc} for the definition of $M_{ \Sigma_{\scriptstyle \pm}, C_{\pm}}$).

Now we classify eigenvalues of $\wh T_\Sigma^*$.
Let us introduce the following mutually disjoint sets:
\begin{gather*}
\mathfrak{A}_0 (\Sigma) = \left\{ \lambda \in \sigma_c (Q_\Sigma) \ :
\ \int_\R |t-\lambda|^{-2} d\Sigma (t) = \infty \right\} , \\
\mathfrak{A}_r (\Sigma) = \left\{ \lambda \not \in \sigma_p (Q) \ :
\ \int_\R |t-\lambda|^{-2} d\Sigma (t) <\infty \right\}, \qquad
\mathfrak{A}_p (\Sigma) = \sigma_p (Q_\Sigma).
\end{gather*}
Observe that $\C = \mathfrak{A}_0 (\Sigma) \cup \mathfrak{A}_r (\Sigma) \cup
\mathfrak{A}_p (\Sigma) $ and
\begin{gather}
\mathfrak{A}_0 (\Sigma) = \left\{ \lambda \in \C \ :
\ \ker (A_\Sigma^* - \lambda) = \{ 0 \} \ \right\} , \notag \\
\mathfrak{A}_r (\Sigma) = \left\{ \lambda \in \C :
\ \ker (A_\Sigma^* - \lambda) =
\{ c(t-\lambda)^{-1}, \; c \in \C \} \right\}, \label{e Arker} \\
\mathfrak{A}_p (\Sigma) = \left\{ \lambda \in \C \ :
\ \ker (A_\Sigma^* - \lambda) = \{ c\chi_{ \{\lambda\} } (t), \; c \in \C \}
\right\} \label{e Apker}.
\end{gather}

The following theorem gives a description of the point spectrum of
$\widehat A $.

\begin{theorem}[\cite{KarDis, KarKROMSH05}] \label{t s p}
Let $\widehat A $ be given by \eqref{e A hat}.
\begin{description}
\item[1)]
If $\lambda \in \mathfrak{A}_0 (\Sigma_+) \cup \mathfrak{A}_0 (\Sigma_-)$,
then $ \lambda \not \in \sigma_p (\widehat A)$.
\item[2)]
If $\lambda \in \mathfrak{A}_p (\Sigma_+) \cap \mathfrak{A}_p (\Sigma_-)$,
then
\begin{description}
\item[(i)] $\lambda$ is an eigenvalue of $\widehat A$; the geometric multiplicity of $\lambda$ equals 1;
\item[(ii)] the eigenvalue $\lambda$ is simple (i.e., the algebraic and geometric multiplicities are equal one) iff at least one of the following conditions is not fulfilled:
\begin{gather} \label{e S-S=S-S}
\Sigma_- (\lambda+ 0 ) - \Sigma_- (\lambda -0) =
\Sigma_+ (\lambda + 0 ) - \Sigma_+ (\lambda -0) ,
\\      
\label{e Int+<Inf}
\int_{\R \setminus \{\lambda \}} \frac 1{|t-\lambda |^2} \ d \Sigma_+ (t) < \infty ,
\\      
\label{e Int-<Inf}
\int_{\R \setminus \{\lambda \}} \frac 1{|t-\lambda |^2} \ d \Sigma_- (t) < \infty ;
\end{gather}
\item[(iii)] if conditions \eqref{e S-S=S-S}, \eqref{e Int+<Inf} and
\eqref{e Int-<Inf} hold true, then the algebraic multiplicity of $\lambda$ equals
the greatest number $k$ ($2 \leq k \leq \infty$) such that the following conditions
\begin{gather} \label{e frj<inf}
\int_{\R \setminus \{\lambda \}} \frac 1{|t-\lambda|^{2j}} \ d \Sigma_- (t) < \infty ,
\qquad
\int_{\R \setminus \{\lambda \}} \frac 1{|t-\lambda|^{2j}} \ d \Sigma_+ (t) < \infty ,
\\            
 \label{e frj=frj}
\int_{\R \setminus \{\lambda \}} \frac 1{(t-\lambda )^{j-1}} \ d\Sigma_- (t) =
\int_{\R \setminus \{\lambda \}} \frac 1{(t-\lambda )^{j-1}} \ d\Sigma_+ (t) ,
\end{gather}
are fulfilled for all $\ j \in \N \cap [ 2 , k-1 ]$.
\end{description}
\item[3)]
Assume that
$\lambda \in \mathfrak{A}_r (\Sigma_+) \cap \mathfrak{A}_r (\Sigma_+)$.
Then $\lambda \in \sigma_p (\widehat A)$ iff
\begin{gather} \label{e fr1=fr1}
\int_\R \frac 1{t-\lambda} d\Sigma_+ (t) =
\int_\R \frac 1{t-\lambda} d\Sigma_- (t) \ .
\end{gather}
If \eqref{e fr1=fr1} holds true, then the geometric multiplicity of $\lambda$
is one and the algebraic multiplicity is the greatest number $k$
($1 \leq k \leq \infty$) such that
the following conditions
\begin{gather}
\int_{\R } \frac 1{|t-\lambda |^{2j}}
\ d \Sigma_- (t) < \infty ,
\qquad
\int_{\R } \frac 1{|t-\lambda |^{2j}} \ d \Sigma_+ (t) < \infty ,
\label{e frj<inf r}
\\
\int_{\R } \frac 1{(t-\lambda )^{j}} \ d\Sigma_- (t) =
\int_{\R } \frac 1{(t-\lambda )^{j}} \ d\Sigma_+ (t)
\label{e frj=frj r}
\end{gather}
are fulfilled for all $\ j \in \N \cap [1, k] $.
\item[4)]
If $\lambda \in \mathfrak{A}_p (\Sigma_+) \cap \mathfrak{A}_r (\Sigma_-)$
or $\lambda \in \mathfrak{A}_p (\Sigma_-) \cap \mathfrak{A}_r (\Sigma_+)$,
then $\lambda \not \in \sigma_p (\widehat A)$.
\end{description}
\end{theorem}

It follows from Theorem \ref{t s p}
(as well as from Proposition \ref{prop2.4}) that
\begin{gather}
 \left\{ \lambda \in \rho ( Q_{\Sigma_+} ) \cap \rho ( Q_{\Sigma_-} ) \ : \
M_+ (\lambda) = M_- (\lambda) \right\} =
\sigma (\widehat A) \cap \rho ( Q_{\Sigma_+} \oplus Q_{\Sigma_-})
\subset \sigma_p (\widehat A). \label{e s in rho}
\end{gather}
It is easy to see that \eqref{e s in rho} and Theorem \ref{t s p}
yield the following description of the essential and discrete
spectra.

\begin{proposition}[\cite{KarDis,KarKROMSH05}] \label{p e and d}
\begin{description}
\item[1)] $ \sigma_{ess} (\wh A) = \sigma_{ess}
(Q_{\Sigma_{+}})
\cup \sigma_{ess} (Q_{\Sigma_{-}}) ;$
\item[2)] $ \sigma_{disc} (\wh A) = \left( \sigma_{disc}
(Q_{\Sigma_{+}}) \cap \sigma_{disc} (Q_{\Sigma_{-}}) \right) \
\cup \
\{ \lambda \in
\rho (Q_{\Sigma_+}) \cap \rho (Q_{\Sigma_-})
\ : \ M_{+} (\lambda) = M_{-} (\lambda) \} $;
\item[3)] the geometric multiplicity equals 1 for all
eigenvalues of $\wh A$;
\item[4)] if $\lambda_0 \in  \left( \sigma_{disc}
(Q_{\Sigma_{+}}) \cap \sigma_{disc} (Q_{\Sigma_{-}}) \right)$,
then the algebraic multiplicity of $\lambda_0$
is equal to the multiplicity of $\lambda_0$ as a zero
of the holomorphic function
$\displaystyle \frac 1{M_{+} (\lambda)} - \frac 1{M_{-} (\lambda)}$;
\item[5)] if $\lambda_0 \in \rho (Q_{\Sigma_{+}}) \cap
\rho (Q_{\Sigma_{-}})$ then the algebraic multiplicity
of $\lambda_0$ is equal to the multiplicity of
$\lambda_0$ as zero of the holomorphic function
$M_{+} (\lambda) - M_{-} (\lambda)$.
\end{description}
         \end{proposition}
          \begin{proposition} \label{p Res} 
Let  $A$ be the operator defined  by  \eqref{e A} and  $\lambda_0 \in \C \setminus \R$.Then

\begin{description}
\item[(i)] \  $\rho (A)\not =\emptyset$  and  $\la_0 \in \rho (A)\cap \C_{\pm}$ if and
only if $M_+ (\la_0) \neq M_- (\la_0 )$.

\item[(ii)] The resolvent of A has the following form
          \begin{equation} \label{e RA}
\Rs_A (\la ) f(\cdot) = \Rs_{A_0^-\oplus A_0^+}(\la) f(\cdot) + G^-(\la)\psi_-(\cdot,
\la) + G^+(\la)\psi_+(\cdot, \la) ,
       \end{equation}
\begin{equation} \label{e G}
G^- (\la)=G^+ (\la ) = \frac{1}{M_+ (\la) - M_- (\la)}
\int_\R \frac{g^-(t)d\Sigma_-(t)-g^+(t)d\Sigma_+(t)}{t-\la} ,
               \end{equation}
where $g^\pm (t)=(\F_\pm f_\pm) (t)$, $f_\pm := P_{\pm}f \in L^2(\R_\pm)$, and  
$P_{\pm}$ is the     orthoprojection in $L^2(\R)$  onto  $L^2(\R_\pm).$
      \end{description}
               \end{proposition}

\begin{proof}
(i) This statement  has already been proved in Proposition \ref{prop2.4}.

(ii)  Let now  $\la \in \rho (A)$ and
 $y(\cdot, \la) = (A-\la I)^{-1}f(\cdot) $. It  means that
$y \in \dom (A_{min}^*)$ and $y$ is
a solution of the equation
\begin{equation} \label{e SLy}
(\sgn x)(-y^{\prime\prime} (x)+q(x)y(x))-\la y(x)=f(x)
\end{equation}
subject to "glue" boundary conditions
\begin{gather}
\label{e BoundCondEVA}
y(-0)=y(+0), \qquad y^\prime (-0)=y^\prime (+0) \ .
\end{gather}
Hence,
\begin{equation*}
y(x, \la)= (\Rs_{A_0^- \oplus A_0^+} (\la )f)(x) + G^-(\la)\psi_-(x,\la) +
G^+(\la)\psi_+(x,\la) ,
\end{equation*}
where $G^\pm (\la)$ are the scalar functions. 
It is clear  that
         \begin{equation*}
y(\pm 0, \la)= (\Rs_{A_0^\pm} (\la )f_\pm)(\pm 0) + G^\pm (\la)\psi_\pm (0,\la) .
\end{equation*}
By \eqref{e def psi}, we get
\begin{equation} \label{e psi0}
\psi_\pm (0,\la) =M_\pm (\la), \qquad
\frac{d}{dx} \psi_\pm (0,\la)=-1 .
\end{equation}

Resolvent representation \eqref{e ResRepF} yields
\begin{gather*} 
(\Rs_{A_0^\pm} (\la ) f_\pm ) (\pm 0)=
\int_\R \frac{ g^\pm (t) d\Sigma_\pm (t)}{t-\la} .
\end{gather*}
It follows from $ \Rs_{A_0^\pm} (\la ) f_\pm \in D(A_2^\pm)$
and \eqref{e DA2} that
$\frac{d}{dx}(\Rs_{A_0^\pm} (\la ) f_\pm)_{x=\pm 0} = 0$.
Taking into account \eqref{e psi0},
we see that conditions \eqref{e BoundCondEVA} take
the form
\[
\begin{cases}
\dint_\R \frac{ g^- (t) d\Sigma_- (t)}{t-\la} + G^- (\la )M_-(\la )=
\dint_\R \frac{ g^+ (t) d\Sigma_+ (t)}{t-\la} + G^+ (\la )M_+(\la ) \\
G^- (\la ) =G^+ (\la )
\end{cases} .
\]
Since $M_+ (\la ) \neq M_- (\la )$, problem \eqref{e SLy}-\eqref{e BoundCondEVA} has the
unique solution $y \in \dom (A_{min}^*)$ and  it  admits  a   representation \eqref{e
RA}-\eqref{e G}.
         \end{proof}

Next we clarify Proposition   \ref{p Res} in the  case  of  $J$-nonnegative operator $A$,
i.e., if $L \geq 0.$
          \begin{proposition} \label{p spJ-nonneg}
If the operator $ L = - d^2/dx^2 + q(x)$ is nonnegative,
then the spectrum of the operator $A=JL$ is real. 
          \end{proposition}
\begin{proof}
Since $L\ge 0$ we have $A^+_{\min}=L^+_{\min}\ge 0$ and $A^-_{\min}=-L^+_{\min}\le 0$. It
is known that the Friedrichs extension $L^{\pm}_F$ of $L^{\pm}_{\min}$ is generated by
the Dirichlet boundary value problem,  that is
     \begin{equation}
L^{\pm}_F=(L^{\pm}_{\min})^*\lceil\dom(L^{\pm}_F), \quad
\dom(L^{\pm}_F=\{f\in\dom(L^{\pm}_{\min})^*:f(0)=0\}.
      \end{equation}
Setting $\Gamma^{\pm}_0 f=f(\pm 0)$ and $\Gamma^{\pm}_1 f=\pm f'(\pm 0)$ we obtain a
boundary triplet $\Pi_{\pm}=\{{\C},\Gamma^{\pm}_0,\Gamma^{\pm}_1\}$ for
$(L^{\pm}_{\min})^*$ such that $\ker\Gamma^{\pm}_0=\dom(L^{\pm}_F)$. Therefore the
corresponding Weyl function  $m^{\pm}_F$ belongs to the Krein-Stielties class $S^-$ (see
\cite{DM95}). Hence, it admits the following integral representation (see \cite{KacKr}).
      \begin{equation}\label{3.23A}
m^{\pm}_F(\lambda)=C_{\pm}+\lambda\int^{\infty}_0\frac{d\sigma_{\pm}(t)}{t-\lambda},
\qquad   \int^{\infty}_0\frac{d\sigma_{\pm}(t)}{1+t}<\infty,
        \end{equation}
with $C_{\pm}\le 0$. On the other hand, it follows from definitions      that
       \begin{equation}\label{3.24A}
-M^{-1}_+(\lambda)=-m^{-1}_+(\lambda)=m^+_F(\lambda), \qquad
M^{-1}_-(\lambda)=-m^{-1}_-(\lambda)=m^-_F(-\lambda)
       \end{equation}
Combining  \eqref{3.23A} and  \eqref{3.24A} we get
      \begin{multline}
M^{-1}_-(\lambda)-M^{-1}_+(\lambda)=m^-_F(-\lambda)+m^+_F(\lambda)\\
=\la \bigl[\frac{C_- +
C_+}{\lambda}+\int^{\infty}_0\frac{d\sigma_+(t)}{t-\lambda}-\int^{\infty}_0\frac{d\sigma_-(t)}{t+\lambda}\bigr]
=:\lambda{\widetilde M}(\lambda), \qquad \qquad
      \end{multline}
where ${\widetilde M}(\cdot)\in(R)$ since $C_{\pm}\le 0$. To complete the proof it
remains to note that
       \begin{equation}
M_+(\lambda)-M_-(\lambda)= M_+(\lambda)\cdot [M^{-1}_-(\lambda)-M^{-1}_+(\lambda)]\cdot
M_-(\lambda) = M_+(\lambda)\cdot\lambda{\widetilde M}(\lambda)\cdot M_-(\lambda)\not =0
      \end{equation}
for $\lambda\in{\C}_{\pm}$, since $M_{\pm},{\widetilde M}\in(R)$.
         \end{proof}
      \begin{remark}
      \begin{description}
\item[(i)] Statement (i) of Proposition   \ref{p Res}  is implied by  \eqref{e s in rho}.
However, we presented an elementary proof based on Proposition \ref{pI.2}.

\item[(ii)]  Note that Proposition \ref{p spJ-nonneg}  follows immediately from Proposition \ref{p Res} and Proposition
\ref{p spJ>0}. However, we presented another  proof that is in a spirit of our paper and
demonstrates applicability of Weyl function technic. 
Note also  that in turn,
Proposition \ref{p spJ>0} can be proved by using Weyl function technic similarly to the
proof of Proposition \ref{p spJ-nonneg}.
\end{description}
       \end{remark}

\section{Similarity conditions for the operator $A$. General case. }

\subsection{Similarity criterion  in terms of  Weyl functions.}

In the sequel we write  $\la = \eta + i \ep$, that is  $\eta = \re \la$, $\ep = \im \la$.

Combining Proposition~\ref{p SimCrJ} and Proposition~\ref{p Res}, we arrive at the
following criterion.

          \begin{theorem} \label{th CrSLSim}
The operator $A= (\sgn x) (-d^2/dx^2 +q(x))$
is similar to a selfadjoint operator if and only if
for all $\ep >0 $ and $g^\pm \in L^2 (\R,d\Sigma_\pm)$
the following inequalities hold:
\begin{gather} \label{e Eneq g-}
\int\limits_{\eta \in \R}
\frac{\im M_\pm (\eta+i\ep)}{|M_+ (\eta+i\ep) - M_- (\eta+i\ep)|^2}
\left|
\int_\R \frac{g^- (t)d\Sigma_- (t)}{t-(\eta+i\ep) }
\right|^2 d\eta \leq K^- \| g^- \|_{L^2(d\Sigma_-)}^2 , \\
\label{e Eneq g+}
\int\limits_{\eta \in \R}
\frac{\im M_\pm (\eta+i\ep)}{|M_+ (\eta+i\ep) - M_- (\eta+i\ep)|^2}
\left|
\int_\R \frac{g^+ (t)d\Sigma_+ (t)}{t-(\eta+i\ep) }
\right|^2 d\eta \leq K^+ \| g^+ \|_{L^2(d\Sigma_+)}^2 ,
\end{gather}
where $K^\pm$ are constants independent of $\ep$ and $g^\pm$.
\end{theorem}

\begin{proof}

It is known (see  \cite{RS2}) that for any  selfadjoint  $B=B^*$  with resolution of
identity $E_t^B $ the following identity holds
      \begin{align} \label{e EstSelfadj}
\ep\cdot \int\limits_{\eta\in\R} \| \Rs_B(\eta+i\ep ) f \|^2 \ d\eta = \pi \|f \|^2,
\qquad \ep > 0, \quad  f\in \Hsp.
         \end{align}
It follows from \eqref{e RA} that
\[
\| \Rs_A (\la ) f \|^2  -  2\| \Rs_{A_0^-\oplus A_0^+} f \|^2 \leq 2 \|
G^-(\la)\psi_-(\la) + G^+(\la)\psi_+(\la) \|^2
\]
\[
 \leq
4\| \Rs_A (\la ) f \|^2 + 4\| \Rs_{A_0^-\oplus A_0^+} f \|^2 .
\]
On the other hand, it follows from  \eqref{e EstSelfadj} with    $B$     replaced     by
$ A_0^-\oplus A_0^+$ that
        \begin{align} \label{e EstRes1}
& \frac {\ep}{2} \intl_{\eta\in\R} \| \Rs_A (\eta+i\ep) f \|^2 d\eta - \pi\|f \|^2 \leq
\notag \\ & \leq \ep \intl_{\eta\in\R}\| G^-(\eta+i\ep)\psi_-(\eta+i\ep) +
G^+(\eta+i\ep)\psi_+(\eta+i\ep) \|^2 d\eta \leq \notag \\ &\leq 2\ep \intl_{\eta\in\R} \|
\Rs_A (\eta+i\ep) f \|^2 d\eta + 2\pi\|f \|^2.
            \end{align}
Since $\psi_\pm \in L^2 ( \R_\pm ,dx)$ and $\| \psi_\pm (\cdot, \la ) \|_{L^2 (\R_\pm)}^2
= \im M_\pm (\la)/\im \la $ (see \cite{LevSar}), we have
         \begin{align*} 
&\| G^-(\la)\psi_-(\cdot, \la) + G^+(\la)\psi_+(\cdot, \la) \|^2 =\\ &= | G^-(\la)|^2 \|
\psi_-(\cdot, \la) \|^2 + |G^+(\la)|^2
\| \psi_+(\cdot, \la) \|^2 = \\
&= \frac{1}{|M_+ (\la )- M_- (\la )|^2}
\left| \,
\intl_\R \frac{ g^- (t)d\Sigma_- (t)-g^+ (t)d\Sigma_+ (t)}
{t-\la }
\right|^2 \frac{\im M_+ (\la)+\im M_- (\la)}{\im \la } .
\end{align*}
Combining this relation with \eqref{e EstRes1}  
one concludes  that \eqref{e ine1J} is equivalent to the following condition
\begin{align} \label{e Eneq gpm}
& \int\limits_{\eta \in \R}
\frac{\im M_+ (\eta+i\ep) + \im M_- (\eta+i\ep)}
{|M_+ (\eta+i\ep) - M_- (\eta+i\ep)|^2}
\left| \int_R
\frac{g^- (t)d\Sigma_- (t) -g^+ (t)d\Sigma_+ (t)}{t-(\eta+i\ep) }
\right|^2 d\eta
\leq 
C_1 \| f \|^2 ,
\end{align}
where $C_1$ is a constant independent of $f$ and $\ep$.

By definition, $\| g^\pm \|_{L^2(d\Sigma_\pm)} = \| f_\pm \|_{L^2(\R_\pm)}$, where
$f_\pm= P_{\pm}f.$ 
Thus, condition \eqref{e Eneq gpm} holds iff both \eqref{e Eneq g+} and \eqref{e Eneq g-}
are satisfied.
         \end{proof}

\subsection{Necessary conditions of  similarity in terms of the Weyl functions and Hilbert
transforms.}

Let $\Sigma_\pm = \Sigma_{ac\pm} + \Sigma_{s\pm} = \Sigma_{ac\pm} + \Sigma_{sc\pm} +
\Sigma_{d\pm} $ be the Lebesgue decomposition of the measure $\Sigma_\pm$ into a sum of
absolutely continuous, singular continuous, and pure point measures (see, for example,
\cite{RS1}).

Denote by $S'_{ac}(\Sigma_\pm)$ and $ S'_s(\Sigma_\pm)$
mutually disjoint (not necessarily topological) supports of measures $\Sigma_{ac\pm}$ and
$\Sigma_{s\pm}$, respectively.

Note that for almost all $t\in \R$ the nontangential limit
\[
\lim\limits_{\substack{\la \ra t \\ \sphericalangle}} M_\pm (\la) =: M_{\pm}(t)
\]
exists (see \cite{Gar}). Since $M_+ (\la) \not \equiv M_{-} (\la)$ on $\C_+$, we see, by
the Luzin-Privalov uniqueness theorem  (see e.g. \cite{Koosis}), that
\begin{equation} \label{e MneM}
M_+ (\eta ) \neq M_- (\eta) \quad \text{ a.e.  on } \ \R \ .
\end{equation}

\begin{theorem}  \label{th MNessH}
Let the operator $A$ be similar to a selfadjoint operator. Then, the following
inequalities hold
         \begin{gather}\label{4.1Hilbert+}
\int\limits_\R \frac{\im M_{\pm}(t)}{|M_+ (t) - M_- (t)|^2} \left|
  g^+(t)\Sigma_{ac+}^\prime (t) + (H (g^+\cdot d\Sigma_+)(t) \right|^2 dt
\leq K_1^+ \int_\R | g^+ (t)|^2 d{\Sigma_{+}(t)} ,
\\  \label{4.1Hilbert-}
\int\limits_\R \frac{\im M_{\pm}(t)}{|M_+ (t) - M_- (t)|^2}
\left|
  g^-(t)\Sigma_{ac-}^\prime (t) + (H (g^-\cdot d\Sigma_-)(t) \right|^2 dt
\leq K_1^- \int_\R | g^- (t)|^2 d{\Sigma_{-}(t)} ,
\end{gather}
with constants  $K_1^+$ and $K_1^-$ independent of $g^\pm \in L^2 \left(\R , d\Sigma_{\pm}\right).$
       \end{theorem}
                \begin{proof}
Applying  Fatou's theorem and using   \eqref{2.7.2} we get
          \begin{equation}\label{4.12}
\lim_{\ep\downarrow 0}\int\frac{g^{\pm}(t)}{t-(\eta + i\ep) }d\Sigma_{\pm}(t) = \pi \cdot
[g^{\pm}(\eta)\Sigma'_{\pm}(\eta) + i H(g^{\pm}d\Sigma_{\pm})(\eta)]
  \end{equation}
Passing to the limit in \eqref{e Eneq g-} (resp.,  \eqref{e Eneq g+}) as $\ep \to 0$ and
taking \eqref{4.12} into account we arrive   at the  inequality     \eqref{4.1Hilbert+}
(resp., \eqref{4.1Hilbert-}).
                   \end{proof}
       \begin{corollary}  \label{th MNess}
Let the operator $A$ be similar to a selfadjoint operator. Then
        \begin{gather}\label{4.10}
\frac{ \im M_{\pm}(t)}{M_+ (t)-M_-(t)}  \in L^{\infty } (\R) \ . \qquad
     \end{gather}
          \end{corollary}
        \begin{proof}
Let $A$ be similar to a selfadjoint operator. Then inequalities \eqref{4.1Hilbert+} and
\eqref{4.1Hilbert-} hold.  By Fatou Theorem  $\pi \Sigma_{ac\pm}^\prime (t) =  \im M_\pm
(t + i0) =: \im M_\pm (t)$ \  for a.e. $t\in \R.$
Taking this  relation  into account and   substituting in    
\eqref{4.1Hilbert+}  (resp. \eqref{4.1Hilbert-})  any   real-valued $g^+_{ac}$ (resp.
$g^-_{ac}$)   with $g^{\pm}_{ac}(t)=0 $ for $ t \in S'_s(\Sigma_{\pm}),$  we easily get
       \begin{gather*} 
\int\limits_{\eta \in \R} \frac{ (\im M_\pm (\eta))^2} {|M_+ (\eta) - M_- (\eta)|^2}
|g_{ac}^{\pm} (\eta)|^2 \cdot \Sigma_{ac\pm}^\prime (\eta)   d\eta \leq K_1^-\int_\R
|g_{ac}^{\pm} (t)|^2 \cdot \Sigma_{ac\pm}^\prime (t) dt.
      \end{gather*}
Since this  inequality holds for any $g_{ac}^{\pm} \in L^2 (\R,d\Sigma_{ac\pm})$, we have
\begin{gather} \label{e CF in+}
       \frac{( \im M_\pm (t))^2}{|M_+ (t)-M_-(t)|^2} \in L^\infty (S'_{ac}(\Sigma_{\pm})).
        \end{gather}
Inequality  \eqref{e CF in+}  yields  \eqref{4.10} since   $\im M_\pm (t)=0$ for a.e.
$t\in \R\setminus  S'_{ac}(\Sigma_{\pm})$.
           \end{proof}
         \begin{corollary}\label{cor5.2}
Let the operator $A$ be similar to a selfadjoint operator. Then, for all
\[ h^\pm \in
L^2(\R) \cap L^2 \left(\frac 1{\Sigma_{ac\pm}^\prime (t)}, \R \right) , \quad h^\pm (t) =
0 \quad \text{for} \quad t  \in  S'_s(\Sigma_{\pm}) \ ,
\]
the following inequalities hold:
       \begin{gather}\label{weakHilbert+}
\int\limits_\R \frac{\im M_{\pm}(t)}{|M_+ (t) - M_- (t)|^2} \left| ( H h^+) (t) \right|^2
dt \leq K_1^+ \int_\R | h^+ (t)|^2 \frac 1{\im M_{+}(t)} dt ,
\\ \label{weakHilbert-}
\int\limits_\R \frac{\im M_{\pm}(t)}{|M_+ (t) - M_- (t)|^2} \left| ( H h^-) (t) \right|^2
dt \leq K_1^+ \int_\R | h^- (t)|^2 \frac 1{\im M_{-}(t)} dt ,
\end{gather}
where $K_1^+$ and $K_1^-$ are constants independent of $h^\pm$.
           \end{corollary}
\begin{proof}
Inequality  \eqref{4.1Hilbert+} yields
         \begin{equation} \label{4.1AHilbert+}
\int\limits_\R \frac{\im M_{\pm}(t)}{|M_+ (t) - M_- (t)|^2} \left|
 (H (g^+\cdot d\Sigma_+)(t) \right|^2 dt
\leq K_1^- \int_\R | g^+ (t)|^2 d{\Sigma_{+}(t)} ,
          \end{equation}
Choosing any    $g^+_{ac}$  with $g^{+}_{ac}(t)=0 $ for $ t \in S'_s(\Sigma_{+}),$ and
setting in  \eqref{4.1AHilbert+}  $h^\pm  := g^\pm \cdot (\Sigma_{ac\pm}^\prime )$ we
arrive at the inequality  \eqref{weakHilbert+}. The inequality  \eqref{weakHilbert-} is
implied by \eqref{4.1Hilbert-} in just the same way.
      \end{proof}

               \begin{corollary}
Let $E_{\pm}=\supp\Sigma_{ac\pm}^{\, \prime}$ be the topological supports of measures
$\Sigma_{ac\pm}$. If the operator $A$ is similar to a selfadjoint operator, then
       \begin{equation} \label{Mac2}
\sup_{\I}\left(\frac{1}{|\I\cap E_{\pm}|}\int_{\I}\frac{\im
M_{\pm}(t)}{|M_+(t)-M_-(t)|^2}dt\right)\cdot \left(\frac{1}{|\I\cap E_{\pm}|}\int_{\I}
\im M_{\pm}(t)dt\right)<\infty.
    \end{equation}
         \end{corollary}
       \begin{proof}
If $A$  is similar to a selfadjoint operator, then by Corollary  \ref{cor5.2} two-weight
estimates  \eqref{weakHilbert+} and  \eqref{weakHilbert-}  for the Hilbert transform are
valid. Due to \eqref{twoA2} the  result  is immediate  from   \eqref{weakHilbert+} and
\eqref{weakHilbert-}.
        \end{proof}

Due to Lebesgue theorem  inequality \eqref{Mac2} yields  \eqref{4.10} and therefore gives
another proof of Corollary  \ref{th MNess}. In fact, it gives a new necessary condition
of similarity to a selfadjoint operator and  is stronger than \eqref{4.10}.

The following corollary gives one more necessary condition of similarity.
        \begin{corollary}
Let  $A$ be  similar to a selfadjoint operator and  let
$$
w_{1\pm}(t) := \frac{\im M_{\pm}(t)}{|M_+(t) - M_-(t)|^2}.
$$
Then
        \begin{equation}\label{PoissonCond}
\sup_{\la\in{\C}_+} {\mathcal P}_{\la}(w_{1\pm}) \cdot \im M_{ac\pm}(\la) = C<\infty,
      \end{equation}
where $ M_{ac\pm} (\lambda ) := \displaystyle \int_{\R} \frac{d \Sigma_{ac\pm} (t)}{t-\lambda} $,  $\lambda \in \C_+$.

           \end{corollary}
\begin{proof}  Note  that  $\im M_{\pm}(t)$ is finite for a.e. $t\in \R$ and
${\mathcal P}_{\la}(\im M_{\pm}) = \im M_{ac\pm}(\la).$ We complete the proof by combining
Corollary  \ref{cor5.2}  with Proposition \ref{prop2.5B}.
          \end{proof}
Inequality \eqref{2.40C}  shows that condition \eqref{PoissonCond} is stronger than
\eqref{Mac2}.

%
%

\begin{conjecture}
 We conjecture  that under the condition $\sigma_{disc} (A) = \emptyset $ inequalities  \eqref{4.1Hilbert+} and
\eqref{4.1Hilbert-} are also sufficient for  the operator $A$ to be similar  to a
selfadjoint operator. Therefore  inequalities   \eqref{4.1Hilbert+} and
\eqref{4.1Hilbert-}  reduce the similarity problem to two weight estimates for the
Hilbert transform.
\end{conjecture}

\begin{conjecture}   Suppose that $\sigma_{disc} (A) = \emptyset $ and both measures $d\Sigma_+$ and $d\Sigma_-$ are absolutely
continuous, $\Sigma_{\pm}=\Sigma_{ac\pm }$.  Then conditions  \eqref{weakHilbert+} and
\eqref{weakHilbert-} are sufficient for $A$ to be similar to a selfadjoint operator.
\end{conjecture}

\subsection{Sufficient conditions of similarity  in terms of  Weyl functions.}

Consider an operator $\wt A$ given by
$\ \wt A = A_{min}^*  \upharpoonright \! \dom (\wt A) $,
\begin{equation} \label{e dom tilA}
\dom( \wt A) =\{ y \in \dom (A_{min}^*) : \
y(+ 0)=y(-0), \ y^\prime (+ 0) = - y^\prime (-0) \} \ .
\end{equation}

      \begin{proposition} \label{p ResAtil}
The operator $\wt A$ is selfadjoint.
For $\la \in \C\setminus\R$ the resolvent of $
\wt A$ has the form
\begin{equation} \label{e RAtil}
\Rs_{\wt A}(\la ) f= \Rs_{A_0^-\oplus A_0^+} (\la ) f +
\wt G^-(\la)\psi_-(\la) + \wt G^+(\la)\psi_+(\la) ,
\end{equation}
\begin{equation} \label{e Gtil}
\wt G^+ (\la)=-\wt G^- (\la ) = \frac{1}{M_+ (\la) + M_- (\la)}
\int_\R \frac{g^-(t)d\Sigma_-(t)+g^+(t)d\Sigma_+(t)}{t-\la} ,
\end{equation}
where $g^\pm (t)=(\F_\pm f_\pm) (t)$, $f_\pm := P_\pm f \in L^2(\R_\pm)$.
       \end{proposition}

\begin{proof}
Let  $\Pi=\{{\C}^2,  \Gamma_0, \Gamma_1\}$ be a boundary triplet for   $S^* :=
A_{min}^* $ defined by  \eqref{boundtriple}. Clearly, the extension  $\wt A $ of $A_{min}$
determined by \eqref{e dom tilA}, admits the following representation
    \begin{equation}
\wt A= S^*|\dom \wt A, \qquad  \dom \wt A=\ker(\Gamma_1-B\Gamma_0), \qquad \text{where}
\qquad B=
\begin{pmatrix}
0&-1\\
-1&0
\end{pmatrix}.
    \end{equation}
Thus, $\wt A$  is selfadjoint because so is  $B.$

The representation  \eqref{e RAtil}  for the resolvent $\Rs_{\wt A} (\la )$ can be
obtained in just the same way as representation for  $\Rs_{A} (\la ) $ in
Proposition~\ref{p Res}.
       \end{proof}
         \begin{theorem} \label{th MpmCond}
Suppose that
         \begin{equation}\label{sufficientcond}
\sup_{\la\in \C_+} \dfrac{|M_+ (\la ) + M_- (\la)|}{| M_+ (\la ) - M_- (\la)|} \ < \
\infty.
        \end{equation}
Then the operator $A$ is similar to a selfadjoint operator.
            \end{theorem}
\begin{proof} 
Since $\wt A$ and $ A_0^-\oplus A_0^+$ are selfadjoint operators, we obtain    from
\eqref{e RAtil}  and  \eqref{e EstSelfadj}

       \begin{equation} \label{4.19A}
\ep \intl_{\eta \in \R} \| \wt G^-(\eta+i\ep) \ \psi_-(\eta+i\ep) + \wt G^+(\eta+i\ep)\
\psi_+(\eta+i\ep) \|^2 d\eta \leq 4\pi\| f\|^2 .
         \end{equation}
On the other hand,  it follows from   \eqref{e Gtil}  with $f=f_{\pm}$   that
\begin{gather} \label{4.20A}
\|  \wt G^{\pm}(\eta+i\ep)\psi_{\pm}(\eta+i\ep) \|^2 = \\
= \frac{\im M_+ (\la ) + \im M_- (\la )}{\im \la \cdot | M_+ (\la ) + M_- (\la )|^2}
\left| \int_\R \frac{g^{\pm}(t)d\Sigma_{\pm}(t)}{t-\la} \right|^2 .
    \end{gather}

Combining   \eqref{4.19A}  with    \eqref{4.20A} we arrive at the following inequalities
         \begin{gather*} 
\intl_{\eta \in \R} \frac{\im M_{\pm} + \im M_{\mp} (\la ) }{| M_+ (\la ) + M_- (\la
)|^2} \left| \int_\R \frac{g^{\pm}(t)d\Sigma_{\pm}(t)}{t-\la} \right|^2 d\eta \leq 4\pi\|
f_{\pm} \|^2 = 4\pi\| g^{\pm} \|^2 .
           \end{gather*}

Combining these inequalities with    \eqref{sufficientcond}
we arrive at estimates \eqref{e Eneq g-} and \eqref{e Eneq g+}. Thus, by Theorem \ref{th
CrSLSim}, $A$ is similar to a selfadjoint operator.
         \end{proof}

\begin{remark}
The condition \eqref{sufficientcond} is not necessary for similarity to selfadjoint operator (see Remark \ref{r M+M})). 
\end{remark} 

      \begin{remark}
Note that sufficient condition  \eqref{sufficientcond} for similarity is  weaker than
either conditions  \eqref{2.55A}--\eqref{2.56A} or conditions   \eqref{2.58A} obtained
from Theorem \ref{t MMM00}  and Theorem \ref{t Sah}, respectively.  While these conditions
guarantee a stronger result: similarity of $A$ to an operator $B=B^*$ with absolutely
continuous spectrum.
       \end{remark}



Finally, we apply Theorems ~\ref{th MNess} and \ref{th MpmCond} to the case of the
operator $A$ with constant potential. Consider a family of such operators
    \begin{equation}\label{SLconst}
A(a):=(\sgn x) (-d^2/dx^2 +a),  \qquad  a \in \R,
        \end{equation}
depending on a parameter $a.$
        \begin{proposition}[\cite{KarMFAT00},\cite{KarMN00}] \label{c ConstPot}
(i)\  The operator  $A(a)$ is similar to a selfadjoint operator if and only if $a\geq 0$.

(ii) The operator $A(0)$ is similar to the multiplication operator $Q:\  f\to xf(\cdot)$
in  $L^2(\R).$
       \end{proposition}
        \begin{proof}
(i) In the case under consideration the functions $M_\pm (\la)$ are given by
      \begin{equation} \label{e Mpm const}
M_\pm (\la ) = \pm \frac{i}{\sqrt{\pm\la-a}}\  .
         \end{equation}
Since
\[
M_+ (\la ) -M_- (\la ) = \frac {i}{(\la-a)^{1/2}}+ \frac {i}{(-\la-a)^{1/2}} \neq 0
\qquad \text{for}\quad  \la \notin \R,
\]
Proposition \ref{p Res} yields that the spectrum of $A(a)$ is real for any $a \in \R$ (see
also \cite{CurNaj96}). It is clear that $M_+ $ and $M_-$ are holomorphic on $\C\setminus
[a,+\infty)$ and $\C\setminus (-\infty,-a]$, respectively. Hence, by Proposition, we have 
\ref{prop2.4}  (iv) $\sigma (A(a))= (-\infty,-a] \cup [a,+\infty)$, that is  $\sigma
(A(a))= \R$ for $a\le 0$ and $\sigma (A(a))=\R\setminus (-a,a)$ for  $a>0$.

If $a\geq0$, then the function
$$
\frac{M_+(\la ) + M_- (\la )}{ M_+(\la ) - M_- (\la )}
$$
is bounded in $\C_+$. Thus, by Theorem~\ref{th MpmCond}, $A$ is similar to a selfadjoint
operator.

Now let $a<0$.  Setting $\la=i\varepsilon$ and   $i\varepsilon - a= \rho e^{i\phi}$ we
get
$$
 M_+ (i\varepsilon)-M_-(i\varepsilon)  = i\rho^{-1/2}\cdot [e^{-i\phi/2} -  e^{i\phi/2}]
 = 2\rho^{-1/2}\sin(\phi/2),
$$
and
$$
\im M_+(i\varepsilon) = \im (i\rho^{-1/2}e^{i\phi/2}) =  \rho^{-1/2}\cos(\phi/2).
$$
Hence
$$
{\im M_{+}(i\varepsilon)}(M_+ (i\varepsilon)-M_-(i\varepsilon))^{-1} = 2^{-1}\cot(\phi/2)
$$
is  unbounded in any neighborhood of zero. Thus, by Corollary ~\ref{th MNess} the
operator $A$ is not similar to a selfadjoint operator.

(ii)\  Let now $A=A(0).$ Substituting expressions  \eqref{e Mpm const} in formula
\eqref{formulaforcharfunc} for $\theta_A(\cdot)$ and using  the relation $\sqrt{\la}/
\sqrt{-\la}= -i$,  we arrive at  the following formula for the  characteristic function
       \begin{equation}\label{formulaforcharfunc2}
\theta_A(\la) =
\begin{pmatrix}
-i& (i-1)/\sqrt{-\la}\\
(i-1)\sqrt{\la}&-i
\end{pmatrix}\  .
    \end{equation}
It follows that $\theta_A(\cdot)$ is unbounded only near zero and infinity. Since the
operator $A$ has no eigenvalues, then by Proposition  \ref{prop2.5abscont}  (or by
Corollary   \ref{cor2.5abscont}) it is similar to a selfadjoint operator  $T_0 = T_0^*$
with absolutely continuous spectrum, $\sigma (T_0)= \sigma_{ac} (T_0) =\R, \  \sigma_{s}
(T_0) = \sigma_{p} (T_0)= \emptyset$. It is easily seen that the multiplicity of spectrum
is one. Therefore  $T_0$ is unitarily  equivalent  to the multiplication operator $Q.$
         \end{proof}
\begin{remark}
Using the Krein-Langer spectral theory of definitizable operators in Krein spaces \'Curgus
and Langer \cite{CL89} investigated  the critical point $\infty$ of differential operators with an indefinite weight. Their results  imply similarity of the operator $A(a)$  to a
selfadjoint one if only  $a>0$.

The case $a=0$ is more complicated since  $A(0)$ has two critical points:
zero and infinity. Similarity of  $A(0)$ to a selfadjoint operator  was established  by
\'Curgus and Najman \cite{CN95} in the framework of  Krein space approach.

Other proofs of the latter result have been obtained by several authors  (see
\cite{KarKROMSH98,KarMFAT00,FSh00,Kap}).
In full generality statement (i) of Proposition \ref{c ConstPot} has originally been
proved by one of the authors  \cite{KarMN00,KarMFAT00}, by using  the resolvent
criterion of similarity (see Theorem \ref{t SimCr}). The proof given above is similar to
that contained in our short communication \cite{KarMal04}.
        \end{remark}

\section{Restrictions of $A$ to invariant subspaces corresponding
to $\sigma_{disc} (A)$ and $\sigma_{ess} (A)$}

Throughout this section we assume additionally that the following assumption is valid.

\begin{assumption} \label{a sdisc}
Suppose that  the set \quad $\sp_{disc} (A)$ \quad is finite.
\end{assumption}


It will be shown in Section \ref{ch SL ss FZ} that this condition is fulfilled if the
potential $q$ is finite-zone.

Since $\dist (\sp_{ess} (A), \sp_{disc} (A) ) >0$, we can apply
the theorem on spectral decomposition (see \cite[Theorem III.6.17]{Kato}).
That is there exists a skew
decomposition $ L^2 (\R ) = {\Hsp} = \Hsp_e \dot + \Hsp_d $ such that
\begin{gather} \label{e skew}
A = \Aess \dot + \Adisc, \quad \Aess = A \upharpoonright ( \dom (A) \cap \Hsp_e),
\quad \Adisc = A \upharpoonright ( \dom (A) \cap \Hsp_d) \\
\text{and} \quad
\sp (\Adisc) = \sp_{disc} (A), \quad \sp (\Aess) = \sp_{ess} (A). 
\notag
\end{gather}
We denote by $P_e$ and $P_d$ the corresponding skew
projections onto $\Hsp_e$ and $\Hsp_d$, respectively.

Since $ \sp_{disc} (A)$ is finite, we see that $\Adisc$ is an operator
in a finite dimensional space $\Hsp_d$. Jordan normal form of
$\Adisc$ is described in Proposition \ref{p e and d} (3)-(5).
By Proposition \ref{p e and d}, we have
$\sp_{ess} (A) = \sp_{ess} (A_2^-) \cup \sp_{ess} (A_2^+)$.
Thus $ \sigma( \Aess) \subset \R$.
This section is devoted to
the question of similarity of $\Aess$ to a selfadjoint operator.

\begin{proposition} \label{p CrSimAe}
Let Assumption \ref{a sdisc}    
be fulfilled.
Suppose $G_d$ be a compact subset of  $\C$ such that $G_d \cap \sp_{ess} (A) = \emptyset$
and all $\la \in \sp_{disc} (A) $ are interior points of $G_d$. Suppose $Q^\pm$ are dense
subsets in $L^2 (\R , d\Sigma_\pm)$.

Then the following conditions are equivalent:
\begin{description}
\item[(i)] the operator $\Aess$ in $\Hsp_e$
is similar to a selfadjoint one;
\item[(ii)] the operator $\Aess P_e$ in $\Hsp = L^2 (\R)$ is similar
to a selfadjoint one;
\item[(iii)] the inequality
\begin{gather} \label{e EneqEss}
\ep \cdot \int_{\R } \| \Rs_A (\eta+i\ep) f_e \|^2 d\eta \leq C_1^e \| f_e \|^2 \
\end{gather}
holds for all $\ep>0$, $f_e \in \Hsp_e$ with
some constant $C_1^e$.
\item[(iv)] for all $\ep>0$ and $g^\pm \in Q^\pm$
the following inequalities hold:
         \begin{gather}\label{e EneqEss2-}
\int\limits_{\substack{\eta \in \R \\ \eta+i\ep\not \in G_d }} \frac{\im M_\pm
(\eta+i\ep)}{|M_+ (\eta+i\ep) - M_- (\eta+i\ep)|^2} \left| \int_\R \frac{g^- (t)d\Sigma_-
(t) } {t-(\eta+i\ep) } \right|^2 d\eta
\leq 
C_2^- \| g^- \|^2_{L^2 (d\Sigma_-)} ,
\end{gather}
\begin{gather}
\int\limits_{\substack{\eta \in \R \\ \eta+i\ep\not \in G_d }}
\frac{\im M_\pm (\eta+i\ep)}{|M_+ (\eta+i\ep) - M_- (\eta+i\ep)|^2}
\left|
\int_\R \frac{g^+ (t)d\Sigma_+ (t) } {t-(\eta+i\ep) }
\right|^2 d\eta
\leq 
C_2^+ \| g^+ \|^2_{L^2 (d\Sigma_+)} ,
          \label{e EneqEss2+}
\end{gather}
where $C_2^\pm$ are constants independent of $\ep$ and $g^\pm$.
\end{description}
\end{proposition}

\begin{proof}
It is clear that $ (i) \Leftrightarrow (ii)$.

Let us show that $ (ii) \Leftrightarrow (iii)$.
It can easily be checked that $\Aess P_e$ is a $J$-selfadjoint operator
(see \cite{Lan82}).
By Proposition~\ref{p SimCrJ}, assertion (ii) holds if and only if
for all $\ep>0$ and $f \in \Hsp$
the following inequality holds
\begin{gather} \label{e EneqPeAPe}
\ep \intl_{\R} \| \Rs_{A P_e} (\eta+i\ep) f \|^2 d\eta \leq C_1 \| f \|^2 \ , \quad
C_1 = const.
\end{gather}
Clearly, \eqref{e EneqPeAPe} is equivalent to \eqref{e EneqEss}.

Now we show that $ (iii) \Leftrightarrow (iv)$.
Let $f \in L^2 (\R)$, $\ f_e = P_e f$, $ \ f_d = P_d f$.
It can be shown in the usual way that there exist constants
$C_2$, $C_3$ such that
$\| \Rs_A f_e \| = \| \Rs_{\Aess} f_e \| \leq C_2 \| f_e \|$ for
$\la \in G_d$, and
$ \| \Rs_A f_d \| = \| \Rs_{\Adisc} f_d \|
\leq \dfrac {C_3}{1+ |\la|} \| f_d \|$
for $\la \in \C\setminus G_d$.
Therefore \eqref{e EneqPeAPe} is equivalent to
\begin{gather} \label{e AeGdRsCr}
\ep \! \! \intl_{\substack{\eta \in \R \\ \eta+i\ep\not \in G_d }}
\| \Rs_A (\eta+i\ep) f \|^2 d\eta \leq C_1 \| f \|^2 , \quad
\forall f \in L^2 (\R) , \forall \ep>0.
\end{gather}
Arguing as in the proof of the Theorem~\ref{th CrSLSim}, we see that condition \eqref{e
AeGdRsCr} is fulfilled iff  the  inequalities \eqref{e EneqEss2-} and \eqref{e EneqEss2+}
hold for all $g^\pm \in L^2 (d\Sigma_\pm)$  and  $\ep >0$.

We show that it suffices to check   \eqref{e EneqEss2-} and \eqref{e EneqEss2+}
only for dense subsets $Q^\pm$.

Let $\ep >0$ be a fixed positive number, $ \I $
an open bounded set in $\R$. Denote \\
$ \I_\ep := \{ \eta + i\ep : \, \eta \in \I \} $. Assume that $ \I_\ep \cap G_d =
\emptyset$. Then $(M_+ (\la) -M_- (\la ) )^{-1}$ is  holomorphic  on $\I_\ep$. By the
Schwarz inequality, the operators
\[
K_{\I_\ep}^\pm: g^+ \mapsto
\frac{(\im M_\pm (\eta+i\ep))^{1/2}}{M_+ (\eta+i\ep) - M_- (\eta+i\ep)}
\int_\R \frac{g^+ (t) d\Sigma_+ (t) } {t-(\eta+i\ep) } , \quad
\]
are bounded  from $L^2 (\R, d\Sigma_+)$ to $L^2 (\I_\ep, d\eta).$

Suppose that  $Q^+$ is dense in $L^2 (d\Sigma_+)$ and \eqref{e EneqEss2+} is fulfilled
for $Q^+$. Then $\| K_{\I_\ep}^\pm \| \leq C_2^+$ for all $\ep >0$ and for all $\I$. This
imply  \eqref{e EneqEss2+}   for all $g^\pm \in L^2 (d\Sigma_\pm)$  and  $\ep >0$. In the
same way we can prove that \eqref{e EneqEss2-} is equivalent to the inequality \eqref{e
EneqEss2-}  for all $g^\pm \in L^2 (d\Sigma_\pm)$.
          \end{proof}


Recall that $\sp_{ac} (T)$ and $\sp_s (T)$ are
the absolutely continuous and singular spectra of a selfadjoint
operator $T$.
Evidently,
\[ \sp_{ac} (A_2^\pm) = \supp d\Si_{ac\pm}, \qquad
\sp_s (A_2^\pm) = \supp (d\Si_{sc\pm} +d\Si_{d\pm} ).
\]
Note that $\sp_{ac} (A_0^\pm) \subset \sp_{ess} (A_0^\pm) $.
Therefore, by Proposition \ref{p e and d}, we have
\begin{gather} \label{e Sac in Sess}
\supp d\Si_{ac-} \cup \supp d\Si_{ac+} = \sp_{ac} (A_0^- \oplus A_0^+)
\subset \sp_{ess} (A) .
\end{gather}

\begin{proposition} \label{p NessAess}
Let Assumption \ref{a sdisc}    
be fulfilled. Suppose the
operator $\Aess$ is similar to a selfadjoint operator. Then
\begin{gather}
\quad \frac{\im M_{ac\pm}(t)}{M_+ (t)-M_-(t)} \in L^{\infty} (\R). \qquad \quad
\end{gather}
\end{proposition}

Taking into account \eqref{e Sac in Sess}, we see that
this theorem can be proved in the same way as Theorem~\ref{th MNess}.

\begin{assumption} \label{a s2}
In what follows we assume that
\begin{gather*} 
d\Si_- = d\Si_{ac-} + d\Si_{d-}, \quad
\supp d\Si_{d-} = \{ \theta_j^- \}_{j=1}^{N_\theta^-},
\quad N_\theta^- <\infty, \\
\text{and} \qquad d\Si_+ = d\Si_{ac+} + d\Si_{d+}, \quad
\supp d\Si_{d+} = \{ \theta_j^+ \}_{j=1}^{N_\theta^+},
\quad N_\theta^+ <\infty .
\end{gather*}
\end{assumption}
Then $M_\pm (\la ) = M_{ac\pm} (\la ) + M_{d\pm} (\la)$,
where
\begin{gather*} 
M_{ac\pm} (\lambda ) = \int_{\R} \frac{d \Sigma_{ac\pm} (t)}{t-\lambda} \ , \qquad
\text{and} \qquad
M_{d\pm} (\lambda ) = \suml_{j=1}^{ N_\theta^\pm}
\frac{ \Sigma_{\pm} (\theta_j^\pm+0) - \Sigma_{\pm} (\theta_j^\pm-0)}
{\theta_j^\pm -\lambda} \ .
\end{gather*}.

Let us introduce the sets
\begin{gather} \label{e wt theta}
\, \{ \wt \theta_j^\pm \}_1^{\wt N_\theta^\pm} :=
\{ \theta_j^\pm \}_1^{N_\theta^\pm} \setminus
\sigma_{disc} (A) ;
\end{gather}
here $\wt N_\theta^\pm < \infty $.
(these sets will be used in Theorem \ref{t uwSufCond}).

Recall that we denote the Smirnov class on $\C_+$
(see Subsection \ref{ss Hp}) by $\Np (\C_+)$.

\begin{theorem} \label{t uwSufCond}
Let Assumptions \ref{a sdisc} and 
\ref{a s2} be fulfilled.
Let $G_d$ be the compact set from Proposition~\ref{p CrSimAe}.

Suppose there exist functions $U_+ (\la)$ and $U_- (\la)$ on $\C_+$ such that
the following conditions hold:
\begin{gather} \label{e Cond u2}
\frac{\im M_{ac\pm} (\la)}{|M_+ (\la) - M_- (\la)|^2}
\leq C_\pm^u | U_\pm (\la)|^2 ,
\quad \la \in \C_+ \setminus G_d \ ,
\end{gather}
\begin{gather} \label{e Cond uN+}
\quad U_\pm (\la) \in \Np (C_+) ,
\end{gather}
\begin{gather}
\frac{U_\pm (t)} {\theta_j^- - t }
\in L^2(\R), \ j=1, \cdots , N_\theta^- ; \quad
\frac {U_\pm (t)}{ \theta_j^+ - t}
\in L^2(\R), \ j=1, \cdots , N_\theta^+ \ , \label{e Cond u3}
\end{gather}
where $C_\pm^u$ are constants.

Suppose there exist functions $w_+ (\cdot )$ and $w_- ( \cdot )$ on $\R$ , $w_\pm (t) > 0$ a.e.,
such that the following conditions hold:
\begin{gather} \label{e Cond w1}
w_\pm (t) \leq C_\pm^w (\Si_{ac\pm}^\prime (t))^{-1}  \quad \text{a.e. on } \
\supp d\Si_{ac\pm} ,
\end{gather}
\begin{gather} \label{e Cond w2}
w_+ (t) \ \text{ and } \ w_- (t) \quad \text{satisfy the} \ (A_2) \
\text{condition (see \eqref{e (Ap)})} ,
\end{gather}
\begin{gather} \label{e Cond w3}
\frac { U_+^2 (t)}{ w_\pm (t)}
\in L^\infty (\R), \qquad
\frac { U_-^2 (t)}{ w_\pm (t)}
\in L^\infty (\R)  \ ;
\end{gather}
where $C_\pm^w$ are constants.

Suppose that for every point
$\wt \theta_j^\pm $ of the set
$ \{ \wt \theta_k^\pm \}_1^{\wt N_\theta^\pm} $,
there exist a function
$U_j^\pm (\la) \in \Np (C_+)$
and a neighborhood $D_j^\pm$ of the point $\wt \theta_j^\pm $ such that
the following conditions hold:
\begin{gather} \label{e C theta u}
\frac{1}{|M_+ (\la) - M_- (\la)|^2}
\im \frac{1}{\wt \theta_j^\pm - \la}
\leq C_\theta^u
| U_j^\pm (\la)|^2  \quad \text{for} \ \la \in
D_j^\pm \cap \C_+ \ ,
\end{gather}
\begin{gather} \label{e C theta w}
\frac { | U_j^\pm (t) |^2}{ w_+ (t)} \in L^\infty (\R), \quad
\frac { | U_j^\pm (t) |^2}{ w_- (t)} \in L^\infty (\R) \ ,
\end{gather}
\begin{gather} \label{e C theta M}
\frac {1}{|\wt \theta_j^\pm - \la| \, |M_+ (\la) - M_- (\la)|}
\leq
C_\theta^M  \quad \text{for} \  \la \in D_j^\pm \cap \C_+ ,
\end{gather}
where $C_\theta^u$ and $C_\theta^M$ are constants.

Then $\Aess$ is similar to a selfadjoint operator.
\end{theorem}

\begin{proof}
Let us show that \eqref{e EneqEss2-} and \eqref{e EneqEss2+} hold.

Let $\la = \eta+i\ep$, $\eta = \re \la$, $\ep = \im \la$.

1) Denote
\[
\I_\pm (\ep) := \int\limits
_{\substack{\eta \in \R \\ \eta+i\ep \not \in G_d}}
\frac {\im M_{ac\pm} (\la)}{|M_+ (\la) - M_- (\la)|^2}
\left| \int_\R \frac{g^+ (s)d\Sigma_+ (s) } {s-\la }
\right|^2 d\eta , \quad g^+ \in L^2 (d\Sigma_+ (t)) \ .
\]
Let
\[
Q^+_{ac} := \{ g^+ \in L^2 (\R, d\Sigma_{ac+} (t)) : \
(g^+ \Sigma_{ac+}^\prime) \in L^2 (\R, dt) \}.
\]
Then the set $Q^+ := Q^+_{ac} \oplus L^2 (\R, d\Sigma_{d+})$
is dense in $L^2 (\R, d\Sigma_+ (t))$.

First we show that
\begin{gather} \label{e Iac}
\I_\pm (\ep) \leq C_2^+ \| g^+ \|_{L^2 (d\Sigma_+)}^2 \quad \text{for}  \ g^+ \in Q^+.
\end{gather}

Let us denote
\begin{gather*}
K_\pm^{ac} (\la ) :=  U_\pm (\la)
\int_\R \frac{g^+ (t)d\Sigma_{ac+} (t) } {t-\la }, \qquad
K_\pm^{d} (\la ) :=  U_\pm (\la)
\int_\R \frac{g^+ (t)d\Sigma_{d+} (t) } {t-\la },
\\
K_\pm (\la ) :=  K_\pm^{ac} (\la) + K_\pm^d (\la) =
U_\pm  (\la) \int_\R \frac{g^+ (t)d\Sigma_+ (t) } {t-\la } .
\end{gather*}
By \eqref{e Cond u2}, we have
\begin{gather} \label{e Iac<u}
\I_\pm (\ep ) \leq \int_\R |K_\pm ( \la) |^2 d\eta .
\end{gather}

It follows from $ U_\pm (\la) \in \Np (\C_+)$ that
$ U_\pm (\la)$ is holomorphic in $\C_+$ and has the
nontangential limit $U_\pm (\eta)$ for almost all $\eta \in \R$
(see \cite{Gar}). Since $g^+ \in Q^+$,
we have
$g^+ (t)\Sigma_{ac+}^\prime (t) \in L^2 (\R, dt)$.
\ Therefore,
\[
 \dint_\R \frac{g^+ (t)d\Sigma_{ac+} (t) } {t-\la } \in H^2 (\C_+) \
.
\]
It follows from \cite[Corollary~II.5.6]{Gar} and
\cite[Corollary~II.5.7]{Gar} that
$K_\pm^{ac} (\la) \in \Np (\C_+)$.
The functions $(\theta_j^+ - \la)^{-1}$ are outer in $\C_+$.
Therefore \cite[Corollary~II.5.6]{Gar} and Lemma \ref{e N++N+}
yield $K_\pm^{d} (\la) \in \Np (\C_+)$.
Hence $K_\pm^{ac} (\la)$, $K_\pm^d (\la)$, and $K_\pm (\la)$
belong to $\Np (\C_+)$ and have the nontangential limits
$K_\pm^{ac} (\eta)$, $K_\pm^d (\eta)$ and $K_\pm (\eta)$ for almost all
$\eta \in \R$. Note also that
\begin{equation} \label{e K=P+H}
K_\pm^{ac} (\eta) := \pi U_\pm (\eta)
\left( g^+(\eta) \Si_{ac+}^\prime (\eta) +
\H (g^+ \Si_{ac+}^\prime ) (\eta) \right) \quad \text{for a.e.}
\ \eta \in \R.
\end{equation}

Assume that the following inequality holds
\begin{gather} \label{e R<fe}
\int_\R \|K_\pm (\eta)\|^2 d\eta \leq C_2^+ \| g^+ \|_{L^2 (d\Sigma_+)}^2 .
\end{gather}
Then, by \cite[Section II.5]{Gar}, we have $K_\pm (\la) \in H^2 (\C_+)$
and for all $\ep>0$
\[
\int_\R \|K_\pm (\eta+i\ep)\|^2 d\eta \leq \|K_\pm (\la)\|_{H^2 (\C_+)}^2
= \| K_\pm (\eta) \|_{L^2 (\R)}^2 \leq C_2^+ \| g^+ \|_{L^2 (d\Sigma_+)}^2 \ .
\]
Combining this with \eqref{e Iac<u}, we see that
\eqref{e R<fe} yields
\eqref{e Iac} with a constant $C_2^+$ independent of
$g^+ \in Q^+$.

Let us prove \eqref{e R<fe}.
By \eqref{e Cond u3}, we have
\begin{gather}
\| K_\pm^d (\eta) \|_{L^2 (\R)}
\leq C_3^\pm \sum_{j=1}^{N_\theta^+}
g^+ (\theta_j^+)
\left( \Si_+ (\theta_j^+ +0) - \Si_+ (\theta_j^+ -0) \right)^{1/2}
\leq \notag \\ \leq
\text{where} \ C_3^\pm \sqrt{N_\theta^+} \ \| g^+ \|_{L^2(d\Si_+)} ,
\label{e Rd<f}
\end{gather}
\begin{gather*}
C_3^\pm
= \max \left\{ (\Si_+ (\theta_j^+ +0) - \Si_+ (\theta_j^+ -0) )^{1/2}
\left\| \frac { U_\pm (\eta) }{ \theta_j^+ -\eta } \right\|_{L^2 (\R)}
       \right\}_{j=1}^{N_\theta^+} < \infty.
\end{gather*}

It follows from \eqref{e Cond w1} that
\begin{gather} \label{e gS<gw}
 \|g^+ (t) \Si_{ac+}^\prime (t) \|_{L^2 (w_+ (t)dt)}^2 \leq
C_+^w \| g^+ \|_{L^2(d\Si_{ac+}) }^2 \leq
C_+^w \| g^+ \|_{L^2(d\Si_+) }^2  \ .
\end{gather}
Since $w_+ (t) \in (A_2)$,
we have
\begin{gather} \label{e Hg<g}
\| \H (g^+ \Si_{ac+}^\prime ) (t) \|_{L^2 (w_+ (t)dt)}^2
\leq C_1
\| g^+ (t) \Si_{ac+}^\prime (t) \|_{L^2 (w_+ (t)dt)}^2 ,
\end{gather}
where $C_1$ is a constant independent of $g^+$.
It follows from \eqref{e Hg<g} and \eqref{e K=P+H} that
\begin{gather}
\int_\R |K_\pm^{ac} (\eta) |^2 d\eta
\leq \notag \\ \leq
\left\| \frac{ U_\pm^2 (\eta)}{w_+ (\eta)}
\right\|_{L^\infty (\R)} \
\int_\R |g^+(\eta) \Si_{ac+}^\prime (\eta) +
\H (g^+ \Si_{ac+}^\prime ) (\eta) |^2 w_+ (\eta) d\eta
\leq \notag \\ \leq
2(1+C_1)
\left\| \frac{ U_\pm^2 (\eta)}{w_+ (\eta)} \right\|_{L^\infty (\R)}
\| g^+ (\eta) \Si_{ac+}^\prime (\eta) \|_{L^2 (w_+(\eta) d\eta)}^2 \ .
\label{e K<U/wg}
\end{gather}
Combining \eqref{e K<U/wg}, \eqref{e Cond w3}, and \eqref{e gS<gw}, we get
\begin{gather*}
\int_\R |K_\pm^{ac} (\eta) |^2 d\eta \leq C_2 \|g^+ \|_{L^2 (d\Si_+)}^2
\ , \label{e Rac<f}
\end{gather*}
where the constant $C_2$ is independent of $g^+$.
Taking into account \eqref{e Rd<f}, we obtain \eqref{e R<fe}.

Let us remember that \eqref{e R<fe} implies \eqref{e Iac}.
Thus \eqref{e Iac} is proved.

2) Denote
\[
\I_{d\pm} (\ep) := \int\limits
_{\substack{\eta \in \R \\ \eta+i\ep \not \in G_d}}
\frac {\im M_{d\pm} (\la)}{|M_+ (\la) - M_- (\la)|^2}
\left| \int_\R \frac{g^+ (s)d\Sigma_+ (s) } {s-\la }
\right|^2 d\eta .
\]
Let us show that
\begin{gather} \label{e Id}
\I_{d\pm} (\ep) \leq C_3 \| g^+ \|_{L^2 (d\Sigma_+)}^2  \quad
\text{for} \  g^+ \in Q^+
\end{gather}
(here and below $C_3$, $C_4$, ...  are some constants).
It is suffices to prove the inequality \eqref{e Id}
for each summand,
i.e.,
\begin{gather} \label{e thetaj int}
\int\limits
_{\substack{\eta \in \R \\ \eta+i\ep \not \in G_d}}
\frac {1}{|M_+ (\la) - M_- (\la)|^2}
\, \im \frac 1{\theta_j^\pm - \lambda} \;
\left| \int_\R \frac{g^+ (s)d\Sigma_+ (s) } {s-\la }
\right|^2 d\eta \leq
C_4 \| g^+ \|_{L^2 (d\Sigma_+)}^2  \
\end{gather}
for $j=1$, \dots, $N_\theta^\pm$.

Assume $\theta_j^\pm \in \sigma_{disc} (A)$. Then
\begin{gather*} 
\im \frac{1}{\theta_j^\pm - \la}
\leq C_5 \im M_{ac\pm} (\la) , \quad
\la \in \C_+ \setminus G_d \ .
\end{gather*}
Thus \eqref{e thetaj int} follows from \eqref{e Iac}.

Assume $\theta_j^\pm \not \in \sigma_{disc} (A)$.
In this case,
$\theta_j^\pm \in \{\wt \theta_k^\pm \}_1^{\wt N_\theta}$.
Let $k$ be such that $\theta_j^\pm = \wt \theta_k^\pm$.
By assumptions of the theorem, conditions \eqref{e C theta u},
\eqref{e C theta w}, and \eqref{e C theta M} hold.
It is easy to see that
\begin{gather*}
\im \frac{1}{\wt \theta_k^\pm - \la}
\leq C_6 \im M_{ac\pm} (\la) , \quad
\la \in \C_+ \setminus D_k^\pm \ .
\end{gather*}
Therefore \eqref{e Iac} implies
\begin{gather}
\int\limits
_{\substack{\eta \in \R \\ \eta+i\ep \not \in \, ( D_k^\pm \cup \, G_d )}}
\frac {1}{|M_+ (\la) - M_- (\la)|^2}
\, \im \frac 1{\wt \theta_k^\pm - \lambda}
\left| \int_\R \frac{g^+ (s)d\Sigma_+ (s) } {s-\la }
\right|^2 d\eta \leq \notag
\\ \leq   \label{e thetaj}
C_4 \| g^+ \|_{L^2 (d\Sigma_+)}^2  \ .
\end{gather}

By \eqref{e C theta M}, we have
\begin{gather}
\int\limits
_{\substack{\eta \in \R \\ \eta+i\ep \in D_k^\pm }}
\frac {1}{|M_+ (\la) - M_- (\la)|^2}
\, \im \frac 1{\wt \theta_k^\pm - \lambda}
\left| \int_\R \frac{g^+ (s)d\Sigma_{d+} (s) } {s-\la }
\right|^2 d\eta \leq \notag \\
\leq
\int\limits
_{\substack{\eta \in \R \\ \eta+i\ep \in D_k^\pm }}
\frac {\ep}{(\wt \theta_k^\pm - \eta)^2 + \ep^2}
\frac {1}{|M_+ (\la) - M_- (\la)|^2 |\wt \theta_k^\pm - \la |^2}
\, \left| (\wt \theta_k^\pm - \la) \int_\R \frac{g^+ (s)d\Sigma_{d+} (s) } {s-\la }
\right|^2 d\eta \leq \notag \\
\leq
C_\theta^M \int\limits
_{\substack{\eta \in \R \\ \eta+i\ep \in D_k^\pm }}
\frac {\ep}{(\wt \theta_k^\pm - \eta)^2 + \ep^2}
\left| (\wt \theta_k^\pm - \la) \int_\R \frac{g^+ (s)d\Sigma_{d+} (s) } {s-\la }
\right|^2 d\eta . \label{e Pois}
\end{gather}
We may assume that
\[
D_k^\pm \cap \left( \{ \theta_j^+ \}_1^{N_\theta^+} \cup
\{ \theta_j^- \}_1^{N_\theta^- }\right) = \{\wt \theta_k^\pm \} \ .
\]
Therefore,
\[
\left| (\wt \theta_k^\pm - \la) \int_\R \frac{g^+ (s)d\Sigma_{d+} (s) } {s-\la }
\right| \leq C_7 \| g^+ \|_{L^2 (d\Sigma_+)}^2 ,
\quad \la \in D_k^\pm  \ .
\]
If we combine this with properties of Poisson kernel
(see \cite[Section I.3]{Gar}),
we get
\begin{gather} \label{e thetaj dP}
\int\limits
_{\substack{\eta \in \R \\ \eta+i\ep \in D_k^\pm }}
\frac {\ep}{(\wt \theta_k^\pm - \eta)^2 + \ep^2}
\left| (\wt \theta_k^\pm - \la) \int_\R \frac{g^+ (s)d\Sigma_{d+} (s) } {s-\la }
\right|^2 d\eta
\leq \pi C_7 \| g^+ \|_{L^2 (d\Sigma_+)}^2  .
\end{gather}
Using \eqref{e thetaj dP} and \eqref{e Pois}, we get
\begin{gather} \label{e thetaj d}
\int\limits
_{\substack{\eta \in \R \\ \eta+i\ep \in D_k^\pm }}
\frac {1}{|M_+ (\la) - M_- (\la)|^2}
\, \im \frac 1{\wt \theta_k^\pm - \lambda}
\left| \int_\R \frac{g^+ (s)d\Sigma_{d+} (s) } {s-\la }
\right|^2 d\eta \leq
\pi C_\theta^M C_7 \| g^+ \|_{L^2 (d\Sigma_+)}^2  \ .
\end{gather}

The inequality
\begin{gather} \label{e thetaj ac}
\int\limits
_{\substack{\eta \in \R \\ \eta+i\ep \in D_k^\pm }}
\frac {1}{|M_+ (\la) - M_- (\la)|^2}
\, \im \frac 1{\wt \theta_k^\pm - \lambda}
\left| \int_\R \frac{g^+ (s)d\Sigma_{ac+} (s) } {s-\la }
\right|^2 d\eta \leq
C_9 \| g^+ \|_{L^2 (d\Sigma_+)}^2
\end{gather}
follows from \eqref{e C theta u}, \eqref{e C theta w}, and
\eqref{e Cond w2}
in the same way as \eqref{e K<U/wg} follows from \eqref{e Cond u2},
\eqref{e Cond u3}, and \eqref{e Cond w2}.

Combining \eqref{e thetaj ac}, \eqref{e thetaj d}
and \eqref{e thetaj}, we get \eqref{e thetaj int}.
Thus \eqref{e Id} is proved.
Inequality \eqref{e EneqEss2+} is proved.
Inequality \eqref{e EneqEss2-} can be shown in the same way.
Thus Proposition \ref{p CrSimAe} yields that
$\Aess$ is similar to a selfadjoint operator.
\end{proof}

\section{Indefinite Sturm-Liouville operators with
finite-zone potentials \label{ch SL ss FZ} }

\subsection{Spectral properties of $\Aess$ and $\Adisc$}

Let $L = - d^2/dx^2 + q(x)$ be a Sturm-Liouville operator with
a finite-zone potential $q$ (see Subsection \ref{ss P FZ}).

In this case, we have
\begin{gather}
\sp (A_0^\pm ) = \sp_{ac} (A_0^\pm)
\cup \sp_{disc} (A_0^\pm),
          \label{e S=S+S FZ}
\end{gather}
\begin{gather}
\sp_{ac} (A_0^+) = - \sp_{ac} (A_0^-) = \sp (L) =
[\wheta_0 , \cheta_1 ] \cup [\wheta_1 , \cheta_2 ]
\cup \cdots \cup [\wheta_N , +\infty ) ,
\label{e sacA2pm}
\end{gather}
\begin{gather}
\sp_{disc} (A_0^\pm) = \{ \pm \tau_j :
\tau_j \not \in \{ \wheta_k\}_0^N \cup
\{\cheta_k \}_1^N , \
Q (\tau_j ) \pm  i \sqrt{R( \tau_j )} \neq 0 \} =:
\{ \theta_k^\pm \}_1^{N_\theta^\pm} .
\label{e SpDisc FZA2}
\end{gather}

Let $M(z)$ be a (multivalued) analytical function. If $M(z)= \suml_{k=-\infty}^{+\infty}
m_k (z-a)^{k/n} $ in a sufficiently small neighborhood of a point $a\in \C$, then  
we say that the number
\[
\frac 1n \inf \{k : m_k \neq 0\} \quad
(-\frac 1n \inf \{k : m_k \neq 0\})
\]
is \emph{the generalized order of a zero (pole) of the function } $M(z)$
\emph{at the point} $a$. Recall that the functions $M_\pm (\la)$
are holomorphic in $\rho (A_0^\pm)$.
For $\eta \in \sigma (A_0^\pm)$, we set
$M_\pm (\eta):= M_\pm (\eta +i0)$.
Note that in the case of a finite-zone potential $q$
the functions $M_\pm (\la)$ can be continued
on $\C$ as multivalued analytical functions with finite number of poles
and finite number of branch points.
Let us denote these continuations by $\widehat M_\pm (\lambda)$.
Then
\begin{description}
\item[] $\{\pm \wheta_j \}_0^N \cup
\{\pm \cheta_j \}_1^N$ are the sets of branch points for $ \widehat M_\pm (\la )$;
\item[]
$\{ \pm \xi_j \}_1^N \cap (\{\pm \wheta_j \}_0^N \cup
\{\pm \cheta_j \}_1^N)$
are the sets of zeroes of the generalized order $1/2$ for
$ \widehat M_\pm (\la )$;
\item[]
$\{ \pm \tau_j \}_0^N \cap (\{\pm \wheta_j \}_0^N \cup
\{\pm \cheta_j \}_1^N)$
are the sets of poles of the generalized order $1/2$ for
$\widehat M_\pm (\lambda)$;
\item[]
$ \{ \theta_j^\pm \}_1^{N_\theta^\pm}$ are the sets of
poles of the first order for $\widehat M_\pm (\lambda)$;
\item[]
$
\, \{ \pm \xi_j :
\xi_j \not \in \{ \wheta_k\}_0^N \cup
\{ \cheta_k \}_1^N , \
Q (\xi_j ) \mp  i \sqrt{R( \xi_j )} \neq 0 \} 
$ are the sets of zeroes of the first order for $\widehat M_\pm (\lambda)$.
\end{description}

We will say that $\lambda_0$ is \emph{a generalized zero (pole)} of $M_\pm $
if the generalized order of a zero (a pole) at $\lambda_0$ is positive.

We denote by $ \Mc_\pm$ the holomorphic continuation of $M_\pm (\la)$ from $\C_+$ to
\[
\C \setminus \left\{ \la : \ \im \la <0 , \ \re \la \in
\{ \pm \wheta_j \}_0^N \cup \{\pm \cheta_j \}_1^N \right\} .
\]

\begin{theorem} \label{t FZ Ad}
Let $L = - d^2/dx^2 + q(x)$ be a Sturm-Liouville operator
with a finite-zone potential $q$.
Let $A = (\sgn x) (- d^2/dx^2 + q(x))$.
Then:
\begin{description}
\item[1)] The operator $A$ has finite number
of eigenvalues,
\begin{gather}
\sp_p (A) =
\left( \{ \theta_j^+ \}_1^{N_\theta^+} \cap \{ \theta_j^- \}_1^{N_\theta^-}
\right) \, \cup \,
\{ \la \in \rho (A_2^+ \oplus A_2^-) :
 M_+ (\la) = M_- (\la) \} .  \label{e SpDisc FZ}
\end{gather}

\item[2)] The eigenvalues of $A$ are isolated and
have finite algebraic multiplicity, 
the geometric multiplicity equals one
for all eigenvalues of $A$.

\item[3)] If
$\la_0 \in \rho (A_2^+) \cap
\rho (A_2^-) $,
then the algebraic multiplicity of
$\la_0$ is equal to the multiplicity
of $\la_0$ as a zero of the
holomorphic function $M_+ (\la ) - M_- (\la )$;
if
$\la_0 \in \{ \theta_j^+ \}_1^{N_\theta^+}
\cap \{ \theta_j^- \}_1^{N_\theta^-} $,
then the algebraic multiplicity of
$\la_0$ is equal to the multiplicity
of $\la_0$ as a zero of the
holomorphic function $\dfrac 1{M_+ (\la )}
- \dfrac 1{M_- (\la )}$.

\item[4)] There exist a skew decomposition $ L^2 (\R ) = \Hsp_e \dot + \Hsp_d $ such that
\begin{gather}
A = \Aess \dot + \Adisc, \quad \Aess = A \upharpoonright ( \dom (A) \cap \Hsp_e),
\quad \Adisc = A \upharpoonright ( \dom (A) \cap \Hsp_d),  \notag \\
\sp (\Adisc) = \sp_{disc} (A), \quad \sp (\Aess) = \sp_{ess} (A).
\label{e skew fz}
\end{gather}
Besides, $\Hsp_d$  is a finite-dimensional space.
\end{description}
\end{theorem}

\begin{proof}
The spectral functions $\Si_\pm $ have the forms
$\Si_\pm (t) = \Si_{ac\pm} (t) + \Si_{d\pm} (t)$, where
\begin{gather}
\Si_{ac\pm} ^\prime (t) = \left\{ \begin{array}{c}
\frac{\sqrt{ R (\pm t)}} { S(\pm t) } ,
\quad t \in \pm \cupl_{j=0}^N (\wheta_j, \cheta_{j+1} )  \\
0 , \quad t \not \in \pm \cupl_{j=0}^N [\cheta_j, \wheta_j ]
\end{array} \right. \ .
          \label{e Sac pr}
\end{gather}
Here the branch of multifunction
$\sqrt{ R (\pm \la)}$ is chosen such that
$\Si_{ac\pm} ^\prime (t) \geq 0$ a.e.
(see \cite{KacKr}, \cite{Lev}).

Consequently,
\begin{gather}
\Si_{ac\pm} ^\prime (t)
\asymp 1 \quad  ( t \ra t_0 ) ,
\quad t_0 \in \pm \cupl_{j=0}^N (\wheta_j, \cheta_{j+1} ) ,
          \label{e Spr as 1}
\end{gather}
\begin{gather}
\Si_{ac\pm} ^\prime (t) \asymp |t-t_0|^{1/2} \chi_\pm  (t-t_0)
\quad (t \ra t_0 ) ,
\quad t_0 \in  \{\pm \wheta_j \}_0^N \setminus \{ \pm \tau_j \}_0^N ,
          \label{e Spr as 1/2wh}
\end{gather}
\begin{gather}
\Si_{ac\pm} ^\prime (t) \asymp |t-t_0|^{1/2} \chi_\mp (t-t_0)
 \quad (t \ra t_0 ) ,
\ t_0 \in  \{ \pm\cheta_j \}_1^N \setminus \{ \pm \tau_j \}_0^N ,
          \label{e Spr as 1/2ch}
\end{gather}
\begin{gather}
\Si_{ac\pm} ^\prime (t) \asymp |t-t_0|^{-1/2}\chi_\pm (t-t_0)
\quad (t \ra t_0 ) ,
\ t_0 \in  \{\pm \wheta_j \}_0^N \cap \{ \pm \tau_j \}_0^N ,
          \label{e Spr as -1/2wh}
\end{gather}
\begin{gather}
\Si_{ac\pm} ^\prime (t)
\asymp |t-t_0|^{-1/2}\chi_\mp (t-t_0)  \quad (t \ra t_0 ) ,
\ t_0 \in  \{ \pm\cheta_j \}_1^N \cap \{ \pm \tau_j \}_0^N .
          \label{e Spr as -1/2ch}
\end{gather}
Therefore,
\[
\intl_{\R \setminus \{\eta_0\}} \frac 1{|t-\eta_0|^2} d\Si_\pm (t) = \infty,
\quad \forall \eta_0 \in \pm (\; \cup_{j=0}^N [\wheta_j, \cheta_{j+1} ] \cup
[\wheta_N, +\infty) \; ) \ .
\]
Combining this with Theorem \ref{t s p} (1), we get
\[
\sp_p (A) \subset \C \setminus (\sp_{ac} (A_2^+) \cup \sp_{ac} (A_2^-))
= \C \setminus \sp_{ess} (A) .
\]
Thus, Proposition \ref{p e and d} yields \eqref{e SpDisc FZ}.

Taking into account \eqref{e M FZ}, we can write
the equation $M_+ (\la) = M_- (\la)$ in the form
\[
\frac {P( \la)}{ Q( \la) - i \sqrt{R( \la)}} =
\frac {P(- \la)}{ Q(- \la) + i \sqrt{R(- \la)}} ,
\]
where $P$, $Q$, and $R$ are polynomials.
Thus the equation $M_+ (\la) = M_- (\la)$
has finite number of solutions.
Therefore the set $\sp_p (A) $ is finite.
Statement (1) is proved.

Statements (2) and (3) follow from Statement (1) and
Proposition \ref{p e and d}.
Statement (4) follows from statements(1), (2), and
\eqref{e skew}.
\end{proof}

           \begin{theorem} \label{t FZ Ae}
Let $L = - d^2/dx^2 + q(x)$ be a Sturm-Liouville operator
with a finite-zone potential, let
$A = (\sgn x) (- d^2/dx^2 + q(x))$.
Then the following statements are equivalent:
\begin{description}
\item[(i)]
The operator $\Aess$ is similar to a selfadjoint operator;

\item[(ii)] The following conditions are satisfied
        \begin{gather} \label{e tAe2+}
        \frac{\im M_{\pm}} {M_+ (t) - M_- (t)} \in L^\infty (\R);
           \end{gather}
\item[(iii)] The function $ \Mc_+ (\la ) - \Mc_- (\la) $ has no generalized zeroes in
\begin{multline*}
(-\infty, -\wheta_N ) \cup (-\cheta_N , -\wheta_{N-1} ) \cup \cdots \cup
(-\cheta_1 , -\wheta_0 )
\cup \\ \cup
(\wheta_0 , \cheta_1 ) \cup (\wheta_1 , \cheta_2 )
\cup \cdots \cup (\wheta_N , +\infty ) \ ,
\end{multline*}
has no zeroes of the generalized order more than
$1/2$ in the set
\[
\left( (\{ \wheta_j \}_0^N \cup \{ \cheta_j \}_1^N) \setminus
\{ \tau_j \}_0^N \right) \cup
\left( (\{ -\wheta_j \}_0^N \cup \{ -\cheta_j \}_1^N) \setminus
\{ -\tau_j \}_0^N
\right) ,
\]
has poles of generalized order greater than or equal to
$1/2$ at the points of the set
\[
\left( (\{ \wheta_j \}_0^N \cup \{ \cheta_j \}_1^N) \cap
\{ \tau_j \}_0^N \right) \cup
\left( (\{ -\wheta_j \}_0^N \cup \{ -\cheta_j \}_1^N) \cap
\{ -\tau_j \}_0^N
\right) .
\]
\end{description}
          \end{theorem}

Combining Theorem  \ref{t FZ Ae}  with  Corollary  \ref{cor2.5abscont} we arrive at the
following result.

 \begin{corollary}
Under the conditions \eqref{e tAe2+} the operator $\Aess$ is similar to a selfadjoint
operator with absolutely continuous spectrum.
            \end{corollary}

     \begin{proof}
Consider the decomposition  \eqref{e skew} and note that the subspace $\Hsp_e$ in
\eqref{e skew} is invariant for the operator  $A,\ \Hsp_e \in \Lat A.$ Alongside the
skew decomposition \eqref{e skew}  we consider the orthogonal decomposition $\Hsp = \Hsp_e \oplus \Hsp_e^\perp.$  According to this decomposition the
characteristic function $\theta_A(\cdot)$ of the operator   $A$ admits the factorization
$\theta_A(\la)=\theta_1(\la)\cdot\theta_2(\la)$ where $\theta_1(\cdot)$ is the
characteristic function of the operator $\Aess=A\lceil \Hsp_e$ and $\theta_2({\cdot})$
is the characteristic function of the operator  $A_2 :=P_2 A\lceil \Hsp^{\perp}_e$, where $P_2$ is the orthoprojection   in  $\Hsp $  onto 
$\Hsp^{\perp}_e$.  Note, that $\theta_2(\cdot) = \theta_{A_2}(\cdot)$ is a finite Blaschke
product since $\sigma(\Adisc)$ is finite.

It follows from \eqref{e M FZ2} that  $M_+(\cdot)$  (resp. $M_-(\cdot)$ ) admits a
continuous extension to the real line with exception  of the set of (real) zeros
$\{s_k\}_1^{N+1}$\  (resp. $\{-s_k\}_1^{N+1}$) of the polynomial $S(\la)$ (resp.
$S(-\la)$). Moreover, it is clear from the formula   \eqref{formulaforcharfunc}  for the
characteristic function $\theta_A(\la)$ that  real singularities (resp. poles) of
$\theta_A(\la)$ coincide with the set of real (resp. non-real) roots of the function
        \begin{equation}
F(\la)=P(\la)Q(-\la)+i P(\lambda)\sqrt{R(-\la)}-P(-\la)Q(\la)+iP(-\la)\sqrt{R(\la)}.
       \end{equation}
In particular, the numbers of real singularities and poles of $\theta_A(\cdot)$  are
finite.

Note, that $A^*_2=A^*\lceil \Hsp^{\perp}_{\theta}$ and
$\theta^{-1}_2(\la)=\theta_{A^*_2}(\la)$ is a finite Blaschke product too. Therefore the
sets of real singularities
 of functions $\theta(\cdot)$ and  $\theta_1(\cdot)=\theta(\cdot)\cdot\theta^{-1}_2(\cdot)$
 coincide. In particular,   $\theta_1(\cdot)$  may have only finite number of
singularities and we can apply Proposition  \ref{prop2.5abscont}. Therefore, combining
Theorem   \ref{t FZ Ae}  with Proposition   \ref{prop2.5abscont} we obtain that $\Aess$ is
similar to a selfadjoint operator with absolutely continuous spectrum.
      \end{proof}
      
             \begin{corollary} \label{c FZ>0}
Let $L= - d^2/dx^2 +q(x)$ be a nonnegative Sturm-Liouville operator with a finite-zone
potential $q$. Then the operator $A=(\sgn x)L$ is similar to a selfadjoint operator with
absolutely continuous spectrum.
             \end{corollary}

\begin{proof}
By \eqref{e Sac pr}, we have
$\supp d\Si_{ac\pm} (t) \subset \R_\pm$ and
$\im M_\pm (t) = \pi \Si_{ac\pm} ^\prime (t) $ for almost all $t \in \R$.
Therefore,
\[
\frac{ |\Si_{ac\pm} ^\prime (t) |^2} {| M_+ (t) - M_- (t)|^2}
\leq
\frac{ | \Si_{ac\pm} ^\prime (t) |^2 }
{ \pi ^2 | \Si_{ac\pm} ^\prime (t) |^2
+ | \re M_+ (t) - \re M_- (t)|^2} \leq \frac 1{\pi^2} ,
\]
for almost all $t \in \R$.
Thus, by Theorem \ref{t FZ Ae}, $\Aess$ is similar to a selfadjoint operator.

The operator $A$ is J--nonnegative. Besides, $L$ has an absolutely continuous spectrum
(see, for example, \cite{Lev}). Hence, $\ker L = 0$. Combining this with Proposition
\ref{p spJ>0}, we see that all the eigenvalues of $A$ are real and simple. Therefore, by
Theorem \ref{t FZ Ad}, the operator $\Adisc$ is similar to a  selfadjoint operator. Thus,
$A$ is similar to a selfadjoint operator.
\end{proof}

\subsection{Proof of Theorem \ref{t FZ Ae}}

The implication $ (i) \Rightarrow (ii)$ follows from Proposition \ref{p NessAess}.

$ (ii) \Rightarrow (iii).$ It follows from \eqref{e Spr as 1} and
\eqref{e tAe2+} that there are no generalized zeroes of the function $\Mc_+ (\la) - \Mc_- (\la)$ in the set $\cup_{j=0}^N (\wheta_j, \cheta_{j+1} )$.
Likewise, it follows from \eqref{e Spr as 1} and \eqref{e tAe2+} that
there are no generalized zeroes of the function $\Mc_+ (\la) - \Mc_- (\la)$ in the set $\cup_{j=0}^N (-\cheta_{j+1}, -\wheta_{j} )$.

It follows from \eqref{e Spr as 1/2wh}, \eqref{e Spr as 1/2ch},
and \eqref{e tAe2+} that
the function $\Mc_+ (\la) - \Mc_- (\la)$ has no zeroes of generalized order greater than $1/2$ in the sets
\[
 (\{ \wheta_j \}_0^N \cup \{ \cheta_j \}_1^N) \setminus \{ \tau_j \}_0^N
\quad \text{and} \quad
 (\{ - \wheta_j \}_0^N \cup \{ - \cheta_j \}_1^N)
\setminus \{ -\tau_j \}_0^N .
\]
It follows from \eqref{e Spr as -1/2wh}, \eqref{e Spr as -1/2ch},
\eqref{e tAe2+}, and \eqref{e tAe2+} that
all the points of the sets
\[
(\{ \wheta_j \}_0^N \cup \{ \cheta_j \}_1^N) \cap \{ \tau_j \}_0^N
\quad \text{and} \quad
(\{ - \wheta_j \}_0^N \cup \{ - \cheta_j \}_1^N) \cap \{ -\tau_j \}_0^N
\]
are generalized poles of $\Mc_+ (\la) - \Mc_- (\la)$.
The generalized orders of these poles are greater than or equal to $1/2$.

$ (iii) \Rightarrow (i).$
By Theorem \ref{t FZ Ad} (1), Assumption  \eqref{a sdisc} is fulfilled for the operator $A$.
It follows from \eqref{e S=S+S FZ} that we can apply Theorem
\ref{t uwSufCond}.

Let Statement (iii) be fulfilled.
We construct the functions $U_\pm (\la)$, $w_\pm(t)$, $U_j^\pm$ and the sets
$G_d$, $D_j^\pm$ such that all the conditions of Theorem \ref{t uwSufCond} hold true.

Let $G_d$ be any compact set such that
$\sp_{ess} (A) \cap G_d = \emptyset$ and
all the points of the set $\sp_{disc} (A) $
are interior points of $G_d$.

The set $\sp_{disc} (A) \cap \C_+$ is finite.
Besides,
\[
\sp_{disc} (A) \cap \C_+ = \{\la \in \C_+ \ : \ M_+ (\la) - M_- (\la) = 0\}.
\]
Let $B_\C (\la)$ be a finite Blaschke product (see \cite{Gar}) with the same zeroes in $\C_+$ as $M_+ (\la) - M_- (\la)$.
Then $ M_+ (\la) - M_- (\la) = B_\C (\la) M_1 (\la)$, the function $M_1 (\la)$ being holomorphic on $\C_+$.
Besides,
\begin{gather}
\frac 1{M_1 (\la)} \in H(\C_+), \qquad
M_1 (\la) \asymp (M_+ (\la)- M_- (\la)) \quad (\la \in \C_+ \setminus G_d^+),
\label{e M1}
\end{gather}
where $G_d^+ $ is any compact subset of $G_d \cap \C_+$ such that all the points of the set $\sp_{disc} (A) \cap \C_+$ are interior points of $G_d^+ $.

The set
$\sp_{disc} (A_2^+) \cap \sp_{disc} (A_2^-) = \{\theta_j^+ \}_1^{N_\theta^+}
\cap \{ \theta_j^- \}_1^{N_\theta^-}$
is finite.
By Theorem \ref{t FZ Ad}, this set is a subset of $\sp_{disc} (A)$.
Let
\[
\; \{ \theta_j \}_{j=1}^{N_\theta} := \sp_{disc} (A_2^+) \cap \sp_{disc} (A_2^-)  ,
\quad N_\theta < \infty .
\]
Each point of the set $ \{ \theta_j \}_{j=1}^{N_\theta} $ is either
a pole of the first order or a removable singularity of the function
$M_+ (\la ) - M_- (\la)$.
By $\kappa_j$ denote the generalized order of a zero of $M_+ (\la ) - M_- (\la)$
at $\theta_j$. Then $\kappa_j \in \{-1,0\} \cup \N $,
$j=1,\dots, N_\theta$.
By Theorem \ref{t FZ Ad}, we have
\[
\,  \{ \wt \theta_j^\pm \}_1^{\wt N_\theta^\pm} =
\{ \theta_j^\pm \}_1^{N_\theta^\pm}\setminus \{\theta_j \}_1^{N_\theta}
\]
(the sets $\{ \wt \theta_j^\pm \}_1^{\wt N_\theta^\pm}$
are defined by \eqref{e wt theta}).

Put
\[
\, \{\wt \theta_j \}_1^{\wt N_\theta} =
(\R \cap \sp_{disc} (A) ) \setminus \{\theta_j \}_1^{N_\theta} \ .
\]
The functions $M_\pm (\la) $ are regular at
$\wt \theta_j$ and $M_+ (\wt \theta_j) - M_- (\wt \theta_j) = 0$,
$j=1,\dots, \wt N_\theta$. Let us denote generalized
order of $\wt \theta_j$ as a zero of $M_+ (\la ) - M_- (\la)$
by $\wt \kappa_j$ (clearly, $ \wt \kappa_j \in N$).

Put  $M_2 (\la) := M_1 (\la)/B_\theta$, where
\begin{gather*}
B_{\theta} (\la) :=
\frac{ \prodl_{j=1}^{N_\theta} (\la -\theta_j)^{\kappa_j+1} }
     { \prodl_{j=1}^{N_\theta} \left( \la -(\theta_j-i\ep_1)
\right)^{\kappa_j+1} } \;
\frac{ \prodl_{j=1}^{\wt N_\theta} (\la -\wt \theta_j)^{\wt \kappa_j} }
     { \prodl_{j=1}^{\wt N_\theta} \left( \la -(\wt \theta_j-i\ep_1)
\right)^{\wt \kappa_j}}.
\end{gather*}
Here and below $\ep_1$ is an arbitrary fixed positive number.
Taking into account \eqref{e M1}, we get
\begin{gather}
\frac 1{M_2 (\la)} \in H(\C_+) , \qquad
M_2 (\la) \, \asymp \, \left( M_+ (\la)- M_- (\la)
\right) \quad (\la \in \C_+ \setminus G_d). \label{e M2}
\end{gather}
Denote
\[
\rho_1 := \rho (L) \cup \rho (-L)  \ .
\]
If $\la_0$ is a generalized zero of $\Mc_+ (\la) - \Mc_- (\la)$ and
$\la_0 \in \rho_1$, then $\la_0 \in \{ \theta_j \}_{j=1}^{N_\theta} \cup
\{ \wt \theta_j \}_{j=1}^{\wt N_\theta}$.
Moreover, it follows from Statement (iii) that the function
$\Mc_+ (\la) - \Mc_- (\la)$ has no generalized zeroes in the set
$ \si_1 := \si_1^+ \cup \si_1^- $,
where
\[
 \si_1^\pm := \pm \bigcup_{j=0}^N (\wheta_j , \cheta_{j+1})
\]
are the sets of interior points of the spectra $ \sp (\pm L)$.
Therefore the definition of $B_\theta $ imply that
\begin{gather}
M_2^{-1} (\la) = O(1) \quad (\la \ra \la_0) \quad
\forall \la_0 \in \rho_1 \cup \si_1 , \label{e M2 O1}
\\
M_2 (\la) \asymp (\theta^\pm_j - \la)^{-1} \quad ( \la \ra \theta^\pm_j), \quad
j=1, \dots , N_\theta^\pm . \label{e asM2theta}
\end{gather}
Let us explain formula \eqref{e asM2theta}.
If $\theta^\pm_j \in \{ \theta_j \}_1^{N_\theta}$, the formula \eqref{e asM2theta} follows from the definition of the function $B_\theta$.
If $\theta^\pm_j \in \{ \wt \theta_j^\pm \}_1^{\wt N_\theta^\pm}$,
the asymptotics
\[
M_\pm (\la) \asymp (\wt \theta^\pm_j - \la)^{-1} , \quad
M_\mp (\la) = O (1) \, (\wt \theta^\pm_j - \la)^{-1/2} \quad
(\la \ra \wt \theta^\pm_j)
\]
and \eqref{e M1} imply \eqref{e asM2theta}.

Let us denote
\[
\, \{ \zeta_j^\pm \}_1^{N_\zeta^\pm} := \left\{ z \in \{ \pm \wheta_j \}_0^N \cup \{ \pm \cheta_j  \}_1^N \ : \ z \ \text{is a generalized zero of }
\Mc_+ (\la ) - \Mc_- (\la ) \right\} .
\]
By Statement (iii), the generalized orders of all
the zeroes $\zeta_j^\pm$ are equal to $1/2$.
It follows from Statement (iii)
and asymptotics for $M_\pm (\la) $ that
\begin{gather}
\; \{ \zeta_j^\pm \} \subset
( \{ \pm \wheta_j \}_0^N \cup \{ \pm \cheta_j  \}_1^N )
\setminus ( \{ \theta_j^\mp \}_1^{N_\theta^\mp } \cup \si_1 ).
\label{e zeta subset}
\end{gather}
Denote
\begin{gather} \label{e wt zeta}
\, \{ \zeta_j \}_1^{N_\zeta } := \{ \zeta_j^+ \}_1^{N_\zeta^+}
 \cap \{ \zeta_j^- \}_1^{N_\zeta^-} ,
\quad
\{\wt {\zeta_j^\pm}_1 \}^{\wt {N_\zeta^\pm}} :=
\{ \zeta_j^\pm \}_1^{N_\zeta^\pm} \setminus
\{ \zeta_j \}_1^{N_\zeta } \ .
\end{gather}
Statement (iii) imply $\wt {\zeta_j^\pm } \not \in \si_1 $.
Besides,
\[
\wt {\zeta_j^\pm } \not \in \{ \zeta_j \}_1^{N_\zeta } =
\{ \zeta_j^\pm \}_1^{N_\zeta^\pm} \cap
( \{ \mp \wheta_j \}_0^N \cup \{ \mp \cheta_j  \}_1^N ) ,
\]
therefore
$\ \wt {\zeta_j^\pm } \in \rho_1^\mp $,
where
\[
 \rho_1^\pm := \pm \rho (L) \;
(=\pm \bigcup_{j=0}^N ( \cheta_j , \wheta_j ) ) \ .
\]
Put
\begin{gather}   \label{e Def u}
u_\pm (\la) := \frac { \sqrt{R(\pm\la)} } {S(\pm\la)} \;
\frac{ \prodl_{j=1}^{N_\theta^\pm} (\la -\theta_j^\pm) }
     { \prodl_{j=1}^{N_\theta^\pm} (\la -(\theta_j^\pm-i\ep_1)) } \;
\frac{ \prodl_{j=1}^{\wt N_\zeta^\pm} (\la -\wt \zeta_j^\mp) }
     { \prodl_{j=1}^{\wt N_\zeta^\pm} (\la -(\wt \zeta_j^\mp-i\ep_1)) } ,
\end{gather}

Now we define $U_\pm$ by
\begin{gather*}
U_\pm := \frac {\sqrt{u_\pm (\la)}}{M_2 (\la)} .
\end{gather*}

Let us check conditions \eqref{e Cond u2}, \eqref{e Cond uN+}, and
\eqref{e Cond u3}.
All the asymptotics given below are considered on $\overline{\C_+}$,
unless otherwise specified.

\begin{lemma}
Let Statement (iii) be true. Then condition \eqref{e Cond u2} is
fulfilled, i.e.,
\begin{gather*} 
\frac{\im M_{ac\pm} (\la)}{|M_+ (\la) - M_- (\la)|^2}
\leq C_\pm^u | U_\pm (\la)|^2 ,
\quad \la \in \C_+ \setminus G_d \ .
\end{gather*}
\end{lemma}

\begin{proof}
By \eqref{e M2}, condition \eqref{e Cond u2} is equivalent to
\begin{gather}
\im M_{ac\pm} (\la) = O(1) \, u_\pm (\la)  , \qquad (\la \in \overline{\C_+ \setminus G_d} ).
\label{e ImM<u}
\end{gather}
Since
\begin{gather}
M_{ac\pm} (\la) \asymp M_\pm (\la) \asymp |\la|^{-1/2} \quad
(\la \ra \infty )  \label{e ImM=M=}
\end{gather}
\begin{gather}
\text{and } \qquad u_\pm (\la ) \asymp |\la|^{-1/2} \quad
(\la \ra \infty ) ,
          \label{e uInf}
\end{gather}
\begin{gather*}
\text{we have } \qquad
\frac{\im M_{ac\pm} (\la ) }{ u_\pm (\la) } = O(1)
\quad (\la \ra \infty ) .
\end{gather*}

If $\la_0 \in \overline{(\C_+ \setminus G_d )} \,
\setminus \, \left( \{ \pm \wheta_j\}_0^N \cup
\{ \pm \cheta_j\}_1^N \cup \{ \wt {\zeta_j^\mp} \}_1^{\wt { N_\zeta^\mp } }
\right) $, then
\begin{gather*}       
\im M_{ac\pm} (\la ) = O (1) \quad (\la \ra \la_0 ), \qquad
u_\pm (\la ) \asymp 1 \quad (\la \ra \la_0 ) .
\end{gather*}

Let $\la_0 \in (\{ \pm \wheta_j\}_0^N \cup
\{ \pm \cheta_j\}_1^N ) \setminus \{ \pm \tau_j\}_0^N $.
Then \eqref{e M FZ2} yields
\begin{gather*}
\im M_{ac\pm} (\la ) \asymp \im M_\pm (\la )
= O (|\la -\la_0|^{1/2}) \quad (\la \ra \la_0 ) ;
\end{gather*}
\begin{gather*}
\text{besides,}
 \qquad
u_\pm (\la ) \asymp |\la -\la_0|^{1/2} \quad (\la \ra \la_0 ) .
\end{gather*}

Let $\la_0 \in (\{ \pm \wheta_j\}_0^N \cup
\{ \pm \cheta_j\}_1^N ) \cap \{ \pm \tau_j\}_0^N $.
Then \eqref{e M FZ2} yields
\begin{gather*}
\im M_{ac\pm} (\la ) \asymp \im M_\pm (\la )
= O ( |\la -\la_0| ^{-1/2} ) \quad (\la \ra \la_0 ) ;
\end{gather*}
\begin{gather*}
\text{besides,}
\qquad u_\pm (\la ) \asymp |\la -\la_0|^{-1/2} \quad (\la \ra \la_0 ) .
\end{gather*}

Let $\la_0 \in \{ \wt {\zeta _j^\mp} \}_0^{\wt {N_\zeta^\mp}} $.
Then \eqref{e zeta subset} and \eqref{e intM}
yield
\begin{gather*}
\im M_{ac\pm} (\la ) \asymp \im M_\pm (\la )= O ( \im \la) = O (\la - \la_0 ) \quad (\la \ra \la_0 ) \ .
\end{gather*}
On the other hand,
\begin{gather*}
u_\pm (\la ) \asymp |\la -\la_0| \quad (\la \ra \la_0 ) .
\end{gather*}
If we combine all these estimates, we get \eqref{e ImM<u}.
Thus \eqref{e Cond u2} is proved.
\end{proof}

\begin{lemma}
Condition \eqref{e Cond uN+} is fulfilled, i.e.,
$U_\pm (\la) \in \Np (C_+) $.
\end{lemma}

\begin{proof}
The functions $U_\pm (\la)$ are holomorphic on $\C_+$ by definition.
Since
\[
M_+ (\la) - M_- (\la) \asymp |\la |^{-1/2}  \quad  \ (\la \ra \infty) \ ,
\]
\eqref{e uInf} imply  the following formula
\begin{gather} \label{e UInf}
U_\pm (\la) \asymp |\la|^{1/4} \qquad (\la \ra \infty).
\end{gather}
Condition \eqref{e Cond uN+} follows from \eqref{e UInf} and
Lemma~\ref{l N+}.
\end{proof}

\begin{lemma}
Let Statement (iii) be true.
Then condition \eqref{e Cond u3} is fulfilled, i.e.,
\begin{gather*}
\frac{U_\pm (t)} {\theta_j^- - t }
\in L^2(\R), \ j=1, \cdots , N_\theta^- ; \quad
\frac {U_\pm (t)}{ \theta_j^+ - t}
\in L^2(\R), \ j=1, \cdots , N_\theta^+ \ . 
\end{gather*}
\end{lemma}

\begin{proof}
The definition of the polynomial $S(\la)$ imply
\begin{gather}
|u_\pm (\la )| = \frac
{\prodl_{(\pm\la_0) \in (\{ \wheta_j \}_0^N \cup \{ \cheta_j \}_1^N) \setminus
\{ \tau_j \}_0^N } |\la - \la_0 |^{1/2}}
{\prodl_{(\pm\la_0) \in \{ \tau_j \}_0^N \cap( \{\wheta_j \}_0^N
\cup \{\cheta_j \}_1^N )} |\la - \la_0 |^{1/2}  \
 \prodl_{j=1}^{N_\theta^\pm} |\la -(\theta_j^\pm -i\ep_1)|}
\notag \times \\ \times
\frac
{ \prodl_{j=1}^{\wt N_\zeta^\mp} |\la -\wt \zeta_j^\mp| }
{  \prodl_{j=1}^{\wt N_\zeta^\mp} |t -(\wt \zeta_j^\mp -i\ep_1)|^{1/2} } .
\label{e |u|}
\end{gather}
It follows from \eqref{e |u|}, \eqref{e M2 O1}, \eqref{e Def u}, \eqref{e asM2theta}, Statement (iii), and the definition of
$ \{ \wt \zeta_j^\mp \}_1^{\wt N_\zeta^\mp }$
that
\begin{gather}
U_\pm (\la) = O(1) \quad (\la \ra \la_0) , \quad \la_0 \in
\sp_1^+ \cup \sp_1^- \cup \rho_1^\pm   \label{e U sp rho}
\\
U_\pm (\la) \asymp (\la - \theta_j^\pm) \quad (\la \ra \theta_j^\pm), \quad
j= 1 , \dots , N_\theta^\pm , \label{e UTheta1}
\\
U_\pm (\la) = O(1)\, |\la - \theta_j^\mp|^{3/4} \quad (\la \ra \theta_j^\mp),
\quad j= 1 , \dots , N_\theta^\mp,        \label{e UTheta2}
\\
U_\pm (\la) = O(|\la - \la_0|^{-1/4}) \ \  (\la \ra \la_0) ,
\quad \la_0 \in \left( \{ \pm \wheta_j \}_{j=0}^{N} \cup \{ \pm \cheta_j \}_{j=1}^{N}
\right) \setminus \{ \pm \tau_j \}_{j=0}^{N} ,  \label{e U mu-tau}
\\
U_\pm (\la) = O(|\la - \la_0|^{1/4}) \ \ (\la \ra \la_0) , \quad
\la_0 \in \left( \{ \pm \wheta_j \}_{j=0}^{N} \cup \{ \pm \cheta_j \}_{j=1}^{N}
\right) \cap \{ \pm \tau_j \}_{j=0}^{N} .    \label{e U mu tau}
\\
U_\mp (\la) = O(|\la - \la_0|^{1/4}) \ \ (\la \ra \la_0) , \quad
\la_0 \in \left( \{ \pm \wheta_j \}_{j=0}^{N} \cup \{ \pm \cheta_j \}_{j=1}^{N}
\right) \cap \{ \pm \tau_j \}_{j=0}^{N} .    \label{e Ump mu tau}
\end{gather}
Therefore, $  \dfrac {U_\pm (t)}{\theta_j^+ - t} \in L^2_{loc} (\R) $ \ and
\ $ \dfrac {U_\pm (t)}{\theta_j^- - t} \in L^2_{loc} (\R) $.
Combining this with \eqref{e UInf}, we get \eqref{e Cond u3}.
\end{proof}

Let $w_\pm (t)$ be defined by
\begin{gather*}
\frac 1{w_\pm (t)} := \left|
\frac{\sqrt{R(\pm t)}}{S(\pm t)} \;
\frac{ \prodl_{j=1}^{N_\theta^\pm} (t -\theta_j^\pm) }
     { \prodl_{j=1}^{N_\theta^\pm} (t -(\theta_j^\pm-i\ep_1)) } \;
\frac{ \prodl_{j=1}^{\wt N_\zeta^\mp} (\la -\wt \zeta_j^\mp)^{1/2} }
     { \prodl_{j=1}^{\wt N_\zeta^\mp} (\la -(\wt \zeta_j^\mp-i\ep_1))^{1/2} }
\right| .
\end{gather*}

Let us check conditions \eqref{e Cond w1},
\eqref{e Cond w2}, and \eqref{e Cond w3}.

Since all the points $\theta_j^\pm$, $\zeta_j^\mp$ belongs to
$\rho_1^\pm (= \R \setminus \supp d\Si_{ac\pm} )$, formulae
\eqref{e Spr as 1}--\eqref{e Spr as -1/2ch} imply
\eqref{e Cond w1}.

\begin{lemma} \label{l A2}
Condition \eqref{e Cond w2}
is fulfilled, i.e., the weights $w_+ $ and
$w_-$ satisfy the $(A_2)$ condition.
\end{lemma}

We give two proofs of this lemma.
The fist proof is based on the Hunt- Muckenhoupt-Wheeden theorem, the second on the Helson-Szeg\"o theorem. Note that \cite[Theorem 4]{HS60} can be used also.

\begin{proof}[Proof 1 of Lemma \ref{l A2} ]
It is clear that all the conditions of Proposition \ref{p A2} is fulfilled for the functions $w_\pm$. Thus, $w_\pm \in (A_2)$.
\end{proof}

\begin{proof}[Proof 2 of Lemma \ref{l A2} 
]
The Helson-Szeg\"o condition (see \eqref{e HSth}) is equivalent to
the $(A_2)$ condition.
Let us prove that condition \eqref{e HSth} is satisfied for $w_+$.

Obviously,
\begin{gather} \label{e w2}
w_+ (t) = \frac{ \prodl
_{\la_0 \in \{ \tau_j \}_0^N \cap( \{\wheta_j \}_0^N \cup \{\cheta_j \}_1^N )}
|t - \la_0 |^{1/2} \
 \prodl_{j=1}^{N_\theta^+} |t -(\theta_j^+ -i\ep_1)| \
 \prodl_{j=1}^{\wt N_\zeta^-} |t -(\wt \zeta_j^- -i\ep_1)|^{1/2} }
{ \prodl_{\la_0 \in (\{ \wheta_j \}_0^N \cup \{ \cheta_j \}_0^N) \setminus
\{ \tau_j \}_0^N } |t - \la_0 |^{1/2}  \
 \prodl_{j=1}^{\wt N_\zeta^-} |t -\wt \zeta_j^-|^{1/2} } .
\end{gather}
Consequently,
\begin{gather}
\log w_+ (t) = \frac 1 2 \suml
_{ \la_0 \in \{ \tau_j \}_0^N \cap( \{\wheta_j \}_0^N \cup \{\cheta_j \}_1^N ) }
\log |t -\la_0 | +
 \suml_{j=1}^{N_\theta^+} \log |t -(\theta_j^+ -i\ep_1)|
+ \notag \\ +
\frac 1 2 \suml_{j=1}^{\wt N_\zeta^-} \log |t -(\wt \zeta_j^- -i\ep_1)|
- \frac 1 2 \suml_{\la_0 \in (\{ \wheta_j \}_0^N \cup \{ \cheta_j \}_0^N) \setminus
\{ \tau_j \}_0^N } \log |t - \la_0 |
- \notag \\ -
\frac 1 2 \suml_{j=1}^{\wt N_\zeta^-} \log |t -\wt \zeta_j^-|
= (\H v_+) (t) + c_1 ,
          \label{e v+ logw}
\end{gather}
where $\H$ is the Hilbert transform (see Subsection \ref{ss Hp}),
\begin{gather*}
v_+ (t) = \frac 1 2 \suml_{\la_0 \in (\{ \wheta_j \}_0^N \cup \{ \cheta_j \}_0^N) \setminus
\{ \tau_j \}_0^N } \arg (t - \la_0 ) +
\frac 1 2 \suml_{j=1}^{\wt N_\zeta^-} \arg (t -\wt \zeta_j^-)
- \notag \\ -
\frac 1 2 \suml
_{ \la_0 \in \{ \tau_j \}_0^N \cap( \{\wheta_j \}_0^N \cup \{\cheta_j \}_1^N ) }
\arg (t -\la_0 ) -
 \suml_{j=1}^{N_\theta^+} \arg (t -(\theta_j^+ -i\ep_1))
- \notag \\ -
 \frac 1 2 \suml_{j=1}^{\wt N_\zeta^-} \arg (t -(\wt \zeta_j^- -i\ep_1)) ;
\end{gather*}
here $c_1 $ is a constant, the branch of $\arg z$ is fixed by $ \ \arg z \in (- \pi , \pi ] $,
$ \ z \in \C$.

The function $v_+$ is bounded and piecewise smooth;
the set of jumps of $v_+$ is
\[
\, \{ \wheta_j \}_0^N \cup \{ \cheta_j \}_1^N \cup
\{ \wt {\zeta_j^- }\}_1^{\wt {N_\zeta^-}}  \ .
\]
The absolute values of all the jumps are equal to $\pi /2$.
Moreover,
\[
v_+ (t) \asymp \arctan \frac 1t \asymp \frac 1t \quad (t \ra +\infty),
\]
\[
v_+ (t) + \frac \pi 2 \asymp \arctan \frac 1{|t|} \asymp \frac 1{|t|}
\quad (t \ra -\infty),
\]
$v_+ (t)$ monotonically increases on
$t \in (-\infty , \wheta_0 )$ and $ ( \wheta_N , +\infty)$.
Therefore, $v_+$ can be represented in the form
\[
v_+ (t) = v_1 (t) + v_2 (t) - \pi /4 ,
\]
where $v_1$ is a piecewise continuous function such that
\begin{gather} \label{e v1}
\| v_1 (t) \|_{L^\infty} < \pi /2 ,
\end{gather}
\[
\, v_1 \ \text{has jumps at the points} \quad
 \{ \wheta_j \}_0^N \cup \{ \cheta_j \}_1^N \cup
\{ \wt {\zeta_j^+ }\}_1^{\wt {N_\zeta^+}} \ ,
\]
$ v_2 $ is a $C^1$ function on $\R$ such that
\begin{gather} \label{e v2}
v_2 (t) = 0 \quad \text{for} \ t \not
\in [\wheta_0 - \delta_2  ,\wheta_N + \delta_2] ;
\end{gather}
here $ \delta_2 $ is a specified positive number.

From \eqref{e v2}, we get
\[
 \quad
(H v_2) (t) \asymp |t|^{-1} \quad (|t| \ra \infty) \ .
\]

It follows from $ v_2 \in C^1 (\R) $ that
$v_2 \in \mathrm{Lip}^\alpha (\I)$ for any compact interval
$\I \subset \R$ and for any $\alpha \in (0,1)$. If we combine this with Privalov's theorem (see \cite{Koosis})
and \eqref{e v2}, we get
$\quad H v_2 \in
\mathrm{Lip}^\alpha (\I)$, $\ 0< \alpha <1 $.
Hence, $H v_2 $ is a continuous function on $\R$ and
\eqref{e v2} imply
$ H v_2 \in L^\infty (\R) \, $. \
Taking into account \eqref{e v+ logw}, we get
$ \
\log w_+ (t) = (H v_1) (t) + (H v_2) (t) +c_1 $,
where $ \ \| v_1 \|_{L^\infty} < \pi/2$, \
$ \ H v_2 + c_1 \in L^\infty (\R)$. \
That is $ w_+ $ satisfy  the Helson-Szeg\"o condition.
The condition \eqref{e Cond w2} is proved for $w_+$.
In the same way we prove \eqref{e Cond w2} for $w_-$.
\end{proof}

\begin{lemma} \label{l Cond w3}
Let Statement (iii) be true. Then
condition \eqref{e Cond w3}
is fulfilled, i.e.,
\begin{gather*} 
\frac { U_+^2 (t)}{ w_\pm (t)}
\in L^\infty (\R), \qquad
\frac { U_-^2 (t)}{ w_\pm (t)}
\in L^\infty (\R)  .
\end{gather*}
\end{lemma}

\begin{proof}
Note that
\begin{gather}
w_+^{-1} (t) \asymp |t|^{-1/2} \quad (|t| \ra \infty) .
          \label{w as inf}
\end{gather}
It follows from \eqref{e UInf}, \eqref{e U sp rho}-\eqref{e Ump mu tau}, \eqref{e w2},
Statement (iii), and \eqref{e M2} that
\begin{gather*}
U_+^2 (t) w_+^{-1} (t) \in L^\infty (\R) \quad \text{and}
\end{gather*}
\begin{gather}
U_-^2 (t) w_+^{-1} (t) = O (1)
\quad (t \ra t_0) \notag \\
\text{for} \quad
t_0 \in \{-\infty \} \cup \{ + \infty \} \cup \rho_1^- \cup \si_1^-
\cup \si_1^+ \cup \{\wheta_j \}_0^N \cup \{ \cheta_j \}_1^N  .
          \label{e u-w+1}
\end{gather}
Note that
\begin{gather}
\R \setminus \left( \rho_1^- \cup \si_1^-
\cup \si_1^+ \cup \{\wheta_j \}_0^N \cup \{ \cheta_j \}_1^N \right) =
( \{ -\wheta_j \}_0^N \cup \{ -\cheta_j \}_1^N) \cap \rho_1^+ .
\label{e R-u-w+set}
\end{gather}
If $\la_0$ is a generalized zero of $\Mc_+ (\la) - \Mc_- (\la)$
and $ \la_0 \in
( \{ -\wheta_j \}_0^N \cup \{ -\cheta_j \}_1^N) \cap \rho_1^+ $,
the definition of the set
$\{ \wt {\zeta_j^-} \}_1^{\wt {N_\zeta^-}}$ imply that
$ \la_0 \in \{ \wt {\zeta_j^-} \}_1^{\wt {N_\zeta^-}}$ and
the generalized order of $\la_0$ equals $1/2$.
Thus, by \eqref{e M2} and the definitions of $w_+$, $U_-$, we have
\begin{gather*}
U_- ^2 (t )w_+^{-1} (t) = O (1)
\quad (t \ra t_0 ), \quad
t_0 \in ( \{ -\wheta_j \}_0^N
\cup \{ -\cheta_j \}_1^N) \cap \rho_1^+ \ .
\end{gather*}
Taking into account \eqref{e u-w+1} and \eqref{e R-u-w+set},
we get
\[
U_-^2 (t) w_+^{-1} (t) \in L^\infty (\R) .
\]
One can prove
$
U_\pm^2 (t) w_-^{-1} (t) \in L^\infty (\R)
$
in the same way.
Thus \eqref{e Cond w3} is proved.
\end{proof}

Let $\wt \theta_j^\pm $ be a point of the set
$\{ \wt \theta_k^\pm \}_1^{\wt N_\theta^\pm}$.
Let $D_j^\pm$ be a sufficiently small neighborhood of
$\wt \theta_j^\pm $ such that
\[
D_j^\pm \cap \left(
\{ \theta_k^\pm \}_1^{N_\theta^\pm} \cup
\{ \pm \wheta_k \}_0^N \cup \{ \pm \cheta_k  \}_1^N \right)
= \wt \theta_j^\pm.
\]

Put
\begin{gather}
U_\theta^\pm := \frac {\sqrt{u_\theta^\pm (\la)}}{M_2 (\la)} , \qquad \text{where}\qquad
u_\theta^\pm (\la) := \frac { \sqrt{R(\pm\la)} } {S(\pm\la)} \;
\frac{ \prodl_{j=1}^{\wt N_\zeta^\pm} (\la -\wt \zeta_j^\mp) }
     { \prodl_{j=1}^{\wt N_\zeta^\pm} (\la -(\wt \zeta_j^\mp-i\ep_1)) } .
\label{e u theta}
\end{gather}
We define $U_j^\pm$ as
$ U_j^\pm := U_\theta^\pm $ for all
$j = 1, \dots , \wt N_\theta^\pm $.

Lemma \ref{l N+} imply that $ U_\theta^\pm \in \Np (\C_+)$.

\begin{lemma}
Let Statement (iii) be true.
Then conditions \eqref{e C theta u}, \eqref{e C theta w}, and
\eqref{e C theta M}
are fulfilled.
That is, for every
$\wt \theta_j^\pm \in \{\wt \theta_k^\pm \}_1^{\wt N_\theta^\pm} $,
the following conditions hold:
\begin{gather*}
\frac{1}{|M_+ (\la) - M_- (\la)|^2}
\im \frac{1}{\wt \theta_j^\pm - \la}
\leq C_\theta^u
| U_\theta^\pm (\la)|^2  \quad \text{for} \ \la \in
D_j^\pm \cap \C_+ \ ,
\end{gather*}
\begin{gather*} 
\frac { | U_\theta^\pm (t) |^2}{ w_+ (t)} \in L^\infty (\R), \quad
\frac { | U_\theta^\pm (t) |^2}{ w_- (t)} \in L^\infty (\R) \ ,
\end{gather*}
\begin{gather*} 
\frac {1}{|\wt \theta_j^\pm - \la| \, |M_+ (\la) - M_- (\la)|}
\leq
C_\theta^M  \qquad \text{for} \  \la \in D_j^\pm \cap \C_+ ,
\end{gather*}
where $C_\theta^u$ and $C_\theta^M$ are constants.
\end{lemma}

\begin{proof}
Note that
\[
M_2 (\la) \asymp M_+ (\la) - M_- (\la) \quad (\la \ra \wt \theta_j^\pm).
\]
Therefore \eqref{e C theta u} is equivalent to
\begin{gather}
\im \frac{1}{\wt \theta_j^\pm - \la}
\leq C_1 | u_\theta^\pm (\la)| \quad \text{for} \ \la \in
D_j^\pm \cap \C_+ \ .
\label{e im theta<u}
\end{gather}
By \eqref{e zeta subset} and \eqref{e wt zeta}, it follows that
$\wt \theta_j^\pm \not \in
\{ \wt \zeta_k^\mp \}_1^{\wt N_\zeta^\pm}$.
Taking into account \eqref{e u theta} and \eqref{e SpDisc FZA2},
we see that $ u_\theta^\pm (\la) $ has a pole of the first order at
$\wt \theta_j^\pm$. This implies \eqref{e im theta<u}. Thus
\eqref{e C theta u} is proved.

Lemma \eqref{l Cond w3} and the definitions of $u^\pm$ and $u_\theta^\pm$
imply that
\begin{gather} \label{}
\frac { | U_\theta^+ (t) |^2}{ | w_\pm (t) |}  \leq C_2 \quad
\text{for} \ t \in \R \setminus \bigcup_{k=1}^{\wt N_\theta^+}
D_k^+ \ .
\end{gather}
Hence, to check condition \eqref{e C theta w} for $U_\theta^+$,
it is suffices to show that
\begin{gather} \label{e U in Dk}
\frac { | U_\theta^+ (t) |^2}{ | w_\pm (t) |}  \leq C_2 \quad
\text{for} \ t \in D_k^+ , \quad k=1, \dots , \wt N_\theta^+ \ .
\end{gather}
It is easy to see that
\begin{gather}
M_2 (\la ) \asymp (\la - \wt \theta_k^\pm)^{-1} , \quad
u_\theta^\pm (\la ) \asymp (\la - \wt \theta_k^\pm)^{-1}
\quad ( t \ra \wt \theta_k^\pm ) \label{e M Dk}
\\
\frac 1{w_\pm (t)} = O(1) \, (t-\wt \theta_k^+)^{1/2}
\quad ( t \ra \wt \theta_k^+ ) \ . \notag
\end{gather}
Combining these formulae, 
we obtain \eqref{e U in Dk}. Thus \eqref{e C theta w} for
$U_\theta^+$ is proved.
The proof of \eqref{e C theta w} for $U_\theta^-$ is similar.

Condition \eqref{e C theta M} follows from \eqref{e M Dk}.
\end{proof}

Since all the conditions of Theorem \ref{t uwSufCond} are fulfilled,
we see that $\Aess$ is similar to a selfadjoint operator.
Theorem \ref{t uwSufCond} is proved.

\subsection{Examples}
\label{ch SL ss Ex}

Let $L = - d^2/dx^2 + q(x)$ be a Sturm-Liouville operator
with a finite-zone potential $q$.
Put
\[
A:=JL=(\sgn x)( - d^2/dx^2 + q(x)) \ .
\]

\begin{definition}
We shall say that a point $a \in \si_{ess} (A) \cup {\infty}$
is \emph{a strong spectral singularity of $\Aess$} if
at least one of the following two functions
\[
 \dfrac{\Si_{ac+} ^\prime (t) } {M_+ (t) - M_- (t)} , \quad
 \dfrac{\Si_{ac-} ^\prime (t) } {M_+ (t) - M_- (t)}
\]
is not essentially bounded in any neighborhood of $a$.
\end{definition}

By Theorem \ref{t FZ Ae},
$\Aess$ is similar to a selfadjoint operator if and only if
$\Aess$ has no strong spectral singularities.
Combining Theorems \ref{t FZ Ae} and \ref{t FZ Ad}, we see that
$A$ is similar to a normal operator if and only if
the following two conditions hold:

1) $\Aess$ is similar to a selfadjoint operator;

2) all eigenvalues of $\Adisc$ are simple.

By $L (\xi,q)$ we denote the Sturm-Liouville operator
with a finite-zone potential $q(x) + \xi $,
\[
L (\xi,q) := -d^2/dx^2 +q(x) + \xi,
\]
where $\xi $ is a real constant.
Put
\[
A (\xi,q) := JL (\xi , q) = (\sgn x)( - d^2/dx^2 + q(x)+\xi)\ .
\]
Let $\Aess(\xi,q_1)$ be the part of $A(\xi,q_1)$ on $\Hsp_e$.

\begin{example} \label{ex1}
Consider the following periodic one-zone potential
\begin{equation} \label{e q1}
q_1(x) = (1-k^2 ) (2\sn ^2 (x,k^\prime ) -1 ), \ \ k\in (0,1),
\ \ k^\prime = \sqrt{1-k^2} ,
\end{equation}
where $ \sn (x,k^\prime)$ is the Jacobi elliptic function.
Then $L(\xi,q_1)$ is a one-zone periodic operator;
$L(\xi,q_1)$ has the gaps
$(-\infty , \xi)$ \ and \ $(k^2+\xi, 1+ \xi).$

The corresponding Weyl functions $M_{\pm}(\lambda)$
has the forms
\[
M_+ (\lambda) = - M_- ( - \lambda) =   i
\frac{\lambda - (\xi+1)}{ \sqrt{ (\lambda - \xi)(\lambda - (\xi+k^2)) }}
\ , \quad  0<k^2<1,
\]
(see \cite[Appendix II]{AG78}).
Theorem \ref{t FZ Ae} imply that
$\Aess(\xi,q_1)$ is similar to a selfadjoint operator
if and only if
\[
\xi\in [-1, -k^2]\cup [0,\infty).
\]
Note that for $\xi\in [-1, -k^2]$ the operator
$L(\xi, q_1)$ is not nonnegative.
If
\[
\xi\in ( -1+\sqrt{1-k^2}, -1-\sqrt{1-k^2} ) ,
\]
then $A(\xi, q_1)$ has exactly two eigenvalues
\[
 \pm\sqrt{(\xi+1)^2 - (1-k^2)} \ ;
\]
these eigenvalues are simple and nonreal.
For sufficiently small $\xi\ge 0$, the potential
$q_1 (x) +\xi $ is not nonnegative, although
$L (\xi, q_1 )\ge 0$.

Spectral properties of $A(\xi, q_1)$ are given in more detail
in the following table.
The abbreviations 'S-A' ('Norm') in the column 'Similarity'
means that
$A (\xi, q_1)$ is similar to a selfadjoint (normal)
operator.
'NonSim' in the column 'Similarity' means that $A (\xi, q_1)$ is not
similar to a normal operator.
We put
$\lambda_{\pm}(\xi):=\pm\sqrt{(\xi+1)^2 - (1-k^2)}$.

\begin{center}
Spectral properties of the operator $A(\xi ,q_1)$ \smallskip

\setlength{\extrarowheight}{4pt}
\begin{tabular}
{|c|c|c|c|}
\hline
Intervals   &  Strong       & Eigenvalues    & Similarity \\
            & spectral      &                &       \\
            & singularities &                &       \\  \hline
$ \xi \in [0,+\infty) $         &  No   &  $\lambda_{\pm}(\xi)$ &  S-A \\ \hline
$ \xi \in (-\frac{k^2}{2} , 0)$ &  $0$   & $\lambda_{\pm}(\xi)$  &  NonSim  \\ \hline
$ \xi = -\frac{k^2}{2}$         & $0$    &  No                  & NonSim  \\ \hline
$ \xi \in (-1+\sqrt{1-k^2}, -\frac{k^2}{2} )$ &
                                  $0$, $\lambda_{\pm}(\xi)$
                                         &  No                  & NonSim  \\ \hline
$\xi = -1+\sqrt{1-k^2} $ & $0$             &  No                  & NonSim  \\ \hline
$\xi \in (-k^2 , -1+\sqrt{1-k^2} )$ &  $0$    & $\lambda_{\pm}(\xi)$  & NonSim  \\ \hline
$ \xi \in [-1, -k^2 ]$       &  No            & $\lambda_{\pm}(\xi)$  & Norm \\ \hline
$ \xi \in (-1-\sqrt{1-k^2},-1 )$ & $0$         & $\lambda_{\pm}(\xi)$  & NonSim  \\ \hline
$ \xi \in -1-\sqrt{1-k^2} $ & $0$             & No                   & NonSim  \\ \hline
$ \xi \in (-\infty, -1-\sqrt{1-k^2} )$ &
                      $0$,  $\lambda_{\pm}(\xi)$ & No        & NonSim   \\ \hline
\end{tabular}
\smallskip

Table \ref{ch SL ss FZ}.1
\end{center}
\end{example}

\begin{remark} \label{r M+M}
Example \ref{ex1} shows that condition \eqref{sufficientcond} is not necessary for similarity of $A$ to a self-adjoint operator. Let us explain this.

Let $\xi >0$. Then 
$A(\xi, q_1)$ is similar to selfadjoint operator, but the function 
$\displaystyle \frac{M_+(\lambda) + M_- (\lambda)}{M_+(\lambda) - M_- (\lambda)}$,
$\la \in \C_+$, 
is unbounded in neighborhoods of the eigenvalues $\lambda_{\pm}:=\pm\sqrt{(\xi+1)^2 - (1-k^2)}$. 
Indeed, the functions $M_\pm$ are holomorphic in points $\lambda_\pm $ and $\lambda_{\pm}$ are zeroes of $M_+ (\cdot) - M_- (\cdot)$.
On the other hand, it is easy to check that 
\begin{gather*}
M_+ (\lambda_+ ) <0 , \quad M_- (\lambda_+ ) <0, \quad M_+ (\lambda_-) >0 , \quad M_- (\lambda_- ) >0.
\end{gather*}
Therefore, $M_+ (\lambda_\pm) + M_- (  \lambda_\pm) \neq 0$.
\end{remark}

\begin{example} \label{ex2}
Consider even periodic potential
\begin{equation*} 
q_2 =- {2 k^2}\bigl(1-(1-k^2) \sn^2(x,k^\prime) \bigr)^{-1} + 1 + k^2 ,
\qquad k\in (0,1) , \quad k^\prime = \sqrt{1-k^2} .
\end{equation*}

The operator $L (\xi, q_2)$ is a one-zone operator with gaps
$(-\infty , \xi)$ and $(k^2+\xi , 1+ \xi)$.
The corresponding Weyl functions $M_{\pm}(\lambda)$ have the forms
(see \cite[Appendix II]{AG78})
\[
M_+ (\lambda) = - M_- (\lambda) =  i \frac{\lambda - (\xi+k^2)}
{\sqrt{ (\lambda - \xi) (\lambda - (\xi+1)) }} \ ,
\quad 0<k^2<1.
\]
The operator $A(\xi, q_2)$ has no eigenvalues for all
$\xi\in \R$. Hence, $\Aess (\xi, q_2) = A(\xi, q_2)$.

Let $0<k^2 \leq \frac12$.
Using Theorem \ref{t FZ Ae}, we get the following result:
The operator $A(\xi, q_2)$ is
similar to a selfadjoint operator if and only if
$\xi\in [-\frac{1}{2}, -k^2]\cup [0,\infty).$
The following table describes spectral properties of
$A (\xi, q_2)$.
\begin{center}
Spectral properties of $A(\xi ,q_2)$, the case $k^2 \in (\, 0, 1/2 \, ]$
\smallskip

\setlength{\extrarowheight}{4pt}
\begin{tabular}
{|c|c|c|}
\hline
Intervals  $\xi $     & Strong spectral singularities  & Similarity \\ \hline
$ \xi \in [ 0 , +\infty ) $     &         No              & S-A      \\ \hline
$ \xi \in ( -k^2 , 0 ) $        &         $0$              & NonSim      \\ \hline
$ \xi \in [ -\frac12 , -k^2 ] $ &         No              &  S-A     \\ \hline
$ \xi \in [ -1 , -\frac12 ) $   & $\pm\sqrt{(\xi+k^2)^2 + k^2 (1-k^2)}$
                                                           &  NonSim    \\ \hline
$ \xi \in ( -\infty , -1 ) $    & $0$, $\pm\sqrt{(\xi+k^2)^2 + k^2 (1-k^2)}$
                                                           &  NonSim     \\ \hline
\end{tabular}
\smallskip

Table \ref{ch SL ss FZ}.2
\end{center}

Assume $k^2 > \frac12$. Then $A (\xi, q_2)$ is similar
to a selfadjoint operator if and only if
$\xi \geq 0$. That is $A (\xi, q_2)$ is similar to a selfadjoint operator
iff $L (\xi, q_2) \geq 0$.
The following table gives a description of spectral properties of
$A (\xi, q_2)$ in this case.

\begin{center}
Spectral properties of $A(\xi ,q_2)$, the case $k^2 \in (\, 1/2 , 1 \, )$
\smallskip

\setlength{\extrarowheight}{4pt}
\begin{tabular}
{|c|c|c|}
\hline
Intervals            & Strong spectral singularities &  Similarity \\ \hline
$ \xi \in [ 0 , +\infty ) $ &          No               &  S-A    \\  \hline
$ \xi \in [ -\frac12 , 0 )$ &          $0$               &  NonSim   \\  \hline
$ \xi \in (  -k^2 , -\frac12 ) $
                            &  $0$, $\pm\sqrt{(\xi+k^2)^2 + k^2 (1-k^2)}$
                                                         &  NonSim   \\ \hline
$ \xi \in [ -1 , -k^2] $    &  $\pm\sqrt{(\xi+k^2)^2 + k^2 (1-k^2)}$
                                                         &  NonSim \\ \hline
$ \xi \in ( -\infty , -1 ) $ & $0$, $\pm\sqrt{(\xi+k^2)^2 + k^2 (1-k^2)}$
                                                         &  NonSim \\ \hline
\end{tabular}
\smallskip

Table \ref{ch SL ss FZ}.3
\end{center}
\end{example}

\begin{example} \label{ex 3}
Let $q_3$ be potential \eqref{e q1} with $k^2 = 1/2$.
Let $\xi \in [-1,-1/2) $.
Then, combining Example~\ref{ex1} with Theorem~\ref{t CrDef},
we see that $A(\xi,q_3)$ is not definitizable,
although $\Aess (\xi,q_3)$ is similar to a selfadjoint operator and
$A (\xi, q_3)$ is similar to a normal operator.
The nonreal spectrum of $A (\xi, q_3)$ consists of
two simple eigenvalues $\lambda_{\pm}(\xi):=\pm\sqrt{(\xi+1)^2 - (1-k^2)}$. The operator
$A (\xi, q_3)$ has no real eigenvalues.
\end{example}

\end{document}